\documentclass{article}
\usepackage[utf8]{inputenc}
\usepackage{graphicx}
\usepackage{multirow}
\usepackage{amsmath,amssymb,amsfonts}
\usepackage{mathrsfs}
\usepackage{xcolor}
\usepackage{booktabs}%
\usepackage{algorithm}
\usepackage{algorithmicx}
\usepackage{algpseudocode}

\usepackage[colorlinks=true, allcolors=blue]{hyperref}
\usepackage[noabbrev,capitalize]{cleveref}
\usepackage{siunitx}
\usepackage{xcolor}
\usepackage[margin=20pt]{subcaption}
\usepackage{babel}[ngerman]
\usepackage{longtable}
\usepackage{tikz}
\usetikzlibrary{decorations.pathreplacing,calc,fit,tikzmark,backgrounds} % Load the decoration library

\usepackage{float}
\usepackage{soul}
\usepackage{caption}
\usepackage{authblk}

\usepackage[sort&compress,square,comma,numbers]{natbib}
\usepackage{array,arydshln}
\usepackage[letterpaper,top=2cm,bottom=2cm,left=3cm,right=3cm,marginparwidth=1.75cm]{geometry}
\usepackage{bbm}

\newcommand{\mpmin}[1]{\underset{#1}{\text{min}}\,}
\newcommand{\mpst}{\text{s.t.}\,}
\newcommand{\mpeq}{{}=}
\newcommand{\norm}[1]{\left\lVert#1\right\rVert}

\usepackage{xcolor}
\definecolor{light-gray}{gray}{0.95}
\newcommand{\code}[1]{\colorbox{light-gray}{\texttt{#1}}}

\makeatletter
\newcommand{\customlabel}[2]{%
   \protected@write \@auxout {}{\string \newlabel {#1}{{#2}{\thepage}{#2}{#1}{}} }%
   \hypertarget{#1}{#2}
}
\makeatother

\crefname{appendix}{Supporting Information}{Supporting Information}
\Crefname{appendix}{Supporting Information}{Supporting Information}

\title{Deterministic Global Optimization of the Acquisition Function in Bayesian Optimization: To Do or Not To Do?}

\author[1,$\dagger$]{Anastasia Georgiou}
\author[2,$\dagger$]{Daniel Jungen}
\author[2]{Luise Kaven}
\author[2]{Verena Hunstig}
\author[3,4]{Constantine Frangakis}
\author[1,5]{Ioannis Kevrekidis}
\author[2,6,7,$*$]{Alexander Mitsos}

\affil[1]{Chemical \& Biomolecular Engineering, Johns Hopkins University, Baltimore, MD 21218, USA}
\affil[2]{Process Systems Engineering (AVT.SVT), RWTH Aachen University, Aachen, Germany}
\affil[3]{Department of Biostatistics, Bloomberg School of Public Health, Johns Hopkins University, Baltimore, MD 21218, USA}
\affil[4]{Department of Medicine, Johns Hopkins University, Baltimore, MD 21218, USA}
\affil[5]{Applied Mathematics \& Statistics, Johns Hopkins University, Baltimore, MD 21218, USA}
\affil[6]{JARA-CSD, 52056 Aachen, Germany}
\affil[7]{Institute of Climate and Energy Systems, Energy Systems Engineering (ICE-1), Forschungszentrum J\"{u}lich GmbH, 52425 J\"{u}lich, Germany}
\affil[$\dagger$]{These authors contributed equally to this work.}
\affil[$*$]{Corresponding author: A. Mitsos, E-mail: \href{mailto:amitsos@alum.mit.edu}{amitsos@alum.mit.edu}}

\date{\today}

\begin{document}%
\maketitle%
%
%%%%%%%%%%%%%%%%%%%%%%%%%%%%%%%%%%%%%%%%%%%%%%%%%%%%%%%%%%%%%%%
% abstract
%%%%%%%%%%%%%%%%%%%%%%%%%%%%%%%%%%%%%%%%%%%%%%%%%%%%%%%%%%%%%%%
\begin{abstract}
    Bayesian Optimization (BO) with Gaussian Processes relies on optimizing an acquisition function to determine sampling.
    We investigate the advantages and disadvantages of using a \textit{deterministic global} solver (MAiNGO) compared to conventional \textit{local} and \textit{stochastic global} solvers (L-BFGS-B and multi-start, respectively) for the optimization of the acquisition function.
    For CPU efficiency, we set a time limit for MAiNGO, taking the best point as optimal.
    We perform repeated numerical experiments, initially using the M\"{u}ller-Brown potential as a benchmark function, utilizing the lower confidence bound acquisition function; we further validate our findings with six alternative benchmark functions.
    Statistical analysis reveals that when the acquisition function is more \textit{exploitative} (as opposed to \textit{exploratory}), BO with MAiNGO converges in fewer iterations than with the local solvers.
    However, when the dataset lacks diversity, or when the acquisition function is overly exploitative, BO with MAiNGO, compared to the local solvers, is more likely to converge to a \textit{local} rather than a \textit{global}ly near-optimal solution of the black-box function.
    L-BFGS-B and multi-start mitigate this risk in BO by introducing stochasticity in the selection of the next sampling point, which enhances the exploration of uncharted regions in the search space and reduces dependence on acquisition function hyperparameters.
    Ultimately, suboptimal optimization of poorly chosen acquisition functions may be preferable to their optimal solution.
    When the acquisition function is more exploratory, BO with MAiNGO, multi-start, and L-BFGS-B achieve comparable probabilities of convergence to a globally near-optimal solution (although BO with MAiNGO may require more iterations to converge under these conditions).
\end{abstract}

%%%%%%%%%%%%%%%%%%%%%%%%%%%%%%%%%%%%%%%%%%%%%%%%%%%%%%%%%%%%%%%
% main
%%%%%%%%%%%%%%%%%%%%%%%%%%%%%%%%%%%%%%%%%%%%%%%%%%%%%%%%%%%%%%%
\section{Introduction}
\label{sec:introduction}
Bayesian Optimization (BO), grounded in early work on Bayesian global search and related global-optimization methods~\citep{Kushner_1964,Mockus_1978,Jones_1998,Mockus_1989,Strongin_2000,Zhigljavsky_2008}, has emerged as a useful framework to optimize expensive-to-evaluate black-box functions~\citep{Shahriari_2016,Archetti_2019,Wang_2023}.
These early studies established the theoretical connection between Bayesian inference, Kriging, and radial basis function models, laying a rigorous foundation for modern BO.
Originally formulated as a deterministic strategy for sequential decision-making under uncertainty, BO has since evolved into a flexible framework for adaptive, data-driven optimization.
The approach employs the construction of a cheap-to-evaluate surrogate model, which also provides uncertainty estimates, starting with an initial (typically small) dataset.
Often a Gaussian Process (GP)~\citep{Rasmussen_2006} is used as the surrogate model.
The initial surrogate model determines the next best sampling point of the black-box function via optimization of a selected acquisition function\footnote{Referred to in early work as a ``criterion,'' including the expected improvement criterion introduced in~\citet{Mockus_1978} and anticipated in the global search procedure of~\citet{Kushner_1964}.}, typically balancing \textit{exploitation} and \textit{exploration}.
This optimization step constitutes the \textit{inner} problem, and the algorithm used to solve it is referred to as the \textit{inner} solver.
The surrogate model is then improved upon with the newly evaluated point; repetition of this process (constituting the \textit{outer} loop) locates a (hopefully \textit{global}ly near-optimal) minimum of the black-box function.
We refer to the process of locating a minimum of the black-box function as the \textit{outer} problem.

A core challenge in BO is that acquisition functions are inherently prone to \textit{local minima}: if the original objective function is multimodal, this characteristic is reflected in the exploitative part of the acquisition function after sufficient samples; moreover, irrespective of the original objective function, if the acquisition function is exploratory, uncertain regions result in additional local/global minima.
Historically, BO was motivated by deterministic global optimization ideas, and several works have stated that global optimization of the acquisition function (thereby avoiding local minima) is important for achieving strong BO performance~\citep{Wilson_2018b,Kim_2021}.

Despite this, practical BO implementations rely on a wide variety of inner solvers~\textemdash{}  ranging from local, gradient-based methods to stochastic global heuristics and deterministic global solvers; c.f., \cref{tab:bo_methods} in \cref{sec:bo_frameworks}.
Among the most commonly used inner solvers are \textit{local, gradient-based solvers} such as the L-BFGS-B algorithm~\citep{Zhu_1997}, which are popular for their speed and ease of implementation.
However, local methods such as L-BFGS-B or trust-region algorithms~\citep{Cartis_2021,Wu_2023} do not guarantee convergence to a global minimum of the acquisition function.
To mitigate this, several alternatives\footnote{Additionally, structural or decomposition-based methods such as additive-structure models~\citep{Oh_2018,Wang_2018} and hybrid or partition-based strategies such as the trust-region-inspired TuRBO framework~\citep{Eriksson_2019} and hierarchical domain partitioning~\citep{Kawaguchi_2015} have also been proposed for use in BO.
Although we do not explicitly consider these approaches here, we note that they likewise do not provide deterministic global guarantees or only ensure convergence in infinite time.} have been proposed, including \textit{stochastic global} approaches such as multi-start strategies or evolutionary algorithms~\citep{Hansen_2006,Bradford_2018,song_2024}.
Nevertheless, these two categories~\textemdash{}  local and stochastic global~\textemdash{}  either lack global guarantees or only guarantee convergence in infinite time~\citep{Rinnooy_1985}.

Using a \textit{deterministic global} solver, as e.g., classical Branch-and-Bound algorithms \citep{Land_1960} or modern global optimization tools like BARON~\citep{Khajavirad_2018}, SCIP~\citep{SureshBolusani_2024}, Gurobi~\citep{gurobi}, or MAiNGO~\citep{maingo},  directly guarantees the avoidance of local minima.\footnote{Similarly, mixed-integer programming is rarely employed in BO with GP surrogates but has found application with non-Gaussian models, such as neural networks~\citep{Papalexopoulos_2022} and decision trees~\citep{Thebelt_2022}.}
Among the deterministic global methods used in BO frameworks, DIRECT~\citep{Jones_1993} is the most widely used.
Still, its use in BO has been limited~\citep{Li_2017b}, and its adoption remains constrained as it becomes computationally prohibitive in high-dimensional settings~\citep{Kim_2021}.

Recent advances attempt to overcome these limitations by replacing the GP kernel with a piecewise-linear approximation, thereby offering, in theory, global optimality guarantees in the BO context while significantly improving computational efficiency~\citep{Xie_2025}.
In particular, \citet{Xie_2025} introduce a piecewise-linear approximation of the GP kernel, which yields an approximated acquisition function whose optimization can be formulated as a mixed-integer quadratic program (PK-MIQP).
\cite{Xie_2025} prove that under suitable assumptions, (a) their approximated kernel mean (and variance) converges in the limit to the true mean (and variance) and (b) give probabilistic bounds on regret for the case where the approximation error is finite.
\cite{Xie_2025} demonstrate the benefit of their PK-MIQP formulation on synthetic functions, constrained benchmarks, and a hyperparameter tuning task and compare the results to various other gradient- and sampling-based methods.
However, their theoretical guarantees do not hold for all kernels, as the kernel must fulfill a smoothness assumption.

This work explores the advantages and disadvantages of employing a deterministic global solver in the context of Bayesian Optimization (BO), without relying on any of the previously mentioned approximations, i.e., we investigate how the choice of the inner solver affects the solution performance of the outer problem.
To this end, we consider an initial small dataset, simulating a real scenario where experiments are costly, and use the Müller-Brown potential as a baseline benchmark problem.
For the inner solver, we choose a representative solver from each solver class: deterministic global, local, and stochastic global.
We aim to answer the following two research questions:
\begin{enumerate}
    \item \label{q:one} How does the choice of the inner solver affect successful convergence to a globally near-optimal solution of the outer problem, and to what extent is that success dependent on the initial dataset?
    \item \label{q:two} Given a sufficiently informative initial dataset (to be discussed in \cref{sec:question_one}), BO with which inner solver converges to a globally near-optimal solution in the fewest iterations?
\end{enumerate}
In essence, these questions aim to delineate the regimes in which each solver class (local, stochastic global, deterministic global) is most effective~\textemdash{} expanding on prior studies that examined a narrow set of regimes that favored global acquisition optimization~\citep{Wilson_2018b,Kim_2021,Xie_2025}~\textemdash{} and determining, across a broader range of conditions, whether deterministic global solvers consistently offer an advantage or whether local or stochastic global methods may be preferable.
We adopt the standard Lower Confidence Bound (LCB) acquisition function and vary its exploration parameter, using primarily fixed and occasionally scheduled values: the resulting insights apply specifically to this acquisition function.
Isolating the effect of the inner solver becomes increasingly important as modern hardware and reduced-complexity surrogate representations (such as ~\citep{Xie_2025}) continue to lower the practical barriers to deterministic global optimization for GP-embedded acquisition functions.

For our investigation, we utilize the L-BFGS-B solver as a representative of the class of \textit{informed} local solvers (ILS).
We adapt the same L-BFGS-B solver by introducing multiple restarts to create an \textit{informed} multi-start (IMS) solver, which is our representative for the class of stochastic global solvers.
The term ``\textit{informed}'' is added in both cases to signify that the solver's initialization is strategically selected, as detailed in~\cref{sec:ils,sec:ims}.
For the deterministic global solver, we use our open-source solver MAiNGO (McCormick-based Algorithm for mixed-integer Nonlinear Global Optimization)~\citep{maingo}, which employs MeLOn (Machine Learning models for Optimization)~\citep{MeLOn} as a submodule to integrate machine learning models (in our case, GPs) into optimization problems.
The parallel computing capabilities of MAiNGO, in conjunction with its reduced space formulation of GPs and custom relaxations for various functions, make it a promising solver, shown to outperform alternatives~\citep{Schweidtmann_2021} and, more importantly, make global optimization with GPs embedded tractable in more realistic settings.
However, computational costs still remain a major challenge.
To maintain an overall tractable computational time of our numerical experiments, we terminate each optimization run in MAiNGO, i.e., the current optimization of the acquisition function, after a pre-defined time limit, \textit{assuming} that the current solution is indeed a globally near-optimal solution of the acquisition function (corroboration, to follow, \cref{sec:maingo}).
An overview of the inner solvers, their categories, and our selection rationale is given in \cref{tab:solver-types}.

\begin{table*}[htbp]
    \centering
    \caption{
        Representative solvers used for the inner optimization of the acquisition function, grouped by solver category.
        The three solvers (MAiNGO, ILS, and IMS) illustrate deterministic global, local, and stochastic global strategies, respectively.
    }
    \label{tab:solver-types}
    \begin{tabular}{p{0.13\textwidth} p{0.10\textwidth} p{0.17\textwidth} p{0.34\textwidth} p{0.10\textwidth}}
        \textbf{Category}   & \textbf{Solver} & \textbf{Underlying Method} & \textbf{Rationale for Selection} & \textbf{Typical in BO?} \\ [0.5ex]
        \hline
        Local               & ILS   & Quasi-Newton (L-BFGS) from single informed initial point      & Serves as a lightweight, gradient-based optimizer that efficiently finds (a local) optimum & Common \\ \\
        Stochastic global   & IMS   & Quasi-Newton (L-BFGS) started from multiple points            & Representative of common practice in BO toolkits & Very common \\ \\
        Global deterministic& MAiNGO& Deterministic branch-and-bound using McCormick relaxations    & Reduced-space formulation enables efficient optimization of models with GPs embedded & Rare\\
        \hline
    \end{tabular}
\end{table*}

For the illustrative examples provided in this work, we use \code{GPYTorch} GPs (as currently required for MeLOn) with the standard \code{BOTorch} \href{https://github.com/pytorch/botorch/tree/v0.9.5}{v0.9.5}~\citep{Balandat_2020} framework, which allows us to switch between various optimizers easily.
We compare the points visited by BO when optimizing the acquisition function using ILS, IMS, and MAiNGO and aim to identify whether the number of black-box function evaluations before BO convergence is affected.

In our pursuit of addressing the two research questions \ref{q:one} and \ref{q:two}, the main contributions of this work can be summarized as follows:
\begin{itemize}
    \item   We explore the use of a deterministic global solver (MAiNGO) for optimizing the acquisition function in BO, extending the one-shot global optimization of GP in \cite{Schweidtmann_2021}.
            Given sufficient (finite) time, MAiNGO guarantees finding an (approximate) global solution of the inner problem, but for CPU efficiency, we impose a time limit, \textit{assuming} optimality from stabilized bounds; manual checks support our assumption (\cref{sec:maingo}).
            Such (strict) time limits would be too restrictive in practice.
    \item   We provide an in-depth analysis of (a) the statistical probability of convergence to a globally near-optimal solution of the black-box function and (b) the statistical distribution of iterations until BO convergence when using ILS, IMS, and MAiNGO, applied to the M\"{u}ller-Brown potential, which we use as our baseline case study.
    \item   We demonstrate how statistical methods can quantify variance in solver performance and analyze the trade-off between selecting a high-quality initial dataset and choosing the solver.
    \item   We further validate our findings on BO using ILS, IMS, and MAiNGO with six alternative benchmark functions characterized by higher problem dimensionality and/or multiple global minima.
\end{itemize}
The paper is organized as follows.
\cref{sec:problem_formulation} introduces the BO framework and the associated formal problem statements of the inner and outer problem.
\cref{sec:approach} presents our numerical study design.
\cref{sec:methods} describes the used inner solvers, BO termination criterion, and case study formulations.
\cref{sec:results} presents our findings, detailing the performance of the various inner solvers \textit{for our baseline case study, the M\"{u}ller-Brown potential, and six additional benchmark functions.}
\cref{sec:discussion} discusses these results, while \cref{sec:recommendations} offers insights and recommendations on using global deterministic solvers in BO.
Finally, we outline future research directions in \cref{sec:future_work}.

\section{Problem Formulation}
\label{sec:problem_formulation}
In the following, we provide a high-level introduction to BO with GPs and more formally introduce the notions of \textit{outer} and \textit{inner} problem.
For a thorough introduction to BO with GPs, the interested reader may refer to readily available literature, e.g.,~\citep{Mockus_1989,Rasmussen_2006,Garnett_2023,Frazier_2018} and the references therein.
\cref{tab:notation_summary} provides an overview of the used symbols.
We use italics for scalar-valued, bold italics for vector-valued symbols, and calligraphic for sets, stochastic processes, or measures.

The GP prior is described by the mean function $\mu$ (which is a function of $\boldsymbol{x}$) and its covariance function/kernel $k$, which is positive semi-definite.
As we do not consider noise in this contribution, $f$ is distributed as a GP, i.e., $f \sim \mathcal{GP}\left(\mu\left(\boldsymbol{x}\right), k\left(\boldsymbol{x},\boldsymbol{x}'\right)\right)$ with $\mu\left(\boldsymbol{x}\right) = \mathbb{E}\left[f\left(\boldsymbol{x}\right)\right]$ and $k\left(\boldsymbol{x},\boldsymbol{x}'\right) =  \mathbb{E}\left[\left(f\left(\boldsymbol{x}\right) - \mu\left(\boldsymbol{x}\right)\right) \left(f\left(\boldsymbol{x}'\right) - \mu\left(\boldsymbol{x}'\right)\right)^{T}\right]$.
We assume the prior mean function to be zero, use a Mat\'{e}rn kernel $(\nu = 5/2)$, and optimize the GP hyperparameters by maximizing the marginal log-likelihood.

The simplified pseudocode of BO with GPs to approximate the global minimizer $\hat{\boldsymbol{x}}^{*}$ of the black-box function $f$ is given in~\cref{alg:bayes_opt}.
We call the optimization problem of finding an approximate global minimizer of the black-box function the \textit{outer} problem, i.e., $\hat{\boldsymbol{x}}^{*} \in \text{arg}\, \min_{\boldsymbol{x}} f\left(\boldsymbol{x}\right)$.
Starting with the initial data $\mathcal{D} = \left( \boldsymbol{x}^{\mathrm{init}}, f\left(\boldsymbol{x}^{\mathrm{init}}\right) \right)$, first the GP is fitted to approximate this data.
In this step, the hyperparameters of the kernel are optimized to best fit the observed data $\mathcal{D}$ and the GP posterior is obtained through conditioning the prior on the data $\mathcal{D}$, i.e., $f \sim \mathcal{GP}\left(\mu\left(\left.\boldsymbol{x}\,\right|\mathcal{D}\right), k\left(\left.\boldsymbol{x},\boldsymbol{x}'\,\right| \mathcal{D}\right)\right)$.
The posterior provides both mean predictions and uncertainty estimates of the black-box function $f$.
Next, an \textit{acquisition function} is defined, which uses the posterior of the GP.
The acquisition function determines the utility of sampling each candidate location, balancing \textit{exploration} of uncertain regions and \textit{exploitation} of promising low-value regions.
We call the optimization of the acquisition function $\alpha$ the \textit{inner} problem, i.e., $\boldsymbol{x}^{t+1} \in \text{arg}\,\max_{\boldsymbol{x}} \alpha\left(\left. \boldsymbol{x} \,\right| \mathcal{D}\right)$.
The solution point $\boldsymbol{x}^{t+1}$ obtained determines the next point at which the black-box function $f$ is evaluated.
The newly acquired sample is then incorporated into the dataset, the surrogate model is updated, and the loop repeats until a convergence criterion is fulfilled.
The termination criterion used and the list of all numerical tolerances will be introduced in \cref{sec:bo_termination_critereon} and \cref{tab:testcase-setup-table}, respectively.
Upon meeting the termination criterion, the best observed point $\hat{\boldsymbol{x}}^{*}$ is returned.
Throughout, we use the same box constraints, i.e., $\boldsymbol{x} \in [\boldsymbol{x}^l, \boldsymbol{x}^u]$, to both the outer and inner optimization problems at every BO iteration, where $\boldsymbol{x}^l$ and $\boldsymbol{x}^u$ are case-study-specific lower and upper bounds on the decision variable.

Note that optimization of the acquisition function is nontrivial, since acquisition functions are typically nonconvex and multimodal.
Multimodality of the acquisition function can be caused by two factors:
(i) the original objective function exhibits multiple (local) optima and (ii) if the acquisition function is exploratory, uncertain regions result in additional local/global minima.
Further note that, even if the \textit{inner} problem is solved to global optimality in all iterations, no guarantee can be given (without additional assumptions) that the \textit{outer} problem is solved to global optimality.

\begin{algorithm}[htbp]
    \caption{
      Pseudocode of BO with GPs.
    }
    \label{alg:bayes_opt}
    \begin{tikzpicture}[overlay, remember picture]
        \node[
            fill=black!07,
            fit={($(20pt,-46pt)$) ($(340pt,-135pt)$)},
            inner sep=4pt
        ] {};
    \end{tikzpicture}
    \begin{algorithmic}[1]
        \Require Termination criterion, initial sampling points $\boldsymbol{x}^{\mathrm{init}}$, and black-box objective function $f$
        \Ensure Approximate global minimizer $\hat{\boldsymbol{x}}^\star$

        \State Initialize dataset with initial samples, i.e., $\mathcal{D} \leftarrow \left( \boldsymbol{x}^{\mathrm{init}}, f\left(\boldsymbol{x}^{\mathrm{init}}\right) \right)$
        \tikz[remember picture] \coordinate (outer-start);
        \While{termination criterion not met}
            \State Fit $\mathcal{GP}$ to data $\mathcal{D}$
            \State Define an acquisition function $\alpha\left(\left. \boldsymbol{x} \,\right| \mathcal{D}\right)$
            \State
                \tikz[baseline=(text.base), remember picture]{
                    \node[fill=gray!40, inner sep=0.75pt, outer sep=1pt, align=left] (text)
                    {Solve $\text{arg}\,\max_{\boldsymbol{x}} \alpha\left(\left. \boldsymbol{x} \,\right| \mathcal{D}\right)$ to select next candidate $\boldsymbol{x}^{t+1}$\hphantom{\hspace{8pt}}};
                    \coordinate (inner-start) at (text.north east);
                    \coordinate (inner-end)   at (text.south east);
                }
            \State Evaluate $f\left(\boldsymbol{x}^{t+1}\right)$ and augment $\mathcal{D}$ with $\left(\boldsymbol{x}^{t+1},f\left(\boldsymbol{x}^{t+1}\right)\right)$
            \State Increment counter $t$
        \EndWhile
        \tikz[remember picture] \coordinate (outer-end);
        \State \tikz[remember picture, baseline=(ret.base)]{
                    \node[inner sep=0pt, outer sep=0pt, anchor=base west] (ret)
                        {\textbf{return} $\hat{\boldsymbol{x}}^{*} \leftarrow \text{arg}\,\min_{\boldsymbol{x} \in \mathcal{D}} f\left(\boldsymbol{x}\right)$};
}
    \end{algorithmic}
    \begin{tikzpicture}[overlay, remember picture]
        \draw[decorate, decoration={brace, amplitude=5pt, mirror}]
            ([xshift=1pt]inner-end) --
            ([xshift=1pt]inner-start)
            node[midway, xshift=5pt,right]{\parbox{20pt}{Inner \\ problem}};
        % helper
        \path (current page.south west) ++(432pt,0) coordinate (braceX);
        \coordinate (outer-brace-top)    at ($(outer-start -| braceX)$);
        \coordinate (outer-brace-bottom) at ($(ret.south east -| braceX)$);
        \draw[decorate, decoration={brace, amplitude=5pt, mirror}]
          ([xshift=1pt,yshift=-4pt]outer-brace-bottom) --
          ([xshift=1pt,yshift=-2pt]outer-brace-top)
          node[midway, xshift=5pt, right]{\parbox{20pt}{Outer \\ problem}};
    \end{tikzpicture}
\end{algorithm}

\section{Approach}
\label{sec:approach}

A key distinction for BO utilizing different inner solvers lies in the reproducibility of a solution trajectory: BO with a deterministic global (inner) solver consistently produces a single, reproducible value for the number of iterations to convergence under fixed inputs ((a) initial dataset; (b) acquisition function; and (c) termination criterion; and no noise).
Further, BO with a deterministic global (inner) solver has the same probability (0 or 1) of convergence to a globally near-optimal solution of the outer problem: if the computations are repeated, a globally near-optimal solution of the outer problem is always (not) found.
In contrast, BO with ILS or with IMS exhibits variability in the number of iterations to convergence and to identify a globally near-optimal solution of the outer problem due to the stochastic initialization of the inner solvers (described in \cref{sec:ils,sec:ims}).
This results in a distribution of iteration counts for repeated runs under identical input conditions and also in a non-trivial probability~\textemdash{}  of generally neither 0 nor 1~\textemdash{}  of convergence to a globally near-optimal solution of the outer problem.
Further variability is introduced through the initial dataset.
Even if a deterministic global (inner) solver is being used (in conjunction with the same (b) acquisition function; and (c) termination criterion; and no noise), the probability of convergence and the number of iterations to converge depends on the initial dataset.

To better understand the impact of both the inner solver and the initial dataset on the performance of BO, we distinguish among \textit{case studies}, \textit{experiments}, and \textit{runs}.
We define a \textit{case study} as characterized by the same underlying black-box function, acquisition function, and termination criterion.
Within each \textit{case study}, we conduct multiple \textit{experiments}, each defined by a unique initial dataset (all datasets of the same size).
A \textit{run} refers to a single complete execution of the BO loop for a given \textit{experiment}: we conduct multiple \textit{runs} of each \textit{experiment} for each \textit{case study}.

Through our design of the empirical study, we establish the following statistical framework for a fixed \textit{case study}.
Let $y^{s}_{dat,r}$ denote the indicator variable for whether solver $s$ ultimately converges to a globally near-optimal solution of the black-box function in \textit{run} $r$ of an \textit{experiment} characterized by initial dataset $dat$.
The statistical representation of the data design is then structured as follows:
\begin{enumerate}
    \item Conditionally on a given dataset $dat$, the convergence behaviors are independent among different solvers;
    \item Conditionally on a given dataset, indicators of $y^{\text{ILS}}_{dat,r=1},..., y^{\text{ILS}}_{dat,r=31}$, are independent Bernoulli given $dat$, with probability $p^{\text{ILS}}_{dat}=\mathbb{E}_{\text{over } r}\{y^{\text{ ILS}}_{dat,r}\mid dat\}$; and similarly for using BO with IMS.
    For BO with MAiNGO, given $dat$, $y^{\text{MAiNGO}}_{dat,r=1}$ also identifies the probability $p^{\text{MAiNGO}}_{dat}$, which is, due to MAiNGO being a deterministic solver, $0$ or $1$.
    \item For solver $s$, the convergence probabilities $p^{\text{s}}_{dat}, dat=1,...$ are independent and identically distributed; their average is denoted below as $p^{\text{s}}:=\mathbb{E}_{\text{over } dat}\{y^{\text{ s}}_{dat}\}$.
\end{enumerate}

Our structure of the empirical study, as well as its associated statistical framework, directly reflects the research questions \ref{q:one} and \ref{q:two} posed in \cref{sec:introduction}.
First, we examine how inner solver choice and initial data influence \textit{where} BO converges~\textemdash{}  that is, whether the trajectory terminates near the global optimum.
Second, conditioned on successful convergence, we assess \textit{how many iterations} BO requires to reach a globally near-optimal solution under each inner solver.
These analyses are carried out across multiple experiments and then multiple runs per experiment, allowing us to quantify the role of both inner solver-induced stochasticity and dataset initialization.
We explore these questions under realistic BO conditions, using only a small number of initial data points and reasonable, but \textit{not} finely-tuned, acquisition-function hyperparameters (as is realistic).
Last but not least, we investigate whether our findings are generalizable to higher dimensional case studies.

\section{Methods}
\label{sec:methods}
In the following, we cover important details of BO with GPs for our investigation: the acquisition function, the solvers employed, and the termination criterion.

\subsection{Acquisition Functions}
\label{sec:acquisition_functions}
Common acquisition functions in BO include the following: probability of improvement, expected improvement, upper / lower confidence bound, entropy-based, Monte-Carlo-based, and feasibility-based functions~\citep{Gan_2021,Wilson_2018b}.
We limit our study to the popular~\citep{Gan_2021} lower confidence bound acquisition function, defined as
\begin{equation}
    \label{eq:lcb}
    \tag{LCB}
    LCB = \mu - \kappa \sigma.
\end{equation}
The hyperparameter $\kappa$ controls the amount of exploitation vs.~exploration: a low value exploits the model by decreasing the relative weight on the mean $\mu$, whereas a high value explores by increasing the weight on the standard deviation $\sigma$.
Prior works~\citep{Srinivas_2012, Kandasamy_2015} have established cumulative regret bounds by dynamically adjusting $\kappa$ over iterations, ensuring a balance between exploration and exploitation.
However, in this study, we mainly study a fixed $\kappa$, the default in many BO frameworks (e.g., \code{BOTorch}).
We aim to minimize \eqref{eq:lcb} with each BO iteration.

\subsection{Solvers Used for Optimizing the Inner Problem}
\label{sec:optimizers}
As we investigate the influence of the solver used for the solution of the \textit{inner} problem on the solution of the \textit{outer} problem, we selected three solvers as representatives of the different solver classes used in the literature: deterministic global, local, and stochastic global solvers, c.f., \cref{tab:solver-types}.

\subsubsection{Informed Local Solver (ILS)}
\label{sec:ils}
We utilize the standard L-BFGS-B solver, available in \code{scipy.optimize}, for our informed local solver (ILS).
The term ``informed'' is added to signify that the solver's initialization is strategically selected.
Solver initialization plays an important role for local optimization solvers.
For nonconvex problems, the initial point significantly impacts the solver's ability to find a global optimum.
Commercial local solvers and optimization experts typically implement an initialization strategy.
We opt for the standard selection process for solver initialization used in \code{BOTorch} \href{https://github.com/pytorch/botorch/tree/v0.9.5}{v0.9.5}~\citep{Balandat_2020}, c.f., \code{botorch/optim/initializers.py} on \href{https://github.com/pytorch/botorch/blob/v0.9.5/botorch/optim/initializers.py#L892}{GitHub}:
First, 20 candidates are drawn using a Sobol sequence~\citep{Sobol_1967} within the search space.
From these 20 candidates, a single one is drawn with a probability proportional to its acquisition function value, normalized by the mean and standard deviation of the acquisition values across all candidates.
Then, the selected candidate is used as the initial point for the L-BFGS-B solver.
This approach prioritizes initialization from a candidate with a better (lower) acquisition function value.

\subsubsection{Stochastic Global Solver (IMS)}
\label{sec:ims}
For stochastic global optimization, we build on ILS, as described in \cref{sec:ils}, but introduce multiple restarts to further explore the search space.
We perform five independent restarts of ILS, where, within each restart, the previously described selection process is used for the initialization of the employed L-BFGS-B solver.
This informed multi-start (IMS) procedure aims to avoid local minima, though it does not guarantee that convergence to a global minimum will be achieved in finite time.

\subsubsection{Deterministic Global Solver (MAiNGO)}
\label{sec:maingo}
We employ our open-source deterministic global solver MAiNGO \href{https://git.rwth-aachen.de/avt-svt/public/maingo/-/tree/v0.8.1/?ref_type=tags}{v0.8.1}~\citep{maingo}, which leverages reduced-space formulations of GPs through \href{https://git.rwth-aachen.de/avt-svt/public/MeLOn}{MeLOn}~\citep{MeLOn} and custom relaxations of the acquisition function (implemented in MC++~\citep{Chachuat_2015}).
Recall that this configuration offers drastic speed-ups of the optimization process compared to full-space formulations with standard deterministic global solvers~\citep{Schweidtmann_2021}, making BO in realistic settings tractable.
Note that for the investigated acquisition function, i.e., \eqref{eq:lcb}, no custom relaxation is implemented in MAiNGO.

We configure MAiNGO to use IBM CPLEX v22.1.1~\citep{cplex} as the lower bounding solver, SLSQP as the upper bounding solver, and IPOPT~\citep{IPOPT} for multi-start.
To reduce the CPU time, we terminate each optimization run in MAiNGO, i.e., the current iterate, after \SI{30}{\minute}, \textit{assuming} that the current solution is indeed a globally near-optimal solution of the acquisition function.
This assumption indicates that we have obtained a globally near-optimal solution through the upper-bounding procedure, but the lower-bounding procedure has not yet (fully) confirmed the solution.
The observation that in these cases, the upper bound remains unchanged while the lower bound continues to increase, further corroborates the assumption.
We validate and confirm our assumption by random manual inspection of the acquisition function.
If we terminate an iterate and that iterate has a remaining relative optimal gap greater than $1$ order of magnitude higher than the set relative optimality tolerance \code{epsilonR} in MAiNGO (e.g., the true gap is $\ge 0.1$ when $\code{epsilonR=0.01}$), we loosen \code{epsilonR} by one order of magnitude for subsequent iterates (e.g., now $\code{epsilonR=0.1}$).
Within each BO run, we start with $\code{epsilonR}=0.01$.
All other MAiNGO settings have been left at the default values.

\subsection{BO Termination Criterion}
\label{sec:bo_termination_critereon}
An important consideration, but still an open question, is when to appropriately terminate the loop of the outer problem.
The most common approach is a resource-based termination, where a fixed budget~\textemdash{}  typically defined by the total number of iterations or the total cost of experiments~\textemdash{}  is prespecified.
While this approach does not guarantee convergence to a minimum, it is straightforward to implement and practical in many applications.

Another category of termination criteria evaluates the progress of the BO.
A commonly used metric, especially used in benchmarking, is a simple regret, defined as $\norm{f^{*} - f^{t}} < \varepsilon$, with $f^{*}$ being the globally optimal objective value of the black-box function $f$, $f^{t}$ being the objective value for the current iterate $\boldsymbol{x}^{t}$, and $\varepsilon$ representing the absolute tolerance.
However, this method requires prior knowledge of the globally optimal objective value: we refrain from using such knowledge-dependent criteria to ensure the study remains both practical and realistic (though we provide such simple regret plots in \cref{sec:regret-plots}).\footnote{A more general alternative to simple regret is the expected minimum simple regret~\citep{Ishibashi_2023}, which removes the need for prior knowledge of the globally optimal objective value.
Other progress-based approaches include monitoring local regret within a convex region~\citep{Li_2023}, detecting stagnation in the best-observed minimum over a specified number of iterations, or terminating when the expected improvement or probability of improvement falls below a threshold~\citep{Nguyen_2017}.
Probabilistic, confidence-based measures have also been proposed in the literature~\citep{Dai_2019b,Makarova_2022,Wilson_2024,Bolton_2004}.}

Herein, we adapt a simple and practical progress-based termination criterion.
Specifically, BO is terminated when either the distance between the current candidate $\boldsymbol{x}^{t}$ and any prior candidate $\boldsymbol{x}^{k}$, $k=\left\{1,...,t-1\right\}$, falls below the predefined threshold $\varepsilon_{x,1}$; or again when $\boldsymbol{x}^{t}$ and any prior candidate $\boldsymbol{x}^{k}$, $k=\left\{1,...,t-1\right\}$, falls below the predefined threshold $\varepsilon_{x,2}$ in conjunction with the objective function value of the current candidate $f^{t}$ being close to the best-found solution $f^{*,t-1} = \min_{k=\left\{1,..,t-1\right\}} f^{k}$:
\begin{equation}
    \label{eq:termination_critereon}
    \tag{TC}
    \begin{aligned}
        & \underbrace{
            \left( \min_{k = \left\{1,..,t-1\right\}} \norm{\boldsymbol{x}^{t} - \boldsymbol{x}^{k}} < \varepsilon_{x,1} \right)
          }_{\customlabel{eq:TC-1}{\text{TC-1}}} \\
        & \quad \lor
        \left(
            \underbrace{
                \min_{k = \left\{1,..,t-1\right\}} \norm{\boldsymbol{x}^{t} - \boldsymbol{x}^{k}} < \varepsilon_{x,2}
            }_{\customlabel{eq:TC-2.1}{\text{TC-2.1}}}
            \land
            \underbrace{
                \left(
                    \norm{ f^{t} - f^{*,t-1} } < \varepsilon_{f,r} \cdot f^{*,t-1}
                    \lor
                    \norm{ f^{t} - f^{*,t-1} } < \varepsilon_{f,a}
                 \right)
             }_{\customlabel{eq:TC-2.2}{\text{TC-2.2}}}
        \right),
    \end{aligned}
\end{equation}
with $\varepsilon_{x,1} < \varepsilon_{x,2}$ and $\varepsilon_{f,r}, \ \varepsilon_{f,a} > 0$.
The thresholds $\varepsilon_{f,r}$ and $\varepsilon_{f,a}$ represent the relative and absolute tolerance in objective function value, respectively.
The key idea behind this easy-to-implement termination criterion is that for systems without any (measurement) noise, as considered in our investigation, we would not repeat the same experiment multiple times (the candidate is in close proximity to another one; \ref{eq:TC-1} in \eqref{eq:termination_critereon}).
Furthermore, we would terminate once the candidate is in proximity to another one, and the current iterate is comparable (in terms of objective value) to the best-found objective function value (\ref{eq:TC-2.1} and \ref{eq:TC-2.2} in \eqref{eq:termination_critereon}, respectively).
A complete table of all numerical tolerances used across case studies is provided in \cref{tab:testcase-setup-table}.
Future work will explore the effects of alternative stopping criteria.

\subsection{Case Studies}
\label{sec:case_studies}
In the following, we briefly review the benchmark functions that we use in our investigation.
We choose the M\"{u}ller-Brown potential function (\cref{sec:muller_brown_potential}) as our baseline benchmark problem to evaluate the performance of ILS, IMS, and MAiNGO within BO.
In \cref{sec:additional_test_cases}, we introduce six additional
benchmark functions characterized by higher dimensions or multiple global minima, which we use to validate our findings on BO using ILS, IMS, and MAiNGO.
A summary of all benchmark functions and their key characteristics is provided in \cref{tab:benchmarks}.

\begin{table*}[htbp]
    \centering
    \caption{
        Benchmark problems used to evaluate BO with various solvers.
    }
    \label{tab:benchmarks}
    \begin{tabular}{p{0.14\textwidth} p{0.05\textwidth} p{0.28\textwidth} p{0.40\textwidth}}
        \hline
        \textbf{Benchmark} & \textbf{Dim.} & \textbf{Description / Source} & \textbf{Rationale / Special Features} \\ [0.5ex]
        \hline
        Camelback & 2 & Test function with two global minima and one local maximum; standard in optimization literature~\citep{Surjanovic_2013}. & Serves as an easily visualized 2D case for functions with multiple global minima. \\ \\
        GKLS-2D & 2 & GKLS generator~\citep{Gaviano_2003} instance producing a quadratic base function perturbed by random local minima. & Provides a controlled multimodal landscape with a known global optimum; used for tunable difficulty and ability to systematically increase dimensions. 10 local minima exist. \\ \\
        GKLS-3D & 3 & GKLS generator~\citep{Gaviano_2003} instance with moderate nonconvexity in three variables. & Represents increased difficulty over GKLS-2D. 10 local minima exist.\\ \\
        Ackley & 3 & Smooth, multimodal benchmark with an exponentially damped cosine structure~\citep{Surjanovic_2013}. & Tests solver behavior on continuous, oscillatory landscapes with many shallow minima and one deep global minimum. \\ \\
        GKLS-4D & 4 & GKLS generator~\citep{Gaviano_2003} instance with increased dimension.  & Represents further increased difficulty for GKLS; 10 local minima exist. \\ \\
    Hartmann-4D & 4 & Standard function provided in~\citep{Surjanovic_2013}. & Provides a smooth landscape with multiple local minima and a single narrow global basin in 4D; often used in BO benchmarks. \\
        \hline
    \end{tabular}
\end{table*}

\subsubsection{M\"{u}ller-Brown Potential}
\label{sec:muller_brown_potential}
The canonical M\"{u}ller-Brown potential~\citep{Muller_1979} is defined as
\begin{equation}
    \label{eq:muller-brown}
    U(\boldsymbol{x}) = \sum_{i = 1}^4 A_i \, \exp \left(
           a_i \left(x_{1} - w^1_{0i}\right)^{2 }
         + b_i \left(x_{1} - w^1_{0i}\right) \left(x_{2} - w^2_{0i}\right)
         + c_i \left(x_{2} - w^2_{0i}\right)^{2} \right),
\end{equation}
with input $\boldsymbol{x}\in \mathbb{R}^2$ and the coefficients as defined in \cref{tab:mb}.
In our considered domain of $x_1\in [-1.5, 1], x_2\in [-0.5, 2]$, $U$ has two local minima and one global minimum.
The surrogate model used in BO is a GP with Mat\'{e}rn kernel $(\nu = 5/2)$ and marginal log-likelihood.
The lengthscale hyperparameters are locally optimized with each BO iteration, and we assume the data points have no associated noise.
Input features are scaled to be between 0 and 1, and the output is scaled to have (approximately) zero mean and unit variance.
We optimize \eqref{eq:lcb} with each BO iteration.

\begin{table*}[htbp]
    \centering
    \caption{
        Coefficients of the M\"{u}ller-Brown potential.
    }
    \label{tab:mb}
    \begin{tabular}{c c c c c c c}
         $i$ & $A_{i}$ & $a_{i}$ & $b_{i}$ & $c_{i}$ & $w^{1}_{0i}$ & $w^{2}_{0i}$\\ [0.5ex]
         \hline
         1 & $-200$ & $-1$      & 0     & $-10$     & 1         & 0 \\
         2 & $-100$ & $-1$      & 0     & $-10$     & 0         & 0.5 \\
         3 & $-170$ & $-6.5$    & 11    & $-6.5$    & $-0.5$    & $1.5$ \\
         4 & 15.0   & 0.7       & 0.6   & 0.7       & $-1$      & 1 \\ [1ex]
    \end{tabular}
\end{table*}

In the baseline case study of the M\"{u}ller-Brown potential, $\kappa=2$ (the hyperparameter in~\eqref{eq:lcb}) is selected for the acquisition function to balance exploitation and exploration.
This value falls between two proposed scheduling approaches.
\cite{Srinivas_2012} proposes a schedule that provides cumulative regret guarantees.
The schedule is defined as:
\begin{equation}
    \label{eq:kappa_t_Srinivas_2012}
    \tag{$\kappa_{t,S}$}
    \kappa_{t,S} = \sqrt{2\log{\left(\frac{ M \, t^{2} \, \pi^{2}}{6 \, \delta}\right)}},
\end{equation}
with $M$ being the (finite) number of points, the search domain is discretized into, $t$ the current timestep of the outer BO loop, and $\delta$ a small positive number.
\cite{Srinivas_2012} note that the schedule may be too conservative (exploratory) in practice and, in their experiments, scale $\kappa_{t,S}$ down by a factor of $\sqrt{5}$ for improved performance.
For $\delta=0.1$, $M=\num{1E6}$ (near infinite), and after scaling down by $\sqrt{5}$, this would start around $\kappa_{1,S}=2.6$ and increase to $\kappa_{30,S} = 3.1$ over $30$ iterations: we call this schedule \ref{eq:kappa_t_Srinivas_2012}.
\cite{Kandasamy_2015} also argue this is still too conservative (exploratory) in practice and offer a schedule starting at $\kappa_{1,K}=0.3$ and increasing to $\kappa_{30,K}=1.6$, following
\begin{equation}
    \label{eq:kappa_t_Kandasamy_2015}
    \tag{$\kappa_{t,K}$}
    \kappa_{t,K} = \sqrt{0.2 D \log{(2t)}},
\end{equation}
with $D$ being the dimension of the search domain and $t$ again the current timestep of the outer BO loop.
Therefore, selecting $\kappa=2$ for our baseline case study provides a reasonable initial estimate, particularly given that the optimal balance is typically unknown \textit{a priori}.

Each experiment begins with three randomly chosen data points ($N=3$), sampled using Latin hypercube sampling (LHS)~\citep{McKay_1979}.
This approach mimics data collection following a real experimental setup, where only limited (but hopefully) well-distributed data are available at the start.
We terminate BO using the previously introduced termination criterion \eqref{eq:termination_critereon} with $\varepsilon_{x,1}=0.001$, $\varepsilon_{x,2}=0.05$, $\varepsilon_{f,r}=0.01$, and $\varepsilon_{f,a}=0.5$ (summarized in \cref{tab:testcase-setup-table}).

\subsubsection{Additional Case Studies}
\label{sec:additional_test_cases}
We further evaluate the performance of BO using ILS, IMS, and MAiNGO on the 2D Camelback benchmark function, selected for its two global minima, along with the 3D Ackley and 4D Hartmann functions, both characterized by a single global minimum and increasing problem dimensionality.
These benchmarks are standard functions from the Bingham test set~\citep{Surjanovic_2013} commonly used in BO literature.
Additionally, we include include instances from the GKLS test-function generator \citep{Gaviano_2003} in 2D, 3D, and 4D.
Although GKLS functions are not typically used in BO~\textemdash{}  largely because their landscapes feature abrupt minima and the functions are often not well-behaved~\textemdash{}  they offer tunable difficulty.

Exact formulations for each function are provided in \cref{sec:test-case-formulation}, with termination criterion thresholds detailed in \cref{sec:set-up-details}.
We use the same GP setup as in our base study (Mat\'{e}rn kernel, marginal log-likelihood, no noise).
Termination thresholds for all benchmark functions are provided in \cref{tab:testcase-setup-table}.

\section{Results}
\label{sec:results}
This section presents our empirical findings on how solver choice and initialization affect BO trajectories.
Our results address the two research questions posed in \cref{sec:introduction}:  (i) to which minimum the ``outer'' BO trajectory ultimately converges, and  (ii) given convergence to a globally near-optimal solution, how many iterations are required to reach it.

We organize the results so as to progressively isolate the influence of solver choice and choice of initial dataset (using \textit{runs}, \textit{experiments}, and \textit{case studies} as previously defined in~\cref{sec:approach} and summarized in~\cref{tab:notation_summary}):

\begin{itemize}
    \item In \cref{sec:single_run}, we compare the trajectories of the BO iterates using MAiNGO, ILS, and IMS in a single run of our baseline case study (defined by the M\"uller--Brown function \eqref{eq:muller-brown}, using the LCB acquisition function \eqref{eq:lcb} with $\kappa = 2$ and $N = 3$ initial points).
    This representative single-experiment, \textit{single-run} example highlights how the specific solvers directly influence the convergence behavior of BO.

    \item In \cref{sec:multipe_runs_same_initial_datasets}, we extend this comparison to multiple runs of the \textit{same experiment} to investigate solver variability under identical initial conditions, given the stochastic nature of ILS and IMS.

    \item We then proceed to investigate multiple runs of \textit{multiple experiments}, still within the same baseline case study, in \cref{sec:varying_initial_datasets} to explore how the initial dataset influences the probability of convergence to a globally near-optimal solution.

    \item In \cref{sec:likelihood_of_success}, we perform statistical testing to quantify the relative influence of solver choice and initial dataset on the likelihood of global success using the previously computed runs.

    \item Finally, in \cref{sec:question_one}, we synthesize our findings to answer our \hyperref[q:one]{first} research question regarding convergence destination.
\end{itemize}

We then turn to iteration efficiency~\textemdash{}  \textit{how fast} the outer BO trajectories reach a globally near-optimal solution.
In \cref{sec:iterations_to_convergence}, we compare the BO runs (across all experiments) that successfully converged to a globally near-optimal solution for all three solvers.
As our performance measure, we chose the number of iterations required to satisfy the defined termination criterion.
For additional insight, particularly for readers working with a fixed iteration budget, we include simple regret plots in \cref{sec:regret-plots}.
These plots illustrate the rate at which outer BO trajectories approach a global minimum, offering a complementary perspective on solver performance.
\cref{sec:more_case_studies} briefly presents additional case studies that explore different black-box functions to assess the generalizability of our findings.

\subsection{Convergence to a Globally Near-Optimal Minimum}
\label{sec:convergence_to_global_minimum}
\subsubsection{Demonstrative Single-Run Example: Analyzing Solver Influence in a Single Experiment}
\label{sec:single_run}
\begin{figure}[htbp]
    \begin{minipage}[t]{0.485\textwidth}
        \centering
        \includegraphics[width=\textwidth,trim={0cm 0cm 0cm 0cm},clip]
        {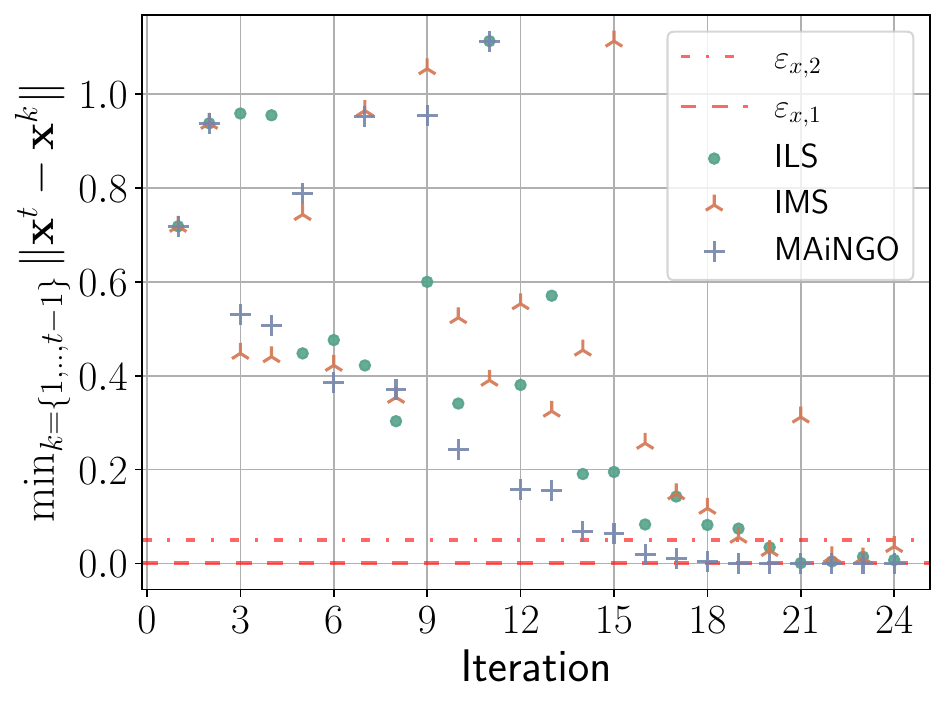}
        \subcaption{
            For this run, BO utilizing MAiNGO fulfills \eqref{eq:termination_critereon} after 16 iterations, BO with ILS after 20, and BO with IMS also after 20.
            If left to run past fulfilling \eqref{eq:termination_critereon}, BO utilizing ILS continues to explore previously uncharted regions.
            Dashed red lines indicate the predefined threshold for the distance termination criteria.
        }
        \label{fig:MB-3X-2Kappa-one_example:candidate_distance}
    \end{minipage}%
    \hfill%
    \begin{minipage}[t]{0.485\textwidth}
        \centering
        \includegraphics[width=\textwidth,trim={0cm 0cm 0cm 0cm},clip]
        {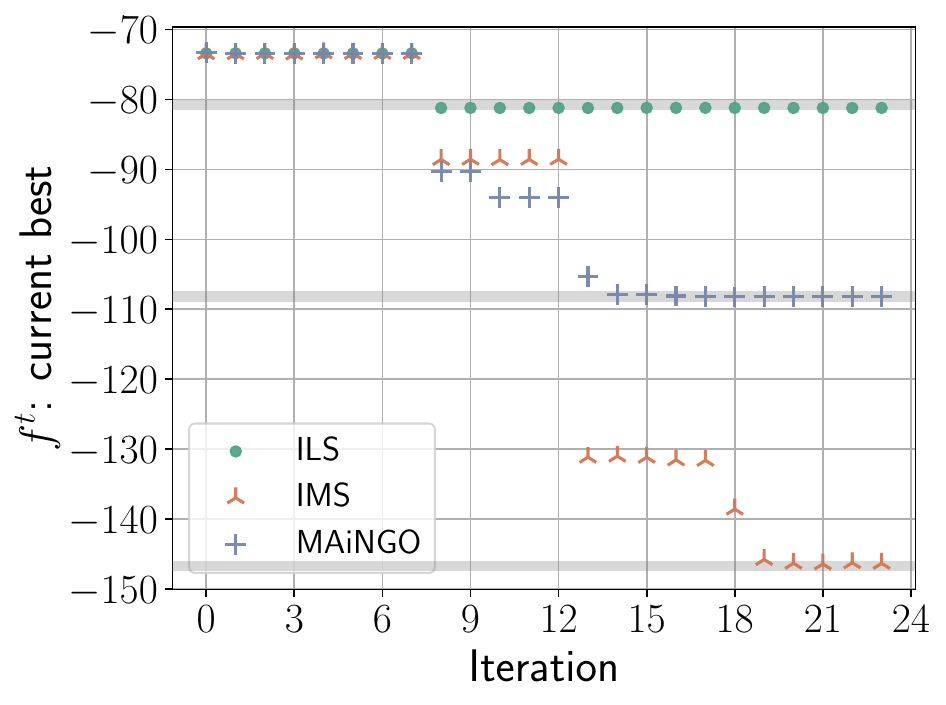}
        \subcaption{
            BO with ILS and with MAiNGO terminate at a local minimum, but BO with IMS finds a globally near-optimal solution.
            If left to run past fulfilling \eqref{eq:termination_critereon}, BO utilizing ILS also finds a globally near-optimal solution.
            Three thick gray horizontal lines mark the three minima.
        }
        \label{fig:MB-3X-2Kappa-one_example:best_minimum}
    \end{minipage}%
    \caption{
        BO trajectories for a single run in the baseline case study, with the black-box function defined by \eqref{eq:muller-brown}, an acquisition function using \eqref{eq:lcb} with $\kappa = 2$, and an initial dataset of $N = 3$ points.
    }
    \label{fig:MB-3X-2Kappa-one_example}
\end{figure}

This section presents a representative single-run example to illustrate \textit{how} and \textit{why} different solvers influence which minimum BO converges to.
\cref{fig:MB-3X-2Kappa-one_example} visualizes the trajectories of the outer BO iterates using MAiNGO, ILS, and IMS in a single run of our baseline case study (characterized by the M\"{u}ller-Brown case study \eqref{eq:muller-brown} and using \eqref{eq:lcb} with $\kappa = 2$ and $N = 3$ initial data points) for a single experiment.
We reiterate that in this baseline case study, BO convergence is defined by \eqref{eq:termination_critereon} with specific values found in \cref{tab:testcase-setup-table}.

\cref{fig:MB-3X-2Kappa-one_example:candidate_distance} tracks the distance between each candidate point and all previously sampled points.
Dashed red lines indicate the predefined threshold for the \textit{distance-based} portions of our termination criterion.
That is, our termination criterion is satisfied when either the bottom dashed line is reached (representing $\varepsilon_{x,1}$ in \ref{eq:TC-1}), or both the top dashed line is reached (representing $\varepsilon_{x,2}$ in \ref{eq:TC-2.1}) and the current candidate is close to the prior best point in \cref{fig:MB-3X-2Kappa-one_example:best_minimum} (representing \ref{eq:TC-2.2}).
In this run, all three BO trajectories converge to different minima: BO with IMS reaches a globally near-optimal solution of~\eqref{eq:muller-brown}, while BO with ILS and BO with MAiNGO terminate at local minima, as shown in \cref{fig:MB-3X-2Kappa-one_example:best_minimum} (thick gray bars denote the three (local and global) minima of \eqref{eq:muller-brown}).

Unlike optimization \textit{of the acquisition function} with a deterministic global solver, its optimization with ILS or IMS can result in locating one of its local minima.
Locating a local minimum of the acquisition function could result in continued exploration of a previously uncharted region of the black-box function (this is also possible \textit{after} the termination criterion is satisfied; c.f.,  \cref{fig:MB-3X-2Kappa-one_example:candidate_distance}).
While this additional exploration increases the likelihood of reaching a globally near-optimal solution (of the outer BO problem), it comes at the cost of additional BO iterations (and potentially past our termination criteria).
Moreover, repeated convergence to similar local minima can lead to indefinite iteration when relying solely on \ref{eq:TC-2.1} and \ref{eq:TC-2.2} in \eqref{eq:termination_critereon}.
To mitigate indefinite iterations, safeguards such as \ref{eq:TC-1} in \eqref{eq:termination_critereon} enforce timely termination, though they risk premature stopping if the new candidate is very close to a previously sampled suboptimal point (c.f.,  \cref{sec:bo_termination_critereon}).

\subsubsection{Multiple Runs of the Same Experiment: Analyzing Solver Variability}
\label{sec:multipe_runs_same_initial_datasets}
We reiterate that, due to the stochastic initialization of IMS and ILS, repeating a single experiment can result in different BO trajectories when using these solvers.
This section considers multiple runs of the same experiment of our baseline case study to demonstrate solver variability under identical initial conditions.
Across $80$ runs using \eqref{eq:lcb} with $\kappa = 2$ and the same $N = 3$ initial data points (i.e., consistent with the setup in \cref{sec:single_run}), BO utilizing MAiNGO never(!) finds a globally near-optimal solution of~\eqref{eq:muller-brown}, while utilizing IMS and ILS will find it with a probability of \SI{19}{\percent} and \SI{29}{\percent}, respectively.
The mean iterations to termination are $16.0 \pm 0$, $19.2 \pm 5.1$, and $21.2 \pm 6.5$ for BO with MAiNGO, IMS, and ILS, respectively ($\pm$ one standard deviation).
Among the runs that converged to a globally near-optimal solution, the mean iterations increased to $25.8 \pm 5.3$ for BO with IMS and $29.3 \pm 5.7$ for BO with ILS.

\subsubsection{Multiple Runs for Multiple Experiments: Analyzing Initial Dataset Influence}
\label{sec:varying_initial_datasets}
In \cref{sec:multipe_runs_same_initial_datasets}, we observed that BO with IMS and BO with ILS, due to their stochastic nature, exhibit different convergence trajectories across runs of a single experiment, while BO with MAiNGO consistently failed to find a globally near-optimal solution.
However, these findings are based on a fixed initial dataset, raising the question: \textit{How does the probability of convergence to a globally near-optimal solution change when the initial dataset varies?}

\begin{figure}[ht]
    \centering
    \includegraphics[width=0.5\textwidth,trim={0cm 0cm 0cm 0cm},clip]
    {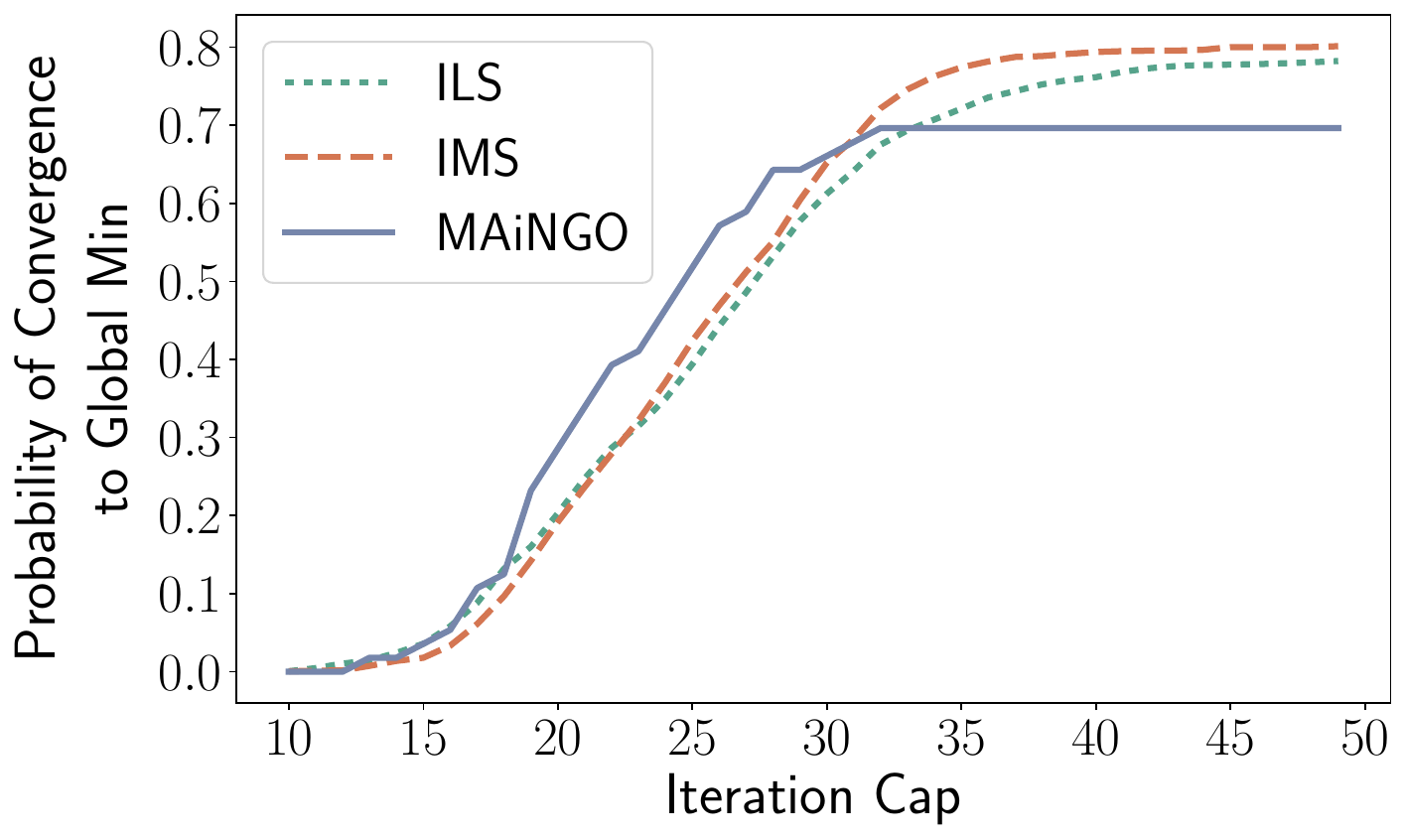}
    \caption{
        The probability of convergence to a globally near-optimal solution of the M\"{u}ller-Brown case study~\eqref{eq:muller-brown} utilizing the specified solver is comparable among all solvers if the maximum number of iterations is limited to $30$; beyond that, IMS dominates.
        The probability is computed as the fraction of runs that successfully converged (\textit{and} converged to a globally near-optimal solution) within the given iteration cap (limit), over the total number of runs for each solver.
    }
    \label{fig:MB-3t-2k-prob-to-success}
\end{figure}

To answer this question, we now investigate solver robustness across multiple experiments for our baseline case study using \eqref{eq:lcb} with $\kappa = 2$ and $N=3$ initial data points.
All initial datasets were generated via LHS.
\cref{tab:mb-prob} and \cref{fig:MB-3t-2k-prob-to-success} summarize the probabilities of convergence to a globally near-optimal solution of~\eqref{eq:muller-brown}.
As shown in \cref{fig:MB-3t-2k-prob-to-success}, all solvers exhibit comparable convergence probabilities when the iteration count is capped at $30$ iterations.
However, when no iteration cap is imposed, BO with IMS and ILS demonstrate a higher probability of reaching a globally near-optimal solution than BO with MAiNGO.
This finding suggests that BO with IMS and ILS often require more iterations than BO with MAiNGO to achieve convergence~\textemdash{}  a trend further explored in \cref{sec:iterations_to_convergence}.
The statistical significance of these probability differences is examined in \cref{sec:likelihood_of_success}.

Beyond solver selection, we found two key factors that influence the probability of convergence to a globally near-optimal solution of the black-box function, c.f., \cref{tab:mb-prob}:
\begin{itemize}
    \item   Increasing the size of the initial dataset provides a more representative surrogate model, i.e., GP, of the black-box function right from the start.
            This increases the likelihood of reaching a globally near-optimal solution of the black-box function.
    \item   Increasing the exploration parameter $\kappa$ in \eqref{eq:lcb} also boosts the probability of convergence to a globally near-optimal solution across all solvers.
\end{itemize}

\begin{table*}[htbp]
    \centering
    \caption{
        Probability of convergence to a globally near-optimal solution of the M\"{u}ller-Brown case study \eqref{eq:muller-brown}, aggregated over multiple experiments and runs.
        The $\kappa$ schedule \ref{eq:kappa_t_Srinivas_2012} is more exploratory than the baseline of $\kappa=2$; \ref{eq:kappa_t_Kandasamy_2015} is more exploitative.
        Details on the number of experiments and runs, and termination criteria are provided in \cref{sec:set-up-details}.
    }
    \label{tab:mb-prob}
    \begin{tabular}{c c c c c c c}
         $\kappa$ & $N$ & ILS & IMS & MAiNGO \\ [0.5ex]
         \hline
         2  & 3  & 0.79 & 0.80 & 0.70 \\
         3  & 3  & 0.94 & 0.97 & 0.96 \\
         2  & 10 & 0.83 & 0.83 & 0.81 \\
         2  & 20 & 0.92 & 0.94 & 0.90 \\
         \ref{eq:kappa_t_Srinivas_2012}  & 3 & 0.92 & 0.98 & 1.0 \\
         \ref{eq:kappa_t_Kandasamy_2015} & 3 & 0.53 & 0.53 & 0.51 \\
    \end{tabular}
\end{table*}

\subsubsection{A Statistical Look: Quantifying the Influence of the Solver and Initial Dataset}
\label{sec:likelihood_of_success}
Here, we address statistically (a) whether the differences seen in \cref{fig:MB-3t-2k-prob-to-success} between solvers are systematic or could be due to chance alone; and (b) how the influence of the solver compares to that of the initial dataset.
Note that the statements are for a given case study (all characterized by the same black-box and acquisition functions), as opposed to across families of functions.
We use the data of \cref{sec:varying_initial_datasets}, which consists of multiple runs across multiple experiments (characterized by different initial datasets).
For the statistical structure of the design of the data refer to \cref{sec:approach}.

\subsubsection*{Do Different Solvers Converge Equally Often, Averaging Over Datasets and Runs?}
We assessed whether $p^{s}$ differs across our three solvers ($s = \text{ILS, IMS, MAiNGO}$), by testing the hypothesis
\begin{equation}
    \label{eq:constancy}
    \tag{H}
    p^{\text{ILS}} = p^{\text{MAiNGO}} = p^{\text{IMS}}.
\end{equation}
To do this, we estimated the conditional logistic regression of the model that assumes
\begin{equation}
    \label{eq:ormodel}
    \begin{aligned}
        \mathrm{odds} \ p^{\text{ILS}}_{dat} / \mathrm{odds} \ p^{\text{MAiNGO}}_{dat} &= \exp(\alpha) \\
        \mathrm{odds} \ p^{\text{IMS}}_{dat} / \mathrm{odds} \ p^{\text{MAiNGO}}_{dat} &= \exp(\beta),
    \end{aligned}
\end{equation}
using conditional Maximum Likelihood Estimation (cMLE)~\citep{Andersen_1970}; where odds $(p):=p/(1-p)$ and $\alpha, \beta$ are the common log-odds ratio of ILS and IMS vs MAiNGO, respectively.
This model is often used to estimate odd ratios
of probabilities between groups when these ratios are assumed to be common across the matching factor (initial datasets).
Estimation of this model is also useful in our case because (a) it addresses the matching of solvers' runs on experiments (datasets); (b) if solvers are overall equivalent, i.e., \eqref{eq:constancy} holds, then the cMLE estimates of $\alpha$ and $\beta$ are known to center at zero (a proof is given in \cref{sec:stats_proof}) even if the TRUE odds ratios in the left-hand sides of \eqref{eq:ormodel} differ arbitrarily across datasets.

The estimated average probabilities $p^{s}$ are \SI{80}{\percent} for BO using IMS, \SI{79}{\percent} for ILS, and \SI{70}{\percent} for MAiNGO
($p-test \text{ value}=0.036$ for the comparison that all three are statistically equivalent).
We answer our question~\textemdash{}  do different solvers converge equally often?~\textemdash{}  by pairwise comparisons, which suggest that the probability of convergence of BO with IMS or ILS are statistically different from BO using MAiNGO, whereas BO with IMS or ILS are not statistically significantly different from each other.

\subsubsection*{Comparing the Importance of the Initial Dataset with that of the Solver}
To compare the influence of the initial dataset with that of the solver, we use minimax tests and refer the reader to the structure of the design described in \cref{sec:approach}.
Consider a list with all possible subsets $ds$, of half (28) of the 56 datasets $dat$.
For each specific such $ds$, and for each solver $s$, we calculate the empirical proportion $\hat{p}^{s}_{ds}$ of convergence over the runs and the datasets in $ds$.
We ask the following questions.
\begin{enumerate}
    \item \textbf{Question:} If we choose a solver $s$, and then ``arbitrarily choose a subset $ds$ from the above list,'' what is the worst that can happen, i.e., what is the subset $ds$ that minimizes the convergence probability $\hat{p}^{s}_{ds}$?
          From an optimization perspective, this question can be formulated as a discrete max-min problem that reads $\max_{s}\min_{ds} \hat{p}^{s}_{ds}$.\\
          \textbf{Answer: }
          We calculated these minima as \SI{68}{\percent} (ILS), \SI{54}{\percent} (MAiNGO), and \SI{69}{\percent} (IMS).
          So, if we can only choose the solver to be used in BO, we should choose IMS to maximize (over solvers) the worst possible convergence (over subsets given solver) and assure convergence is $\geq \SI{69}{\percent}$.
    \item \textbf{Question:} If we choose a subset $ds$, and ``arbitrarily choose a solver,'' what is the worst that can happen, i.e., what is the solver that minimizes the convergence probability?
          Note that this questions corresponds to the discrete max-min problem $\max_{ds}\min_{s} \hat{p}^{s}_{ds}$.\\
          \textbf{Answer: }
          The minimum convergence probability over solvers for a given subset $ds$ ranges from \SI{54}{\percent} -- \SI{86}{\percent} over $ds$, and their maximum \SI{86}{\percent} occurs at a particular dataset, say $ds^{*}$, when using, as it happens, again IMS.
          So, if we can only choose a subset $ds$, we should choose the subset $ds^{*}$ that maximizes (over subsets $ds$) the worst possible convergence (over solvers given subset $ds$) and assure that we have convergence $\geq \SI{86}{\percent}$.
\end{enumerate}
This analysis suggests that the choice of the dataset is at least as important as the choice of the solver to achieve a high probability of convergence.

\subsubsection{Reconciling Solver Influence and Dataset Dependence}
\label{sec:question_one}
We return to our \hyperref[q:one]{first} question from \cref{sec:introduction}: How does the choice of solver affect successful convergence to a globally near-optimal solution, and to what extent is that success dependent on the initial dataset?
Our analysis shows BO utilizing MAiNGO does not always converge to a globally near-optimal solution of the black-box function and can terminate at a local minimum.
Among the solvers tested, BO using IMS or ILS demonstrates statistically significantly higher probabilities of finding a globally near-optimal solution than BO utilizing MAiNGO (c.f.,  \cref{tab:mb-prob} and cMLE results).
Furthermore, the observed heterogeneity highlights that, while BO utilizing MAiNGO and IMS differ statistically in their performance (as compared to each other), the initial dataset plays a critical role in determining to which minimum (local or global) BO converges towards, as we have shown for this specific case study using minimax tests.
It seems reasonable to postulate that convergence to a local minimum of the black-box function occurs when the initial dataset is not representative of the true objective function, leading to a poorly constructed initial GP surrogate, combined with an acquisition function that insufficiently favors exploration.
BO with IMS or ILS can compensate for insufficient exploration and a poor initial surrogate through stochasticity; BO with MAiNGO cannot.
In other words, suboptimal optimization of poorly chosen acquisition functions can be preferable to their optimal solution.

Incorporating stochastic elements in the optimization process parallels  using a stochastic gradient descent solver in the training phase of, e.g., deep-learning \citep{Zhou_2019} and physics-informed neural networks \citep{Markidis_2021}.
In machine learning, stochastic gradient descent solvers are widely used to avoid or escape local minima of the loss function during training \citep{Kleinberg_2018}.
Therein, if only a local solver, such as L-BFGS-B, is used in the training phase, the training is prone to termination at a local minimum \citep{Markidis_2021}.
Therefore, a hybrid approach is typically applied in the training phase of, e.g., physical-informed neural networks: First, a stochastic gradient descent solver is used to avoid (premature) termination at local minima as stochastic gradient descent solver are capable of ``escaping'' local minima \citep{Kleinberg_2018}.
Because stochastic gradient descent solvers converge fast to the vicinity of a global optimum, but the iterates then will stay around this global optimum with constant probability \citep{Kleinberg_2018}, subsequently: second, the solution is refined using a local solver, e.g., L-BFGS-B \citep{Markidis_2021}.

It is important to note that these findings are based on our single baseline case study: various datasets consisting of three initial points and \eqref{eq:lcb} with $\kappa = 2$ as the acquisition function.
Significantly increasing the number of initial data points (making an exploitative acquisition function more appropriate) or increasing the exploration factor $\kappa$ in \eqref{eq:lcb} improves the probability of convergence to a globally near-optimal solution of the black-box function (\cref{tab:mb-prob}).
Unfortunately, selecting an appropriate exploration factor for \eqref{eq:lcb} is challenging without prior domain knowledge, and increasing the number of initial samples is often constrained by physical, time, or cost limitations.
These findings emphasize the importance of developing strategies for selecting (a) initial sampling points to construct a high-quality initial dataset \textit{before} entering the BO loop (see, e.g.,~\citep{bull_2011,Ren_2024,Wang_2024}), (b) more suitable acquisition functions, and (c) methods to detect poorly performing acquisition functions early in the optimization process or ``self-adaptive'' acquisition functions (see e.g.,~\citep{Yan_2022,Chen_2023,Benjamins_2023,Maitra_2024}).

\subsection{Iterations to Convergence}
\label{sec:iterations_to_convergence}
Having established how solver choice and initial dataset impact the probability of convergence to a globally near-optimal solution, we now investigate efficiency (question \hyperref[q:two]{two} posed in \cref{sec:introduction}): Given a well-chosen initial dataset, which solver within BO converges, as defined by~\eqref{eq:termination_critereon},  in the fewest iterations?

The following analysis focuses exclusively on runs that successfully reached a globally near-optimal solution, grouped by experiment (i.e., initial dataset).
Experiments where any solver failed to achieve at least one successful run are excluded.
This approach allows for a controlled comparison of iteration counts until convergence, assuming \textit{a priori} knowledge of success while accounting for dataset effects.
Including all runs~\textemdash{}  regardless of which minimum they reached~\textemdash{} i.e., using a broader dataset, would introduce a bias: solvers with lower convergence probabilities tend to terminate earlier, requiring fewer iterations, while those more likely to reach a globally near-optimal solution explore further and take longer.
For completeness, statistics with this broader dataset are provided in \cref{sec:stats-with-all-min}.

\begin{figure}[htbp]
    \begin{minipage}[t]{0.31\textwidth}
        \centering
        \includegraphics[trim={0cm 0cm 0cm 0cm}, clip,width=\textwidth]
        {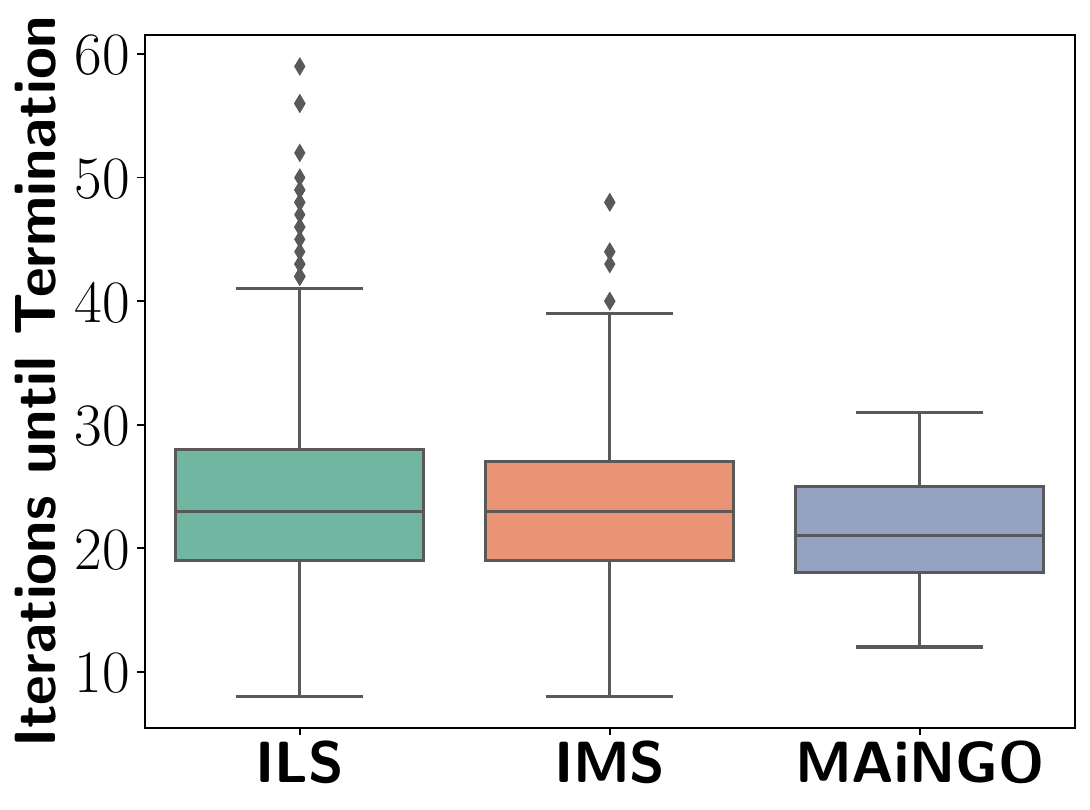}
        \subcaption{
            $\kappa = 2$ and $N = 3$.
        }
        \label{fig:MB-3t-2k-conv:k2_N2}
    \end{minipage}%
    \hspace{\fill}%
    \begin{minipage}[t]{0.31\textwidth}
        \centering
        \includegraphics[trim={0cm 0cm 0cm 0cm}, clip,width=\linewidth]
        {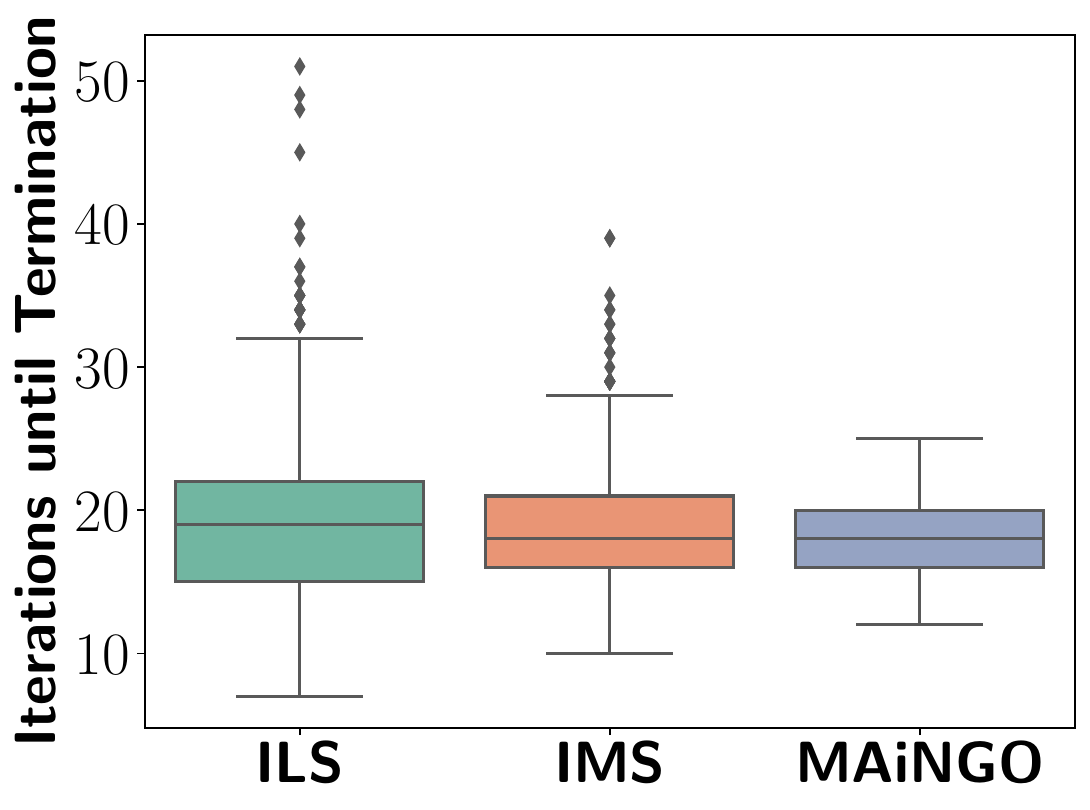}
        \subcaption{
            $\kappa = 2$ and $N = 10$.
        }
        \label{fig:MB-3t-2k-conv:k2_N10}
    \end{minipage}%
    \hspace{\fill}%
    \begin{minipage}[t]{0.31\textwidth}
        \centering
        \includegraphics[trim={0cm 0cm 0cm 0cm}, clip,width=\linewidth]
        {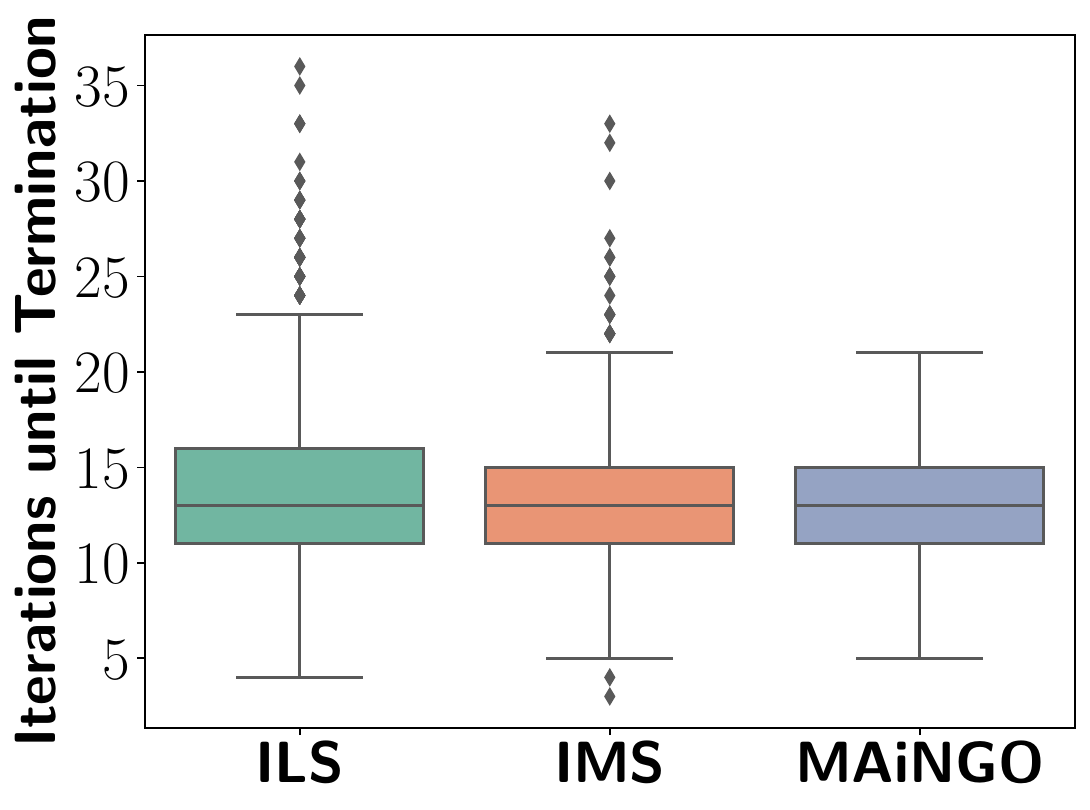}
        \subcaption{
            $\kappa = 2$ and $N = 20$.
        }
        \label{fig:MB-3t-2k-conv:k2_N20}
    \end{minipage}%
    \par\medskip
    \begin{minipage}[t]{0.31\textwidth}
        \centering
        \includegraphics[trim={0cm 0cm 0cm 0cm}, clip,width=\linewidth]
        {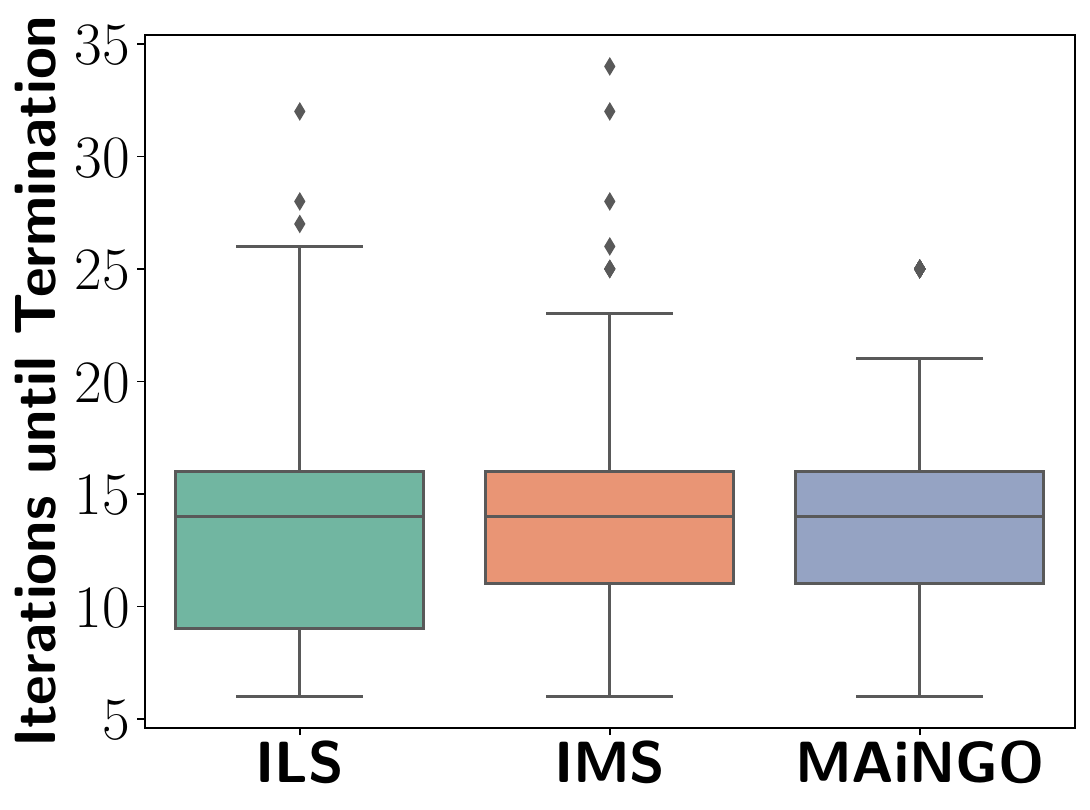}
        \subcaption{
            $\kappa = $~\ref{eq:kappa_t_Kandasamy_2015} and $N = 3$.
        }
        \label{fig:MB-3t-2k-conv:kK_N3}
    \end{minipage}%
    \hspace{\fill}%
    \begin{minipage}[t]{0.31\textwidth}
        \centering
        \includegraphics[trim={0cm 0cm 0cm 0cm}, clip,width=\linewidth]
        {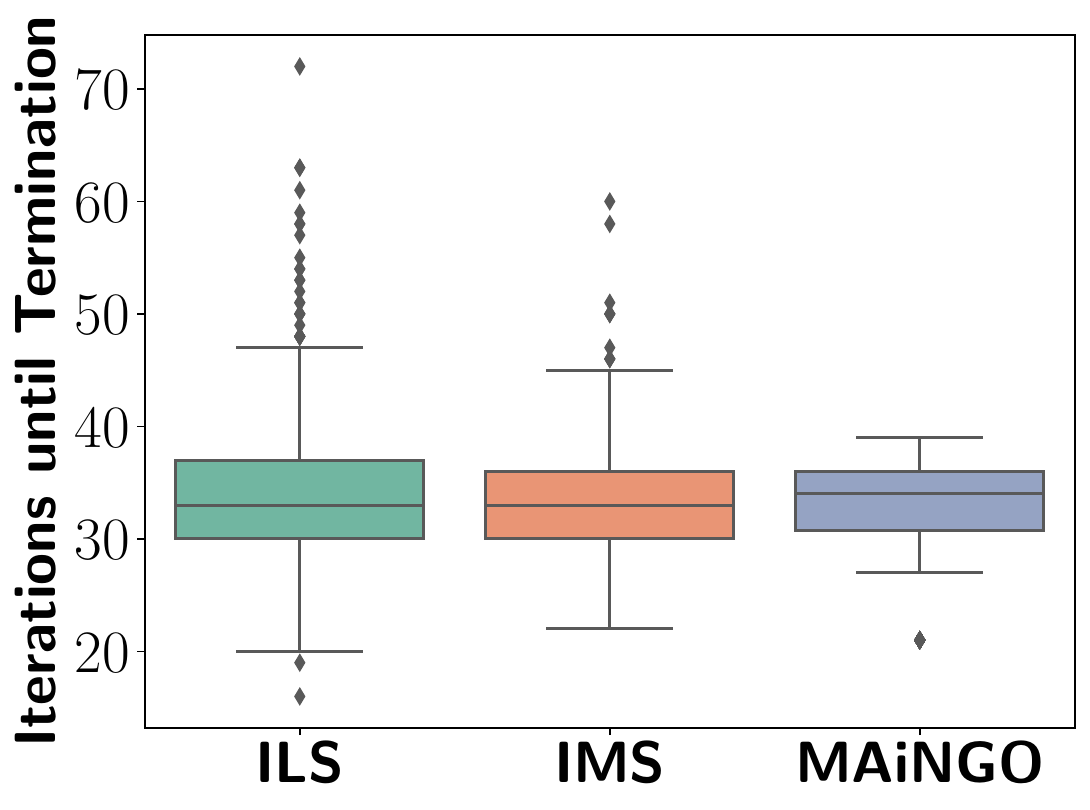}
        \subcaption{
            $\kappa = 3$ and $N = 3$.
        }
        \label{fig:MB-3t-2k-conv:k3_N3}
    \end{minipage}%
    \hspace{\fill}%
    \begin{minipage}[t]{0.31\textwidth}
        \centering
        \includegraphics[trim={0cm 0cm 0cm 0cm}, clip,width=\linewidth]
        {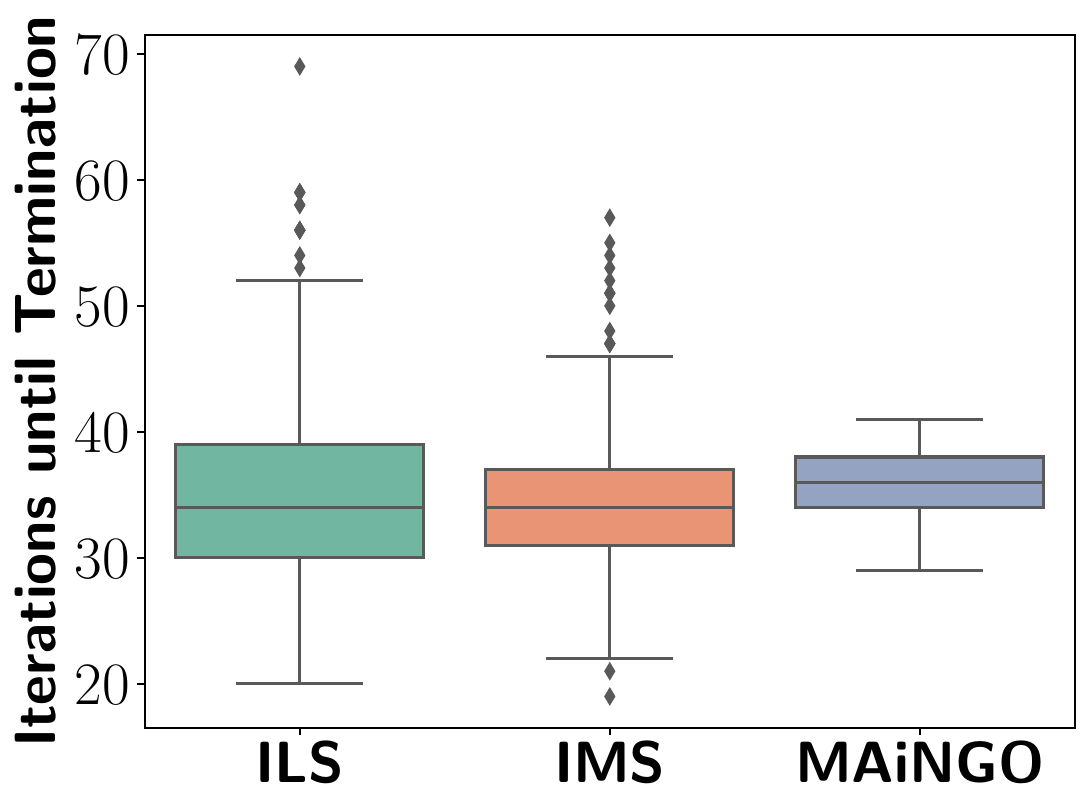}
        \subcaption{
            $\kappa =$~\ref{eq:kappa_t_Srinivas_2012} and $N = 3$.
        }
        \label{fig:MB-3t-2k-conv:kS_N3}
    \end{minipage}
    \caption{
        Distribution of iteration counts for runs that successfully converged to a globally near-optimal solution, considering only experiments where all three solvers have at least one successful run.
        Results are shown for a series of M\"{u}ller-Brown case studies.
        Details on the number of experiments, runs, and termination criteria are provided in \cref{sec:set-up-details}.
    }
    \label{fig:MB-3t-2k-conv}
\end{figure}

In this restricted analysis, conducted for our baseline case study (using \eqref{eq:lcb} with $\kappa = 2$ and $N = 3$ initial data points), BO utilizing MAiNGO required fewer iterations on average to satisfy~\eqref{eq:termination_critereon} compared to both BO utilizing IMS and BO utilizing ILS (\cref{fig:MB-3t-2k-conv:k2_N2}).
Statistical comparisons using one-sided t-tests (paired by experiments) confirm that this difference is significant, demonstrating MAiNGO's advantage in terms of mean BO iterations to convergence (c.f., \cref{tab:mb-fast}).
The same trends hold when the number of initial points in the dataset is increased to $10$ or $20$ (when a more exploitative acquisition function is appropriate), but not when $\kappa$ is increased to $3$ or $\kappa_{t,S}$ is used (more exploratory acquisition functions), c.f., \cref{fig:MB-3t-2k-conv:k3_N3,fig:MB-3t-2k-conv:k2_N10,tab:mb-fast}.
These results are consistent with \cite{Xie_2025}, who demonstrate the effectiveness of PK-MIQP in exploitative scenarios with  $N \geq 10$.
Beyond achieving a lower mean number of iterations as compared to BO with ILS or IMS, BO utilizing MAiNGO also exhibits reduced variability by eliminating the stochasticity associated with the solvers.
This finding is consistent across all case studies (all tested $\kappa$ and initial dataset sizes).

For an alternative perspective, we present simple regret plots, which illustrate how quickly BO trajectories approach a globally near-optimal solution for each solver, in \cref{sec:regret-plots}: these plots reinforce the aforementioned findings.

\begin{table}[ht]
    \centering
    \caption{
        Summary of convergence statistics for each case study, considering only runs that successfully reached a globally near-optimal solution.
        Experiments are included only if at least one run converged to a globally near-optimal solution for all three solvers.
        Number of runs (\textit{count}), mean, median, and standard deviation of iterations to convergence are reported for this subset of data.
        One-sided t-tests, \textit{paired by experiment}, evaluate whether MAiNGO converges in fewer iterations than IMS or ILS.
        MAiNGO converges in fewer iterations than ILS for all case studies and significantly fewer iterations than IMS for all but the case study using the exploratory acquisition function ($\kappa=3$, $\kappa=$~\ref{eq:kappa_t_Srinivas_2012}).
        Details on the total number of experiments, runs, and termination criteria for each case study are provided in \cref{sec:set-up-details}.
    }
    \label{tab:mb-fast}
    \begin{tabular}{c c c c c c c c c}
        $\kappa$ & $N$ & Solver & Count & Mean & Median & Std & t-test & p-value \\ [0.5ex]
        \hline
        2 & 3   & ILS   & 1085  & 23.8 & 23.0 & 7.24 & 4.11 & \num{1.0E-4} \\
        2 & 3   & IMS   & 1140  & 23.3 & 23.0 & 5.47 & 3.84 & \num{2.3E-4} \\
        2 & 3   & MAiNGO& 1209  & 21.3 & 21.0 & 4.50 & - & -\\ [1ex]
        \hline
        3 & 3   & ILS   & 703   & 34.3 & 33.0 & 6.43 & 1.22 & \num{1.2E-1}\\
        3 & 3   & IMS   & 722   & 33.4 & 33.0 & 4.55 & -0.56 & \num{7.1E-1}\\
        3 & 3   & MAiNGO& 744   & 33.2 & 34.0 & 4.14 & -    & -\\
        \hline
        2 & 10  & ILS   & 738   & 19.3 & 19.0 & 5.78 & 2.48 & \num{1.0E-2}\\
        2 & 10  & IMS   & 760   & 18.8 & 18.0 & 4.26 & 1.54 & \num{6.8E-2}\\
        2 & 10  & MAiNGO& 775   & 18.1 & 18.0 & 3.06 & -    & -\\
       \hline
        2 & 20  & ILS   & 1093  & 13.8 & 13.0 & 4.91 & 2.89 & \num{3.2E-3}\\
        2 & 20  & IMS   & 1138  & 13.1 & 13.0 & 3.71 & 0.47 & \num{3.2E-1}\\
        2 & 20  & MAiNGO& 1147  & 13.0 & 13.0 & 3.37 & -    & -\\
        \hline
        \ref{eq:kappa_t_Srinivas_2012} & 3 & ILS   & 404 & 35.3 & 34.0 & 7.28 & -0.73 & \num{7.6E-1}\\
        \ref{eq:kappa_t_Srinivas_2012} & 3 & IMS   & 431 & 33.9 & 34.0 & 5.34 & -2.21 & \num{9.8E-1}\\
        \ref{eq:kappa_t_Srinivas_2012} & 3 & MAiNGO& 441 & 35.4 & 36.0 & 3.46 & -     & -\\
        \hline
        \ref{eq:kappa_t_Kandasamy_2015} & 3 & ILS   & 496 & 13.3 & 14.0 & 4.57 & 0.28 & \num{3.9E-1}\\
        \ref{eq:kappa_t_Kandasamy_2015} & 3 & IMS   & 529 & 13.7 & 14.0 & 4.55 & 2.20 & \num{1.9E-2}\\
        \ref{eq:kappa_t_Kandasamy_2015} & 3 & MAiNGO& 546 & 13.3 & 14.0 & 4.56 & -    & -\\ [1ex]
    \end{tabular}
\end{table}

\subsection{Additional Case Studies}
\label{sec:more_case_studies}
\begin{figure}[htbp]
    \centering
    \begin{minipage}[t]{0.303\textwidth}
        \centering
        \includegraphics[trim={0cm 0cm 0cm 0cm}, clip,width=\textwidth]
        {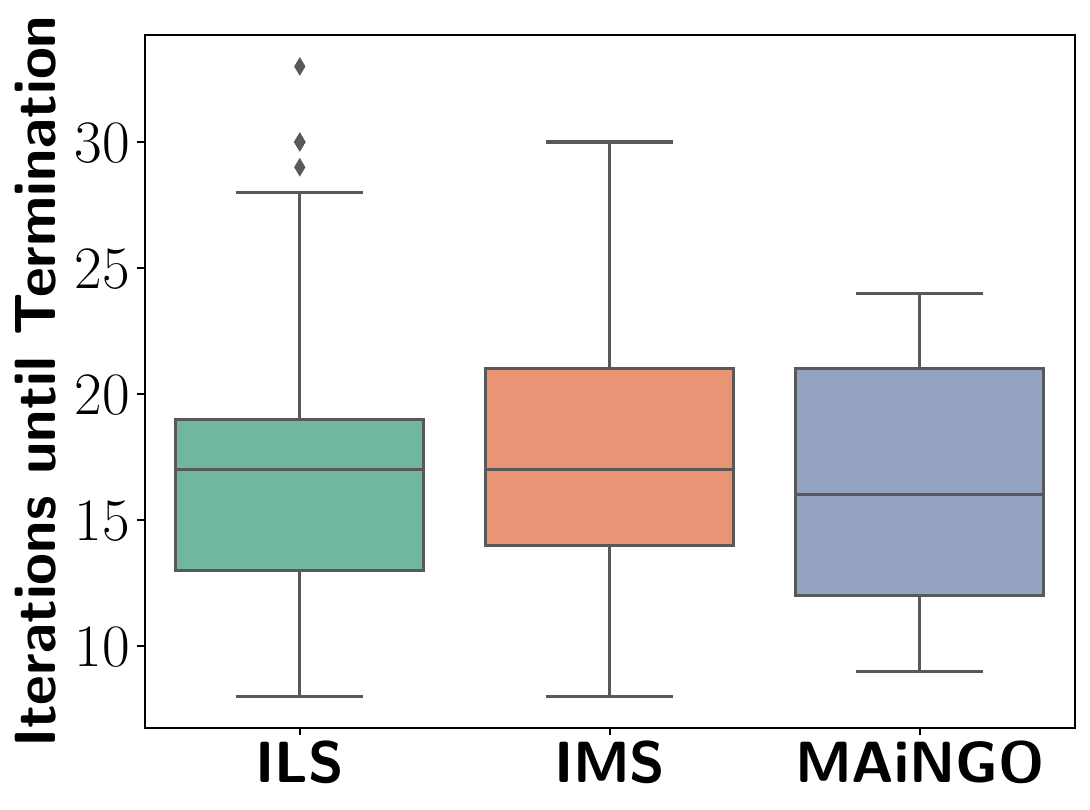}
        \subcaption{
            2D Six-Hump Camelback
        }
        \label{fig:CB:k2_N3}
    \end{minipage}%
    \hfill%
    \begin{minipage}[t]{0.303\textwidth}
        \centering
        \includegraphics[trim={0cm 0cm 0cm 0cm}, clip,width=\linewidth]
        {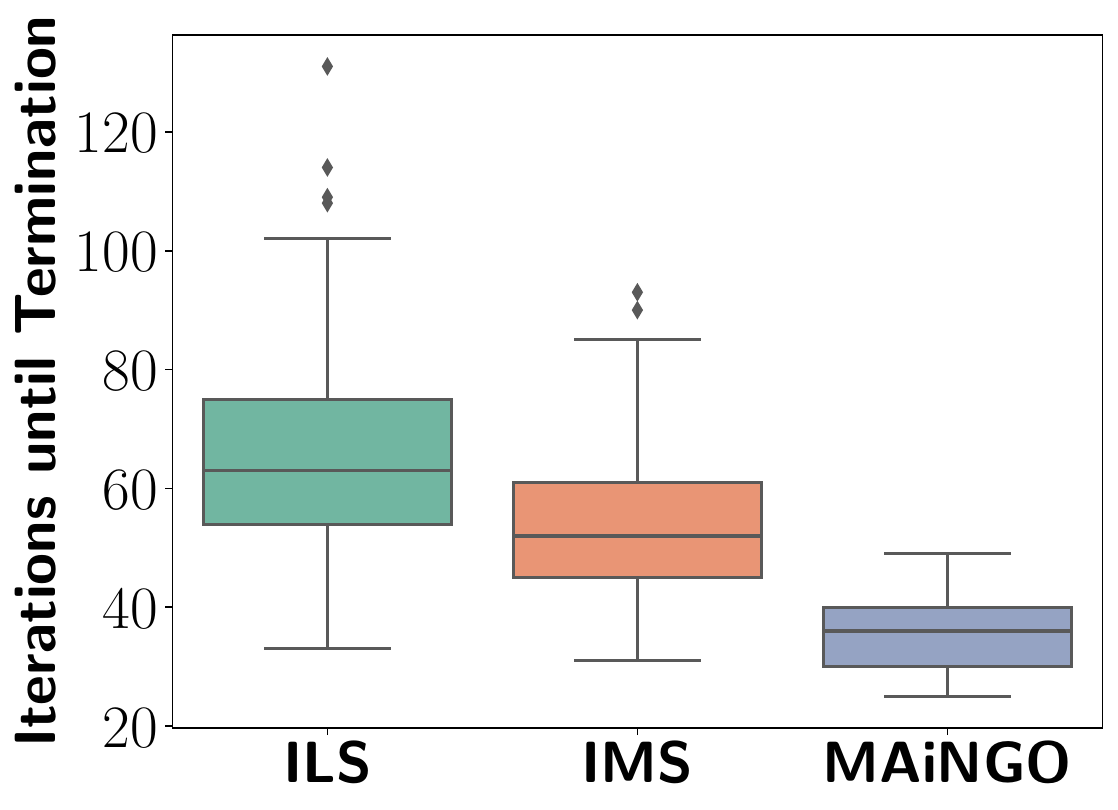}
        \subcaption{
            3D Ackley
        }
        \label{fig:A3D:k2_N4}
    \end{minipage}%
    \hfill%
    \begin{minipage}[t]{0.303\textwidth}
        \centering
        \includegraphics[trim={0cm 0cm 0cm 0cm}, clip,width=\linewidth]
        {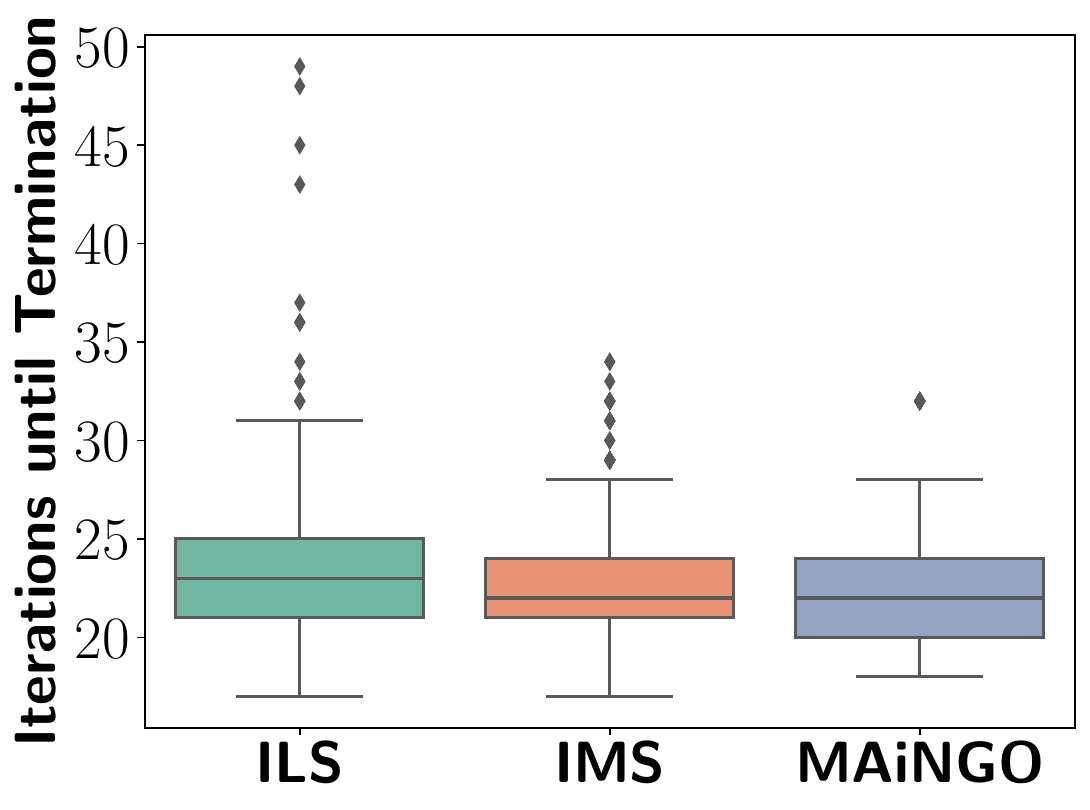}
        \subcaption{
            4D Hartmann
        }
        \label{fig:H4D:k2_N5}
    \end{minipage}%
    \par\medskip
    \begin{minipage}[t]{0.303\textwidth}
        \centering
        \includegraphics[trim={0cm 0cm 0cm 0cm}, clip,width=\linewidth]
        {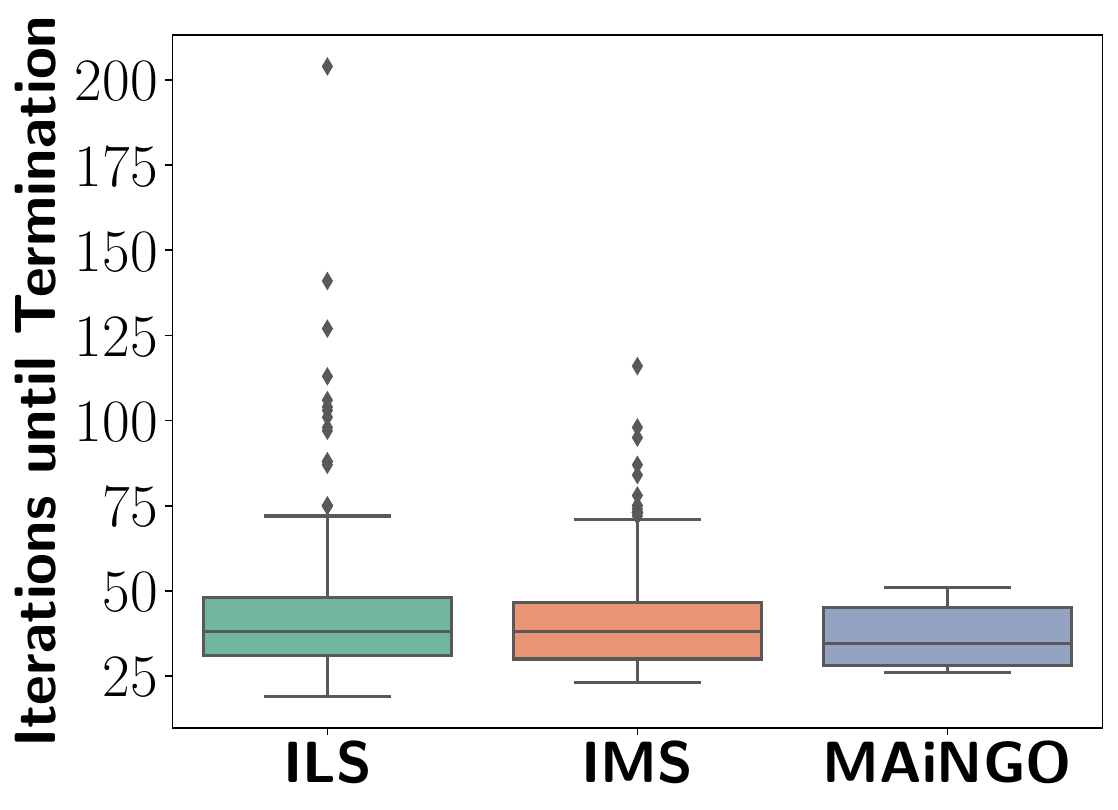}
        \subcaption{
            2D GKLS
        }
        \label{fig:2GKLS:k3_N3}
    \end{minipage}%
    \hspace{1.5em}%
    \begin{minipage}[t]{0.303\textwidth}
        \centering
        \includegraphics[trim={0cm 0cm 0cm 0cm}, clip,width=\linewidth]
        {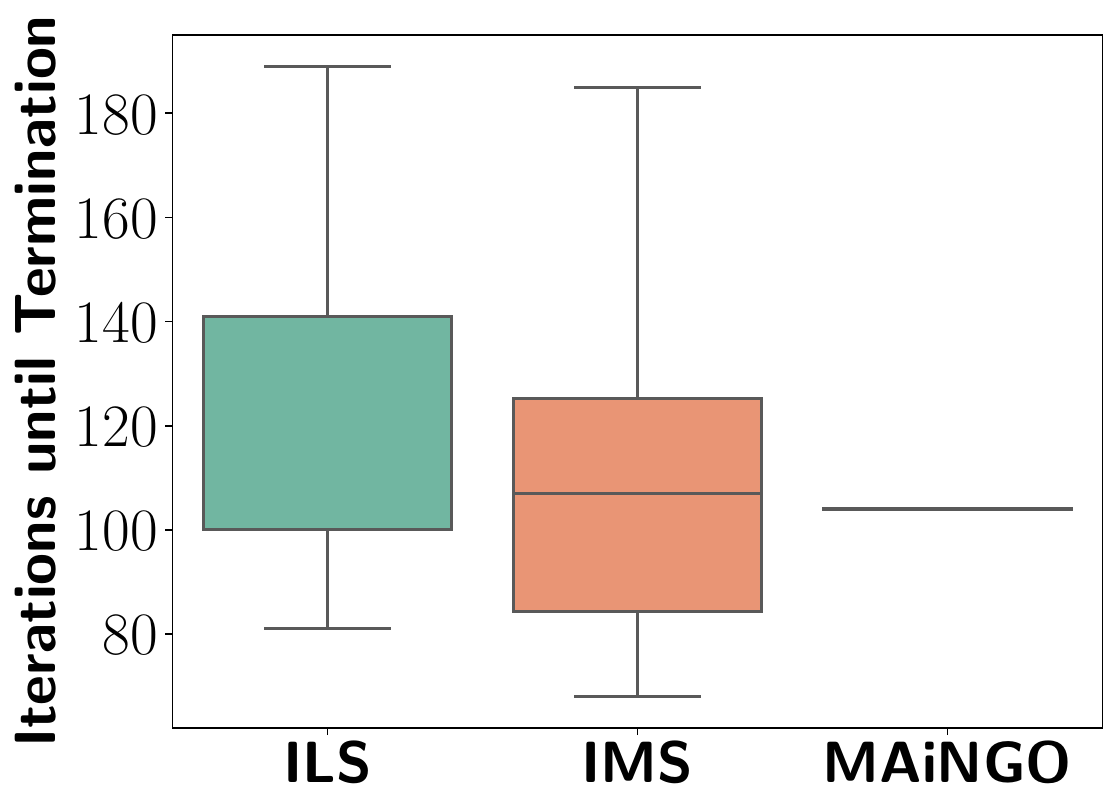}
        \subcaption{
            3D GKLS
        }
        \label{fig:3GKLS:k3_N4}
    \end{minipage}%
    \hfill
    \caption{
        Distribution of iteration counts for runs that successfully converged to a globally near-optimal solution,
        considering only experiments where all three solvers have at least one successful run.
        Results are shown for a series of case studies growing in dimensions.
        Details on the number of experiments, runs, termination criteria, $\kappa$ values, and number of points in the initial dataset are provided in \cref{sec:set-up-details}.
    }
    \label{fig:other-test-cases-conv}
\end{figure}

We tested BO with ILS, IMS, and MAiNGO on six additional benchmark functions, i.e., 2D Six-Hump Camelback, 3D Ackley, and 4D Hartmann as well as GKLS functions with increasing dimensionality to evaluate generalizability.
As in our baseline case, we use \eqref{eq:lcb} as the acquisition function with $\kappa=2$ for all but the GKLS functions, where we use $\kappa=3$.
\cref{sec:test-case-formulation} provides the formulations for these functions, while \cref{sec:set-up-details}  outlines the termination thresholds, number of experiments, and runs per experiment.
\cref{fig:other-test-cases-conv} and \cref{tab:othertestcases-results} present the iterations to convergence distributions and statistical comparisons using one-sided t-tests, paired by experiment.

\begin{table}[ht]
    \centering
    \caption{
        Summary of convergence statistics for various objective function case studies.
        The probability of convergence to a globally near-optimal solution (\textit{prob}) is computed over all experiments and runs.
        The remaining statistics~\textemdash{}  the number of successful runs (\textit{count}), and for those runs, the mean, median, and standard deviation of iterations to convergence~\textemdash{}  are reported only for experiments where at least one run successfully converged to a globally near-optimal solution for all three solvers.
        This approach controls for dataset effects.
        One-sided t-tests, \textit{paired by experiment}, evaluate whether BO with MAiNGO converges in fewer iterations than BO with IMS or ILS.
        BO with MAiNGO converges in fewer iterations than with IMS for all case studies, and significantly fewer iterations than with ILS for all but the 2D Camelback.
        Note that in the 4D GKLS benchmark, BO never found the global optimum, despite which inner solver was used.
        Details on the total number of experiments, runs, and termination criteria for each case study are provided in \cref{sec:set-up-details}.
    }
    \label{tab:othertestcases-results}
    \begin{tabular}{l c c c c c c c c}
         Test Function & Solver & Prob & Count & Mean & Median & Std & t-test & p-value \\ [0.5ex]
         \hline
         \multirow{3}{*}{2D Camelback}
         & ILS  & 1  & 1271 & 16.5 & 17.0 & 3.99 & -0.05 & \num{5.2E-1} \\
         & IMS  & 1   & 1271 & 17.6 & 17.0 & 4.65 & 2.23 & \num{1.6E-2} \\
         & MAiNGO & 1 & 1271 & 16.5 & 16.0 & 4.60 & - & - \\ [1ex]
         \hline
         \multirow{3}{*}{3D Ackley}
         & ILS  & 1    & 465 & 65.0 & 63.0 & 15.33 & 38.6 & \num{6.68E-19} \\
         & IMS  & 1   & 465 & 53.5 & 52.0 & 11.1 & 30.8 & \num{4.34E-17} \\
         & MAiNGO & 1 & 465 & 35.3 & 36.0 & 6.31 & - & - \\ [1ex]
         \hline
         \multirow{3}{*}{4D Hartmann}
         & ILS  & 0.80    & 334 & 23.5 & 23.0 & 4.17 & 2.31 & \num{1.6E-2} \\
         & IMS  & 0.74   & 350 & 23.1 & 22.0 & 3.50 & 2.95 & \num{3.9E-3} \\
         & MAiNGO & 0.73 & 352 & 22.5 & 22.0 & 3.34 & - & - \\ [1ex]
         \hline
         \multirow{3}{*}{2D GKLS}
         & ILS  & 0.574  & 257 & 44.0 & 38.0 & 21.06 & 0.845 & \num{1.99E-1} \\
         & IMS  &  0.593  & 291 & 41.3 & 38.0 & 14.95 & 0.931 & \num{1.83E-1} \\
         & MAiNGO & 0.560 & 350 & 35.9 & 34.5 & 8.86 & - & - \\ [1ex]
         \hline
         \multirow{3}{*}{3D GKLS}
         & ILS  & 0.058  & 5 & 130.4 & 141.0 & 37.48 & - & - \\
         & IMS  &  0.102  & 6 & 110.3 & 107.0 & 31.72 & - & - \\
         & MAiNGO & 0.040 & 25 & 104.0 & 104.0 & 0 & - & - \\ [1ex]
         \hline
        \multirow{4}{*}{4D GKLS}
         & ILS  & 0  & - & - & - & - & - & - \\
         & IMS  &  0  & - & - & - & - & - & - \\
         & MAiNGO & 0 & - & - & - &- & - & - \\ [1ex]
    \end{tabular}
\end{table}

The 2D Six-Hump Camelback function has two global minima, though we did not differentiate between them.
In the BO context, MAiNGO performed comparably to ILS on this function, while IMS required significantly more iterations to reach~\eqref{eq:termination_critereon}.
For higher-dimension functions such as the 3D Ackley and 4D Hartmann, \textit{BO with MAiNGO consistently converged in significantly fewer iterations} than BO with ILS or IMS, \textbf{\textit{suggesting BO utilizing ILS and IMS are more prone to getting stuck in local optima as dimensionality increases}}.
Results for the 4D Hartmann function again highlight the importance of initial dataset quality; BO utilizing MAiNGO, as compared to BO utilizing ILS or IMS, was the least likely to converge to a globally near-optimal solution if the initial dataset (or acquisition function and $\kappa$ value) was poorly selected.

We also evaluated solver performance on the GKLS test functions in 2D, 3D, and 4D.
Note these functions are challenging for BO because the global minimum basin is typically narrow and abrupt (i.e., not a smooth extension of the bulk behavior of the function).
In 2D GKLS, BO with all three solvers achieve moderate probabilities of convergence, with MAiNGO still exhibiting a lower probability of success than ILS and IMS (\cref{tab:othertestcases-results}).
Conditional on success, MAiNGO converges in fewer iterations on average, consistent with the trend observed on the other benchmarks.
The 3D GKLS instance is substantially more challenging.
Across all experiments, BO with MAiNGO has the lowest overall probability of converging to a globally near-optimal solution, as it succeeds only for a single initial dataset (and, being deterministic, succeeds in all 25 runs of that experiment).
In contrast, BO with ILS and IMS succeed in only 5 and 6 runs, respectively, for that same experiment, although, unlike BO with MAiNGO, they are sometimes able to locate the near global optimum in cases where the initialization is poor.
Thus, when~\textemdash{}  and only when~\textemdash{}  the initialization is favorable, BO with MAiNGO is both more reliable and more efficient than the stochastic solvers for the evaluated value of~$\kappa$.
Because only one experiment yields at least one successful run for all solvers, paired statistical tests cannot be computed, and no $t$- nor $p$-values can be calculated for 3D GKLS.
For 4D GKLS, none of the solvers successfully locate a globally near-optimal solution under our experimental protocol, and thus no convergence nor iteration statistics are reported.

\section{Discussion}
\label{sec:discussion}
In this work, we investigated the advantages and limitations of employing the \textit{deterministic global} solver MAiNGO for optimizing the acquisition function in Bayesian optimization (BO), comparing its performance to that of the more conventional stochastic solvers used: \textit{informed} multi-start (IMS), i.e., a \textit{stochastic global} solver, and \textit{informed} \textit{local} solver (ILS).
For our investigation, we used the M\"{u}ller-Brown potential function as our baseline case study with the lower confidence bound \eqref{eq:lcb} acquisition function under varying $\kappa$.

The most important consideration in BO is arguably the probability of convergence to a \textit{global}ly near-optimal solution of the black-box function.
We analyzed this critereon using \eqref{eq:lcb} and varying initial datasets for optimizing the M\"{u}ller-Brown potential.
Our findings reveal that, in \textit{exploitative} settings, BO using MAiNGO (compared to BO using ILS or IMS) is more prone to terminating at a \textit{local} minimum of the black-box function when the acquisition function is insufficiently exploratory and especially when the initial dataset lacks diversity.
In contrast, BO using IMS and ILS have a higher probability of identifying a globally near-optimal solution under these conditions than BO using MAiNGO.
Increasing $\kappa$ (the exploration parameter of \eqref{eq:lcb}) or the size and diversity of the initial dataset boosts convergence to a globally near-optimal solution across all solvers.
We also demonstrated the utility of statistical methods (e.g., conditional logistic regression, minimax testing, and paired t-tests) for quantifying the trade-offs between selecting a high-quality initial dataset and choosing an appropriate solver while accounting for the stochasticity introduced by the solver initialization.

The second most important consideration is the number of iterations needed to reach the termination criterion \eqref{eq:termination_critereon}, i.e., the number of experiments to be performed, especially when experiments are costly.
Our findings reveal that, in \textit{exploitative} settings and given that BO converges to a globally near-optimal solution, BO using MAiNGO generally converges in fewer iterations than BO using ILS or IMS.
We found that when the acquisition function is more \textit{exploratory}, BO using MAiNGO often requires a similar, if not greater, number of iterations to satisfy~\eqref{eq:termination_critereon} than BO with ILS or IMS.

We extended our analysis (of the probability of converging to a globally near-optimal solution and the number of iterations to reach the termination criterion~\eqref{eq:termination_critereon}) to six additional benchmark functions, suggesting that these trends are generalizable.
Additionally, our results indicate that MAiNGO may be rather advantageous in the BO setting with \textit{increasing dimensionality of the black-box function}.
Arguably, this apparent advantage in higher-dimensional problems is the most suggestive finding.

Beyond the aforementioned considerations, the computational cost per BO iteration is also a key consideration in real-world applications.
In our investigation, MAiNGO generally required more time per BO iteration than IMS or ILS, particularly when the number of optimization variables was large or the problem relaxations were weak (resulting in weak lower bounds and, as such, very slow convergence).
Some BO iterations in MAiNGO took seconds, while other iterations exceeded the prespecified time limit for MAiNGO, necessitating an adjustment of optimality tolerances.
As a test, we removed the prespecified time limit in MAiNGO.
Removing this time limit led to some iterations taking as much as \SI{24}{\hour}.
This highlights a trade-off: using deterministic global solvers in the BO context can reduce the number of iterations needed for convergence but may require more computation time per iteration.
Notably, the benefits (fewer iterations) and drawbacks (longer per-iteration times) scale with the dimensionality of the black-box function.
In applications involving costly experiments, allowing a single MAiNGO iteration to run overnight (or even substantially longer) is warranted, given its advantage of minimizing the total number of BO iterations compared to BO with ILS or IMS.
Significant progress has already been made in reducing computation times for optimization  with Gaussian Processes (GP) embedded through advancements in reduced-space formulations (as opposed to conventional full-spaced formulations)~\citep{Schweidtmann_2021} and specialized Branch-and-Bound algorithms~\citep{Tang_2024}.
These developments have made deterministic global optimization for BO increasingly computationally feasible, positioning this work among the first to explore its potential in BO.
However, further reductions in MAiNGO's computation time are still needed to fully realize this potential (see~\cref{sec:future_work}).

\section{Recommendations}
\label{sec:recommendations}
Based on the seven functions considered, we discuss some practical considerations.
Note that our following recommendations depend on the perspective of the decision maker.
A practitioner seeking a good solution may prioritize efficiency over rigor and be satisfied with finding a local minimum, while another may instead value rigor and finding a globally near-optimal solution.
For those seeking rigor, it is important to note that BO inherently lacks complete precision due to potential non-representative initial data and suboptimal acquisition function hyperparameters.
Considering these perspectives and the findings of \cref{sec:discussion}, we provide the following heuristic recommendations:
\begin{itemize}
    \item \textbf{Prioritize initial dataset quality, especially when using deterministic global solvers in BO.}
    A poorly representative initial dataset can produce a misleading initial surrogate model, reducing the effectiveness of BO in finding a globally near-optimal solution of the black-box function.
    While this issue affects all solvers, it is particularly problematic for deterministic global solvers like MAiNGO when paired with insufficiently exploratory acquisition functions.
    Our statistical analysis on our baseline case study indicates that initial dataset selection \textit{is at least as} important as solver choice with regard to the probability of convergence to a globally near-optimal solution.
    A poor initial dataset may also prolong convergence.

    \item \textbf{Consider leveraging deterministic global solvers like MAiNGO when well-tuned acquisition functions are available; otherwise, if tuning is not possible, favor stochastic solvers.}
    This recommendation is also in accordance with our frequent observation that if we optimize correctly, i.e., by using a global deterministic solver, and the results are not satisfactory, it indicates that the objective (acquisition function) may be flawed.
    BO using MAiNGO excels when the acquisition function is well-calibrated for the task (e.g., exploitative or exploratory) and when the initial surrogate GP is representative of the black-box function.
    BO with MAiNGO reliably follows the behavior dictated by the acquisition function.
    For example, if the decision maker envisions highly exploitative behavior (inscribed into the \eqref{eq:lcb} acquisition function with a low value of $\kappa$), BO using MAiNGO is indeed very exploitative.
    This is not the case for BO using ILS or IMS, as these solvers might only report a local solution of the acquisition function, deviating from the intended (highly exploitative) behavior.
    This stochastic damping, caused by ILS or IMS, reduces the direct impact of acquisition function hyperparameters, which can be beneficial when the initial dataset is suboptimal or the acquisition function is insufficiently exploratory.
    Stochasticity can increase the likelihood of finding a globally near-optimal solution by facilitating the exploration of previously uncharted regions: \citep{bull_2011,Ath_2021} introduce further stochasticity by incorporating $\epsilon$-greedy strategies into the BO framework.
    Parallels exist in machine learning, where stochastic gradient descent solvers are used to ``escape'' local minima \citep{Markidis_2021} to further explore previously uncharted regions.

    \item \textbf{If the goal is to find \textit{any} local minimum of the black-box function, consider pairing an exploitative acquisition function with a deterministic global solver, even in the absence of well-tuned (hyperparameters of) acquisition functions or diverse initial datasets.}
    The strict adherence of deterministic global solvers like MAiNGO to the acquisition function enables rapid exploitation of known regions, minimizing the total number of iterations when paired with exploitative acquisition functions.
    This makes it an efficient choice when locating a globally near-optimal solution is not necessary.

    \item \textbf{Use exploratory acquisition functions with deterministic global solvers judiciously.}
    In exploratory settings, BO with MAiNGO achieves comparable or better probabilities of convergence to a globally near-optimal solution but often requires more iterations compared to ILS or IMS.
    We postulate this is because MAiNGO strictly follows the acquisition function, while IMS and ILS may get trapped in local minima of the acquisition function, unintentionally exhibiting (different) exploratory behavior.

    \item \textbf{Utilize hybrid approaches.}
    A practical strategy may be to begin with a stochastic solver or exploratory acquisition function when the representativeness of the initial dataset is uncertain.
    As the dataset becomes sufficiently rich (i.e., when the associated GP encompasses enough of the the black-box function), transitioning to an exploitative function with MAiNGO can capitalize on its efficiency in converging to a globally near-optimal solution.
    This approach contrasts with theoretical approaches in~\citep{Srinivas_2012,Kandasamy_2015}, which advocate increased exploration over time.
    We hypothesize that the proposed approach could perform well in practice and leave its rigorous validation to future research.
    Our hypothesis is supported by the successful application of similar hybrid approaches in machine learning \citep{Markidis_2021}, where in the training phase, first, a stochastic gradient descent solver is used to efficiently navigate the search space and obtain a globally near-optimal solution, and subsequently a local solver is employed to refine this solution.
\end{itemize}
Many of these recommendations require knowledge of what constitutes an exploratory versus exploitative acquisition function (and how to balance them in the absence of knowledge) and when an initial dataset is representative enough.
These factors are challenging to determine \textit{a priori} and often require informed judgment.
Several works are actively addressing these challenges.
\citep{bull_2011,Ath_2021} balance exploitation and exploration by incorporating $\epsilon$-greedy strategies into the BO framework; \citep{Yan_2022,Chen_2023,Benjamins_2023,Maitra_2024} work on ``self-adaptive'' acquisition functions; and \citep{bull_2011,Ren_2024,Wang_2024} provide insights on selecting the initial dataset.

\section{Future Work}
\label{sec:future_work}
Our investigation showed that deterministic global solvers are effective in BO when the acquisition function is well-tuned and the initial dataset is diverse.
Consequently, a promising direction for future research is the development of acquisition functions that can be dynamically tuned on-the-fly or tuned without requiring prior knowledge ~\citep{Yan_2022,Chen_2023,Benjamins_2023,Maitra_2024}.
Additionally, determining whether an initial dataset is sufficiently diverse \textit{before} beginning BO and devising sampling strategies to produce such diverse datasets (e.g.,~\citep{bull_2011,Ren_2024,Wang_2024}) remain critical challenges: these challenges require a deeper understanding of what constitutes a ``diverse-enough'' dataset for initial reliable surrogate modeling.
Advancing these areas can increase the probability of BO converging to a globally near-optimal solution, broadening the scenarios where deterministic global solvers are effective.
The complementary strategy of starting with a stochastic solver and transitioning to MAiNGO also warrants further investigation.

In exploitative settings with a well-tuned acquisition function and diverse initial dataset, BO with MAiNGO can converge in significantly fewer iterations than stochastic solvers.
However, we observed that computation times with MAiNGO could be considerable, necessitating setting time limits.
Specifically, while the upper bound in MAiNGO converges relatively quickly, the lower bound does not, c.f., \cref{sec:maingo}.
This observation aligns with findings of \cite{Tang_2024}, who have shown that McCormick relaxations~\citep{McCormick_1976}, which are also used in MAiNGO for generating a lower bound, and also $\alpha$BB underestimators~\citep{Adjiman_1996} may be weak when optimizing the posterior mean function for GPs (due to a similar structure of \eqref{eq:lcb}, its relaxation may also be weak).
To address this limitation, \cite{Tang_2024} proposed a spatial Branch-and-Bound algorithm for optimizing GP posterior mean functions, achieving significant speed-ups.
Future work should extend this approach to acquisition functions or develop stronger lower-bounding heuristics to optimize acquisition functions more efficiently.
Improving lower-bound convergence would enhance the computational efficiency of MAiNGO, making it even more competitive for BO when appropriate.

To further validate deterministic global solvers in BO, more extensive benchmark testing and real-world experiments are required, particularly under noisy conditions, which were not considered in this study.
Additionally, exploring deterministic global optimization for GP hyperparameter optimization (e.g., during training), instead of the local solver currently used, could enhance BO performance.
These research directions will collectively clarify the scenarios in which deterministic global solvers are most effective, advancing BO’s capability to solve complex optimization problems in science and engineering.

\section{Statements and Declarations}
\label{sec:statements_and_declerations}
\subsection{Competing Interests}
\label{sec:competing_interests}
The authors have no competing interests to declare that are relevant to the content of this article.

\subsection{Data Availability Statement}
The authors declare that the data supporting the findings of this study are available within the paper.
Should any raw data files be needed in another format they are available from the corresponding author upon reasonable request.

\subsection{Funding}
\label{sec:funding}
This research was supported by the National Cancer Institute (NCI) under grant number T32CA153952, the National Science Foundation under Grant Number 143695, and the US Department of Energy.
It has also been funded by the Deutsche Forschungsgemeinschaft (DFG, German Research Foundation) under Germany's Excellence Strategy – Cluster of Excellence 2186 \glqq{}The Fuel Science Center\grqq{} – ID: 390919832.
We also acknowledge support of the Werner Siemens Foundation in the frame of the WSS Research Centre ``catalaix''.

\subsection{Acknowledgments}
\label{sec:acknowledgments}
We thank Clara Witte and Jannik Lüthje for their invaluable discussions and guidance on MAiNGO and MeLOn.

\subsection{Author Contributions}
\label{sec:author_contributions}
All authors contributed to the study's conception and/or design.
Methodology development, software implementation, execution of experiments, and formal analysis were led by Anastasia Georgiou, supported by Daniel Jungen.
Anastasia Georgiou and Daniel Jungen drafted the initial manuscript, with all authors contributing to writing and editing.
Verena Hunstig validated the numerical results and software functionality.
Constantine Frangakis developed the statistical framework and conducted conditional logistic regression and minimax testing.
Luise Kaven contributed to the conceptual framework.
Ioannis Kevrekidis and Alexander Mitsos provided supervision, conceptual guidance, and secured funding.

%%%%%%%%%%%%%%%%%%%%%%%%%%%%%%%%%%%%%%%%%%%%%%%%%%%%%%%%%%%%%%%
% bibliography
%%%%%%%%%%%%%%%%%%%%%%%%%%%%%%%%%%%%%%%%%%%%%%%%%%%%%%%%%%%%%%%
\bibliography{bibliography}
\bibliographystyle{abbrvnat}

%%%%%%%%%%%%%%%%%%%%%%%%%%%%%%%%%%%%%%%%%%%%%%%%%%%%%%%%%%%%%%%
% appendix/SI
%%%%%%%%%%%%%%%%%%%%%%%%%%%%%%%%%%%%%%%%%%%%%%%%%%%%%%%%%%%%%%%
\clearpage
\appendix

\renewcommand{\thepage}{S\arabic{page}}
\renewcommand{\thesection}{S\arabic{section}}
\renewcommand{\thetable}{S\arabic{table}}
\renewcommand{\thefigure}{S\arabic{figure}}
\renewcommand{\theequation}{S\arabic{equation}}

\setcounter{page}{1}
\setcounter{figure}{0}
\setcounter{equation}{0}
\setcounter{table}{0}

%%%%%%%%%%%%%%%%%%%%%%%%%%%%%%%%%%%%%%%%%%%%%%%%%%%%%%%%%%%%%%%
\section*{Supporting Information}
%%%%%%%%%%%%%%%%%%%%%%%%%%%%%%%%%%%%%%%%%%%%%%%%%%%%%%%%%%%%%%%
%
\section{Notation}
\label{sec:notation}
\begin{longtable}{p{0.21\linewidth}p{0.72\linewidth}}
    \caption{Summary of acronyms, notation, and key terms used.} \label{tab:notation_summary} \\
    \toprule
    \textbf{Symbol / Term} & \textbf{Definition}  \\
    \midrule
    \endfirsthead
    \caption[]{continued.}\\
    \toprule
    \textbf{Symbol / Term} & \textbf{Definition}  \\
    \midrule
    \endhead
    \multicolumn{2}{l}{\textbf{Acronyms}} \\
    BO & Bayesian Optimization \\
    cMLE & conditional Maximum Likelihood Estimation \\
    GP & Gaussian Process \\
    ILS & Informed Local Search \\
    IMS & Informed Multi Start \\
    LCB & Lower Confidence Bound acquisition function \\
    MAiNGO & Mixed-Integer Nonlinear Global Optimizer \\
    MeLOn & Machine Learning models for Optimization \\
    MB & Müller-Brown potential \\
    PK-MIQP & Piecewise-linear Kernel Mixed Integer Quadratic Programming \\
    TC & Termination criteria \\
    \\[-0.75em] \midrule
    \multicolumn{2}{l}{\textbf{Bayesian Optimization Notation}} \\
    $\mathcal{GP}$ & Gaussian process \\
    $\mathcal{D}$ & initial dataset \\
    $D$ & Dimensionality of the search space or optimization problem \\
    $N$ & Number of initial data points \\
    $U$ & Müller-Brown potential function \\
    $\boldsymbol{x}\in\mathbb{R}^D$ & Decision variable in $D$ dimensions \\
    $\boldsymbol{x}^t$ & Candidate at iteration $t$ \\
    $\boldsymbol{x}^{*}$ & Global minimizer of $f$ \\
    $\boldsymbol{\hat{x}}^{*}$ & Global minimizer of $f$ \\
    $\boldsymbol{x}^{\mathrm{init}}$ & Matrix of initial candidate points\\
    $f$ & Black-box objective function \\
    $f^* = f(\boldsymbol{x}^{*})$ & Globally optimal objective function value \\
    $f^t = f(\boldsymbol{x}^t)$ & Objective value for current iterate $\boldsymbol{x}^t$ \\
    $f^{*,t-1}$ & Best objective value observed up to iteration $t-1$ \\
    $LCB$ & Lower confidence bound acquisition function \\
    $\nu$ & Matérn kernel parameter \\
    $\mu$ & GP posterior mean \\
    $\sigma$ & GP posterior standard deviation \\
    $k$ & Covariance \\
    $\kappa$ & Exploration–exploitation hyperparameter for LCB \\
    $\kappa_{t,K}$ & Hyperparameter at iteration $t$, scheduled per~\cite{Kandasamy_2015} \\
    $\kappa_{t,S}$ & Hyperparameter at iteration $t$, scheduled per~\cite{Srinivas_2012} \\
    $M$ & Finite discretization size used in $\kappa_{t,S}$ \\
    $\delta$ & Confidence parameter used in $\kappa_{t,S}$ \\
    $\varepsilon$ & Near-optimality tolerance \\
    $\varepsilon_{x,i}$ & Absolute movement tolerance in $x$ ($\varepsilon_{x,1}<\varepsilon_{x,2}$) \\
    $\varepsilon_{f,r}$ & Relative objective-value tolerance \\
    $\varepsilon_{f,a}$ & Absolute objective-value tolerance \\
    \\[-0.75em] \midrule
    \multicolumn{2}{l}{\textbf{Experimental Design Terminology}} \\
    case study & Configuration defined by a specific objective function, acquisition function, and stopping criterion. All experiments and runs under it share these settings. \\
    experiment & A single instantiation of a case study with its own initial dataset (i.e., a unique initial design). \\
    run & One full BO execution performed for a given experiment. \\
    \\[-0.75em] \midrule
    \multicolumn{2}{l}{\textbf{Statistics \& Experimental Indexing}} \\
    $dat$ & Index for a unique initial dataset \\
    $r$ & Run index (repeated executions of the full BO loop) \\
    $s$ & Solver index (e.g., ILS, IMS, MAiNGO) \\
    $ds$ & A possible subset of all experiments (all initial datasets) \\
    $ds^*$ & The subset of experiments that maximizes the worst possible convergence over solvers \\
    $y^{s}_{dat,r}$ & Indicator variable for whether solver \textit{s} converges to globally near-optimal solution of $f$ in run $r$ of experiment with initial dataset $dat$ \\
    $p^{s}_{dat}$ & Probability solver $s$ converges to a globally near-optimal solution on experiment \textit{dat} ($:=\mathbb{E}_{\text{over } r}\{y^{\text{ ILS}}_{dat,r}\mid dat\}$) \\
    $p^{s}$ & Probability solver $s$ converges to a globally near-optimal solution averaged over all experiments $(:=\mathbb{E}_{\text{over } dat}\{y^{\text{ s}}_{dat}\})$ \\
    $\hat{p}^{s}_{ds}$ & Empirical proportion of convergence over the runs and datasets in $ds$ for solver $s$ \\
    $\alpha$ & common log-odds ratio of ILS vs MAiNGO \\
    $\beta$ & common log-odds ratio of IMS vs MAiNGO \\
    \bottomrule
\end{longtable}

\section{BO Frameworks}
\label{sec:bo_frameworks}
\begin{table*}[htbp]
    \caption{
    A non-exhaustive summary of recent Python-based BO framework software.
    BO packages exist in further programming languages, including C++ (e.g., MOE), R (e.g., DiceOptim, \& laCP), and MATLAB (e.g., DACE).
    }
    \label{tab:bo_methods}
    \begin{tabular}{p{0.12\linewidth} p{0.1\linewidth} p{0.16\linewidth} p{0.16\linewidth} >{\raggedright\arraybackslash}p{0.33\linewidth} }
        Software & Updated & GP Implementation & Acquisition Functions & Optimizer \\ [0.5ex] \hline
        \\
        \href{https://botorch.org/}{BOTorch} & 2024 & Wrapper around GPYTorch & PI, EI, UCB, Posterior Mean, MC-Based, + & Modular (scipy.optimize, torch.optim), with L-BFGS-B informed multi-start as default \\
        \\
        \href{https://github.com/secondmind-labs/trieste}{Trieste} & 2024 & GPFlow & PI, EI, LCB, Feasibility Based, MC-Based, + & Modular (tensorflow.optimizers, tf.keras.optimizers, scipy.optimize) with L-BFGS-B as default \\
        \\
        \href{https://github.com/bayesian-optimization/BayesianOptimization}{bayes-optimization} & 2024 & sklearn.GPR & PI, EI, UCB & scipy.optimize, L-BFGS-B informed multi-start \\
        \\
        \href{https://github.com/automl/SMAC3}{SMAC3} & 2024 & Wrapper around sklearn.GPR & PI, EI, LCB, TS & scipy.optimize, custom local \& random search~\cite{Lindauer_2022} \\
        \\
        \href{https://pygpgo.readthedocs.io/en/latest/}{pyGPGO} & 2022 & Native & PI, EI, UCB, Entropy & scipy.optimize, L-BFGS-B (random) multi-start \\
        \\
        \href{https://scikit-optimize.github.io/stable/modules/generated/skopt.gp_minimize.html}{skopt} & 2021 & sklearn.GPR & EI, LCB, PI & L-BFGS-B informed multi-start, parallel \\
        \\
        \href{https://gpyopt.readthedocs.io/en/latest/GPyOpt.methods.html}{GPyOpt} & 2020 (archived) & \href{https://github.com/SheffieldML/GPy/blob/deploy/GPy/models/gp_regression.py}{GPy} & EI, LCB, PI, Entropy & scipy.optimize, L-BFGS-B multi-start, DIRECT, covariance matrix adaptation
    \end{tabular}
\end{table*}

\section{Additional Test Cases}
\label{sec:test-case-formulation}
\subsection{Six-Hump Camelback}
\label{sec:six_hump_camelback}
2D benchmark~\citep{Jamil_2013}.
Six minima, with two global minima located at $(\pm \num{0.0898}, \mp \num{0.7126})$, with objective value $f\left(\boldsymbol{x}^{*}\right) \approx \num{-1.0316}$.
\begin{equation}
    \label{eq:six_humb_camelback}
    \begin{array}{rlrlrr}
        f\left(\boldsymbol{x}^{*}\right) \mpeq \mpmin{\boldsymbol{x}} & \left(4 - 2.1 x_{1}^{2} + \frac{x_{1}^{4}}{3}\right) x_{1}^{2} + x_{1} x_{2} + \left(-4 + 4x_{2}^{2}\right) x_{2}^{2} \\
        \mpst   & x_{1} \in [-3, 3], \ x_{2} \in [-2, 2]
    \end{array}
\end{equation}

\subsection{Ackley}
\label{sec:3d_ackley}
3D benchmark~\citep{Mitchell_1987}.
Many local minima; one global minimum at $\left(0, 0, 0\right)^{T}$, with objective value $f\left(\boldsymbol{x}^{*}\right) = 0$.
\begin{equation}
    \label{eq:3d_ackley}
    \begin{array}{rlrlrr}
        f\left(\boldsymbol{x}^{*}\right) \mpeq \mpmin{\boldsymbol{x}} & -20 \exp\left(-0.2 \sqrt{\frac{1}{3} \sum_{i=1}^{3} x_{i}^{2}}\right)
                    - \exp\left(\frac{1}{3} \sum_{i=1}^{3} \cos(2\pi x_{i})\right) + 20 + e \\
        \mpst   & \boldsymbol{x} \in \left[-5, 5\right]^{3}
    \end{array}
\end{equation}

\subsection{Hartmann}
\label{sec:4d_hartmann}
4D benchmark~\citep{Dixon_1978}.
One global minimum located at $\left(0.1146, 0.5556, 0.8525, 0.8525\right)^{T}$, with objective value $f\left(\boldsymbol{x}^{*}\right) \approx -3.3224$.
\begin{equation}
    \label{eq:4d_hartmann}
    \begin{array}{rlrlrr}
        f\left(\boldsymbol{x}^{*}\right) \mpeq \mpmin{\boldsymbol{x}} & -\sum_{i=1}^{4} \alpha_{i} \exp\left(-\sum_{j=1}^{4} A_{ij} (x_{j} - P_{ij})^2\right) \\
        \mpst   &  \boldsymbol{x} \in \left[0, 1\right]^{4},
    \end{array}
\end{equation}
with
\begin{equation}
    \nonumber
    \begin{array}{c}
    \boldsymbol{\alpha}  = \left(1.0, 1.2, 3.0, 3.2\right)^{T}, \\ \\
    \boldsymbol{A} =
                    \begin{bmatrix}
                        10      & 3     & 17    & 3.5 \\
                        0.05    & 10    & 17    & 0.1 \\
                        3       & 3.5   & 1.7   & 10 \\
                        17      & 8     & 0.05  & 10
                    \end{bmatrix}, \\ \\
    \boldsymbol{P} =
                    \begin{bmatrix}
                        0.1312 & 0.1696 & 0.5569 & 0.0124 \\
                        0.2329 & 0.4135 & 0.8307 & 0.3736 \\
                        0.2348 & 0.1451 & 0.3522 & 0.2883 \\
                        0.4047 & 0.8828 & 0.8732 & 0.5743
                    \end{bmatrix}.
    \end{array}
\end{equation}

\subsection{GKLS Test Functions}
\label{sec:gkls}
2D, 3D, and 4D benchmarks~\citep{Gaviano_2003}.
The GKLS generator creates test functions by smoothly perturbing a quadratic bowl to produce controlled local minima and a known global minimum.
For each selected dimension $D \in \{2,3,4\}$ we generate one GKLS instance using the parameters listed in each subsection.
All functions share a global minimum value of $f^{*} = -1.0$.
Tthe search domain is the unit hypercube $\boldsymbol{x} \in [-1,1]^{D}$.
Parameters follow the widely used GKLS classes reported in \citep{Kudela_2023}, with a larger radius ($r$) of the attraction region of the global minimizer to make the problems more amenable to BO.

\noindent
\subsubsection{2D GKLS}
2D.
Twice differentiable.
Parameters: distance from the global minimizer to the vertex of the quadratic function $d= 0.90$; radius of the attraction region of the global minimizer $r= 0.40$ (increased relative to defaults); number of local minima $h = 10$, random seed $s=12$.
Many local minima; one global minimum at $(-0.7270, 0.1728)^T$ with objective value $f\left(\boldsymbol{x}^{*}\right) = -1.0$.

\subsubsection{3D GKLS}
3D.
Twice differentiable.
Parameters: distance from the global minimizer to the vertex of the quadratic function $d= 0.66$; radius of the attraction region of the global minimizer $r= 0.30$ (increased relative to defaults); number of local minima $h = 10$, random seed $s=12$.
Many local minima; one global minimum at $(-0.4713, 0.0908, 0.3609)^T$ with objective value $f\left(\boldsymbol{x}^{*}\right) = -1.0$.

\subsubsection{4D GKLS}
4D.
Twice differentiable.
Parameters: distance from the global minimizer to the vertex of the quadratic function $d= 0.66$; radius of the attraction region of the global minimizer $r= 0.20$; number of local minima $h = 10$, random seed $s=12$.
Many local minima; one global minimum at $(-0.4278, 0.0678, 0.2233, -0.7765)^T$ with objective value $f\left(\boldsymbol{x}^{*}\right) = -1.0$.

\clearpage  % formatting
\section{Set Up Details}
\label{sec:set-up-details}
\begin{table}[h]
    \centering
    \caption{
        Settings for all case studies, detailing the number of points in the initial dataset ($N$), the $\kappa$ value used in~\eqref{eq:lcb}, the number of experiments and runs, and the termination criterion thresholds.
    }
    \label{tab:testcase-setup-table}
    \begin{tabular}{ c c c c c c c c c }
        Case Study & $N$ & $\kappa$  & Experiments & Runs per Exp. & $\varepsilon_{x,1}$ & $\varepsilon_{x,2}$ & $\varepsilon_{f,r}$ & $\varepsilon_{f,a}$ \\ [1ex]
        \hline
        \multirow{5}{*}{2D Müller-Brown}
        & 3  & 2   & 56 & 31 & \multirow{5}{*}{0.001} & \multirow{5}{*}{0.05} & \multirow{5}{*}{0.01} & \multirow{5}{*}{0.5} \\
        & 3  & 3   & 25 & 31 &  &  &  &  \\
        & 10 & 2   & 31 & 31 &  &  &  &  \\
        & 20 & 2   & 41 & 31 &  &  &  &  \\
        & 20 & 0.5 & 16 & 31 &  &  &  &  \\
        & 3 & \ref{eq:kappa_t_Srinivas_2012}  & 21 & 21 &  &  &  &  \\
        & 3 & \ref{eq:kappa_t_Kandasamy_2015} & 51 & 21 &  &  &  &  \\[1ex]
        \hline
        2D Camelback & 3 & 2 & 41 & 31 & 0.001 & 0.05 & 0.02 & 0.05 \\ [1ex]
        2D GKLS & 3 & 3 & 25 & 25 & 0.001 & 0.05 & 0.02 & 0.05 \\ [1ex]
        \hline
        3D Ackley & 4 & 2 & 31 & 15 & 0.001 & 0.05 & 0.02 & 0.05 \\ [1ex]
        3D GKLS & 4 & 3 & 25 & 25 & 0.001 & 0.05 & 0.01 & 0.02 \\ [1ex]
        \hline
        4D Hartmann & 5 & 2 & 30 & 16 & 0.001 & 0.02 & 0.02 & 0.01 \\ [1ex]
        4D GKLS & 5 & 3 & 25 & 25 & 0.001 & 0.02 & 0.005 & 0.01 \\ [1ex]
    \end{tabular}
\end{table}

\break

\section{Simple Regret}
\label{sec:regret-plots}
\subsection{Considering Only Runs That Reached a Globally Near-Optimal Solution}
\begin{figure}[htbp]
    \centering
    \begin{minipage}[t]{0.47\textwidth}
        \centering
        \includegraphics[trim={0.2cm 0.2cm 0.2cm 0cm}, clip, width=\textwidth]
        {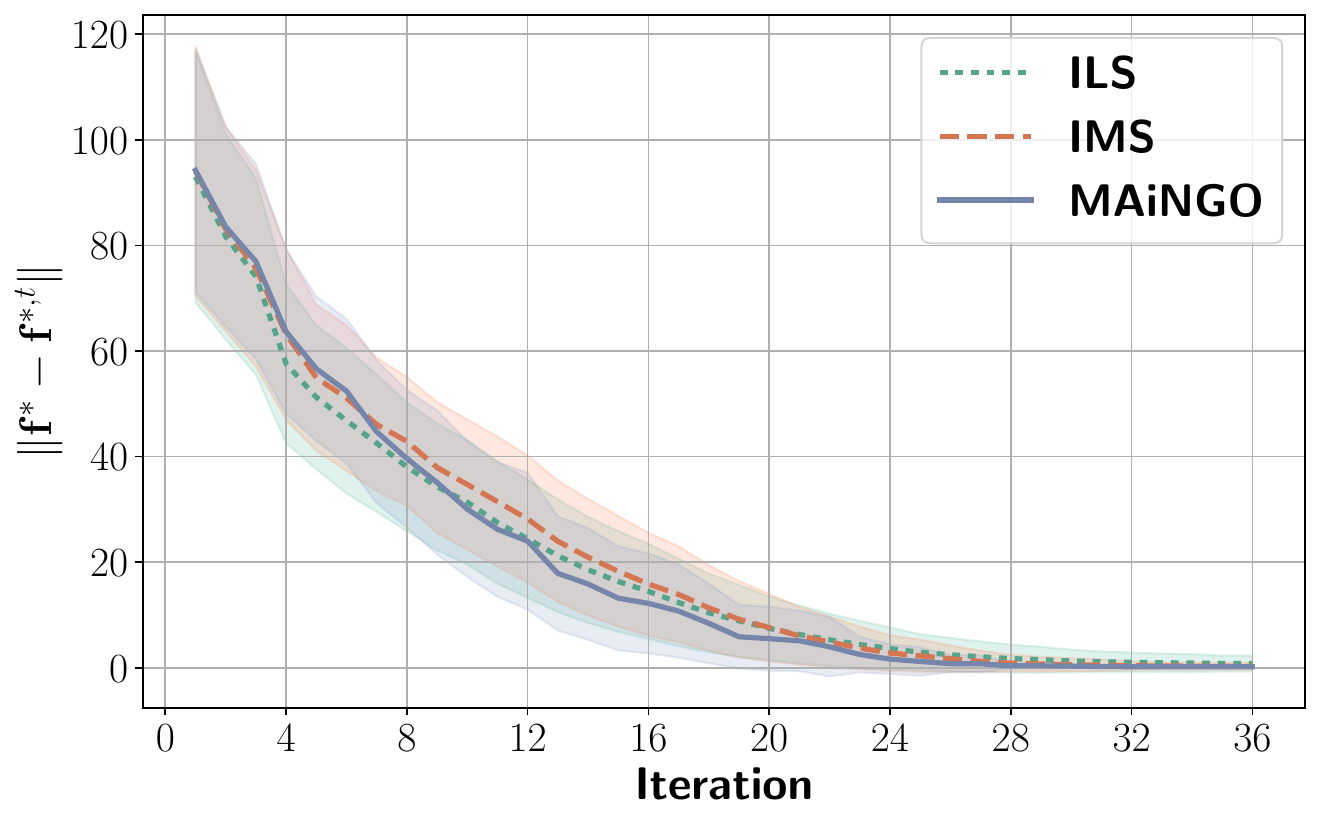}
        \subcaption{Müller-Brown, $\kappa = 2$ and $N = 3$.}
        \label{fig:MB-3t-2k-regret:k2_N2}
    \end{minipage}%
    \hfill%
    \begin{minipage}[t]{0.47\textwidth}
        \centering
        \includegraphics[trim={0.2cm 0.2cm 0.2cm 0cm}, clip, width=\textwidth]
        {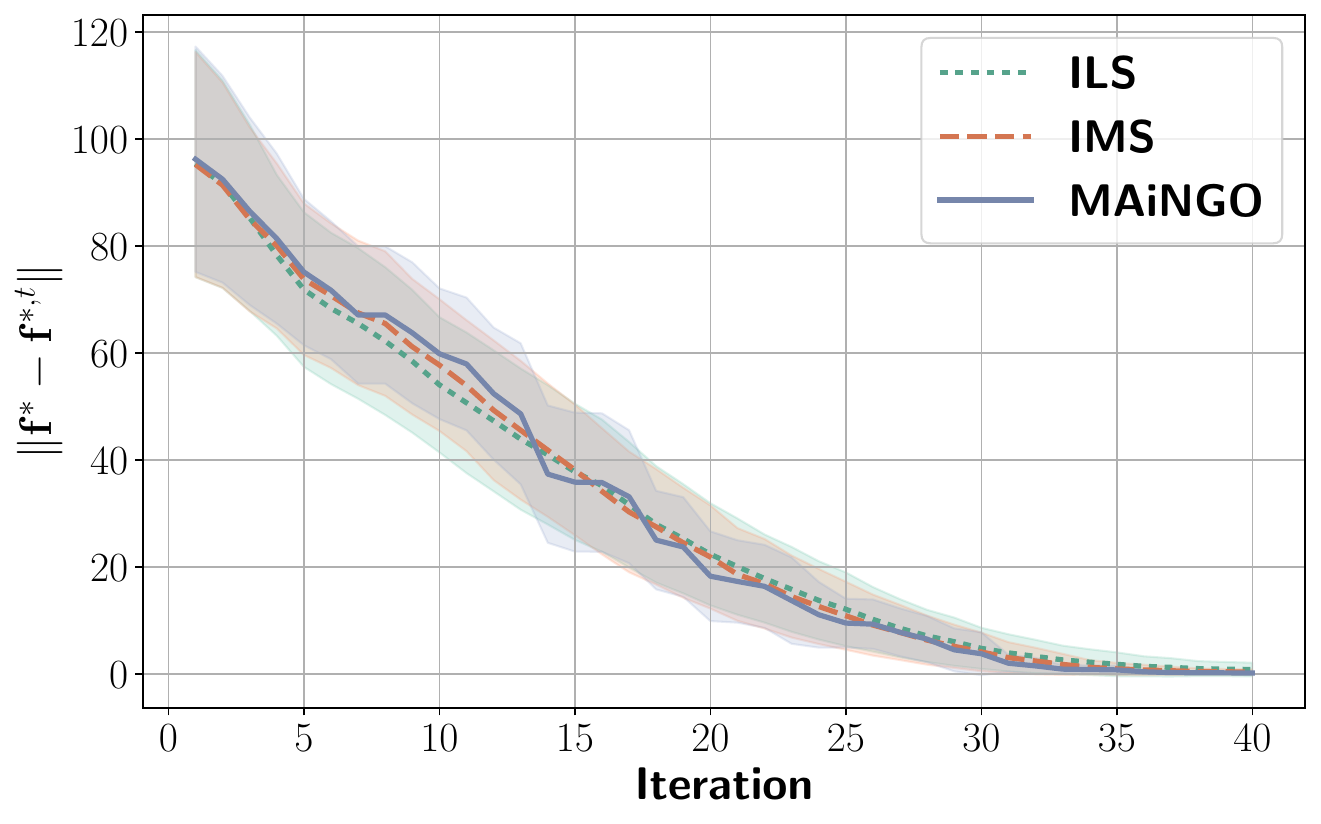}
        \subcaption{Müller-Brown, $\kappa = 3$ and $N = 3$.}
        \label{fig:MB-3t-2k-regret:k3_N3}
    \end{minipage}%
    \par\medskip
    \begin{minipage}[t]{0.47\textwidth}
        \centering
        \includegraphics[trim={0.2cm 0.2cm 0.2cm 0cm}, clip, width=\textwidth]
        {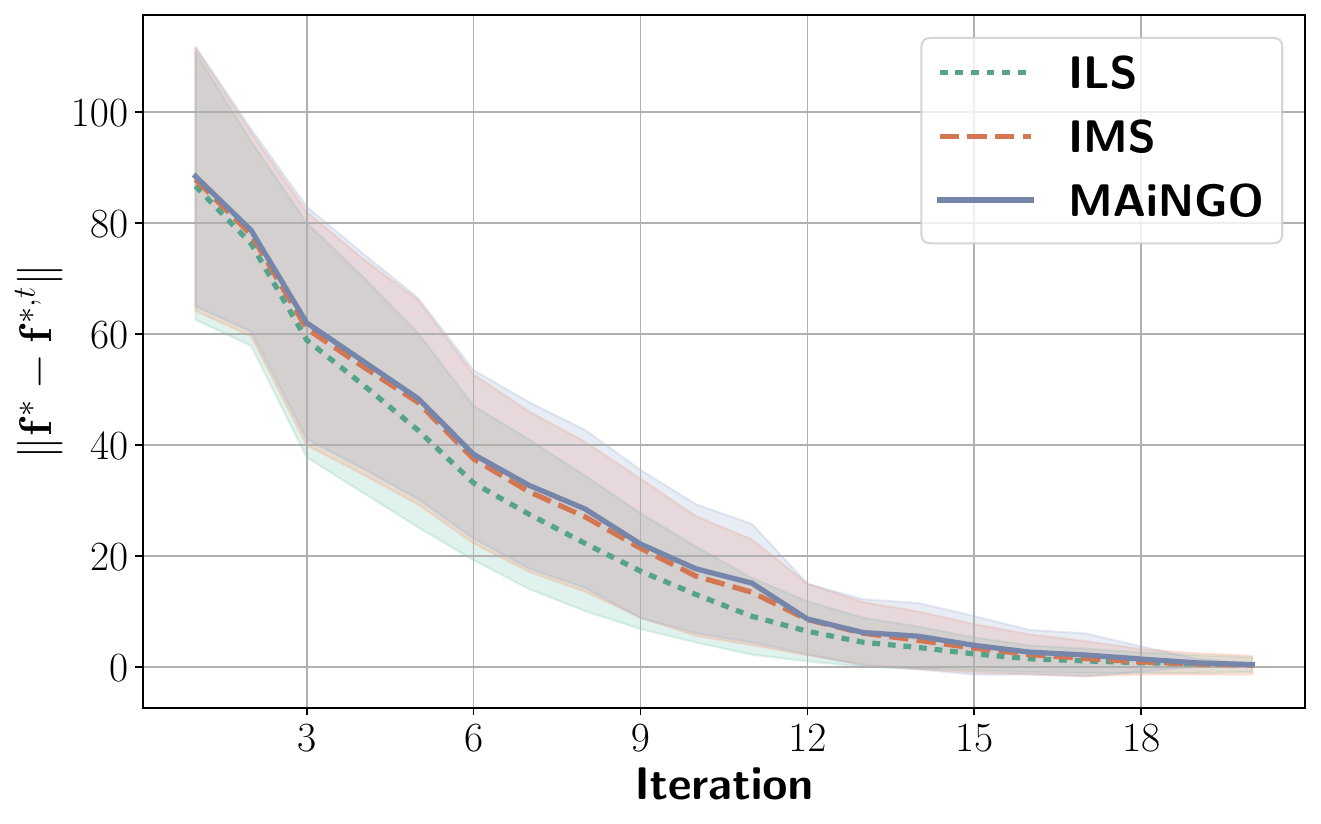}
        \subcaption{Müller-Brown, $\kappa =$ ~\ref{eq:kappa_t_Kandasamy_2015} and $N = 3$.}
        \label{fig:MB-3t-2k-regret:kK_N3}
    \end{minipage}%
    \hfill%
    \begin{minipage}[t]{0.47\textwidth}
        \centering
        \includegraphics[trim={0.2cm 0.2cm 0.2cm 0cm}, clip, width=\textwidth]
        {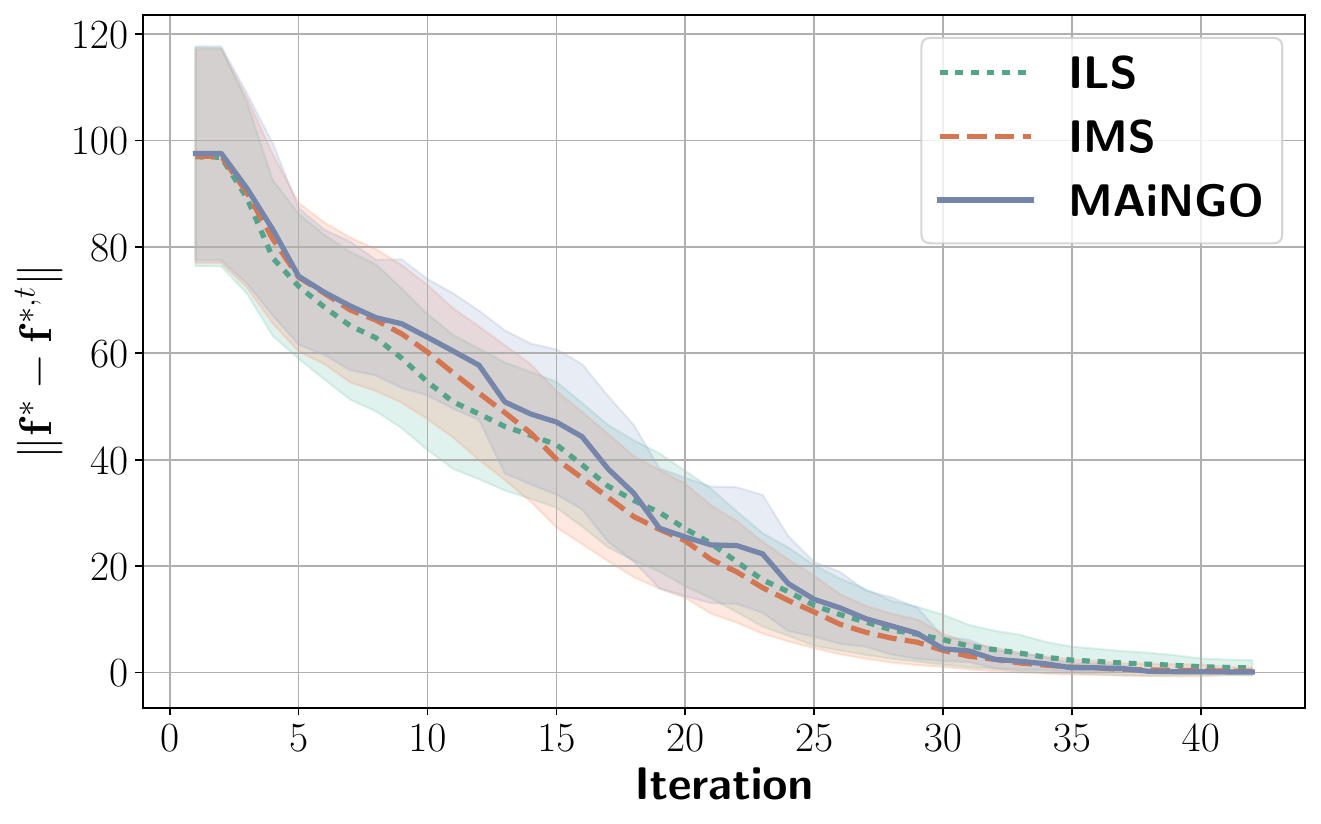}
        \subcaption{Müller-Brown, $\kappa =$ ~\ref{eq:kappa_t_Srinivas_2012} and $N = 3$.}
        \label{fig:MB-3t-2k-regret:kS_N3}
    \end{minipage}%
    \par\medskip
    \begin{minipage}[t]{0.47\textwidth}
        \centering
        \includegraphics[trim={0.2cm 0.2cm 0.2cm 0cm}, clip, width=\textwidth]
        {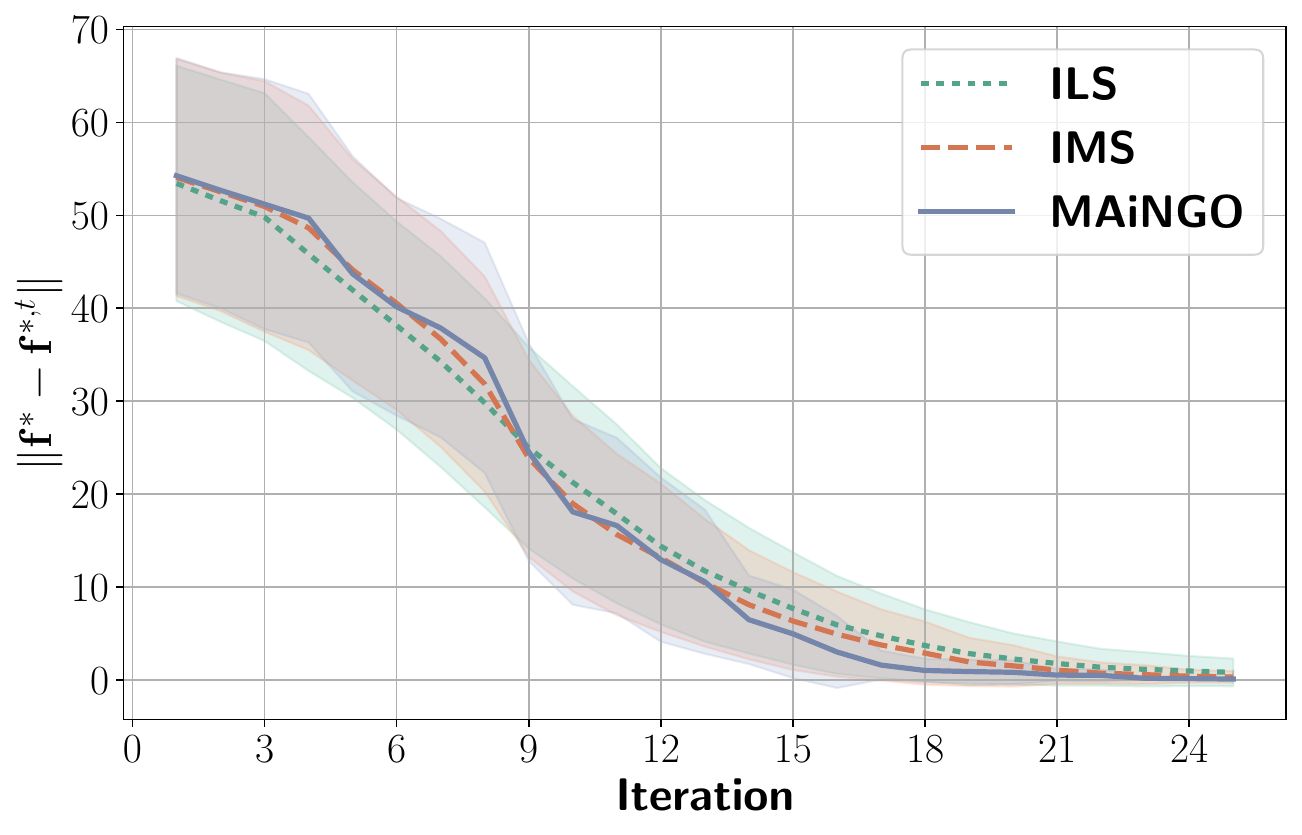}
        \subcaption{Müller-Brown, $\kappa = 2$ and $N = 10$.}
        \label{fig:MB-3t-2k-regret:k2_N10}
    \end{minipage}%
    \hfill%
    \begin{minipage}[t]{0.47\textwidth}
        \centering
        \includegraphics[trim={0.2cm 0.2cm 0.2cm 0cm}, clip, width=\textwidth]
        {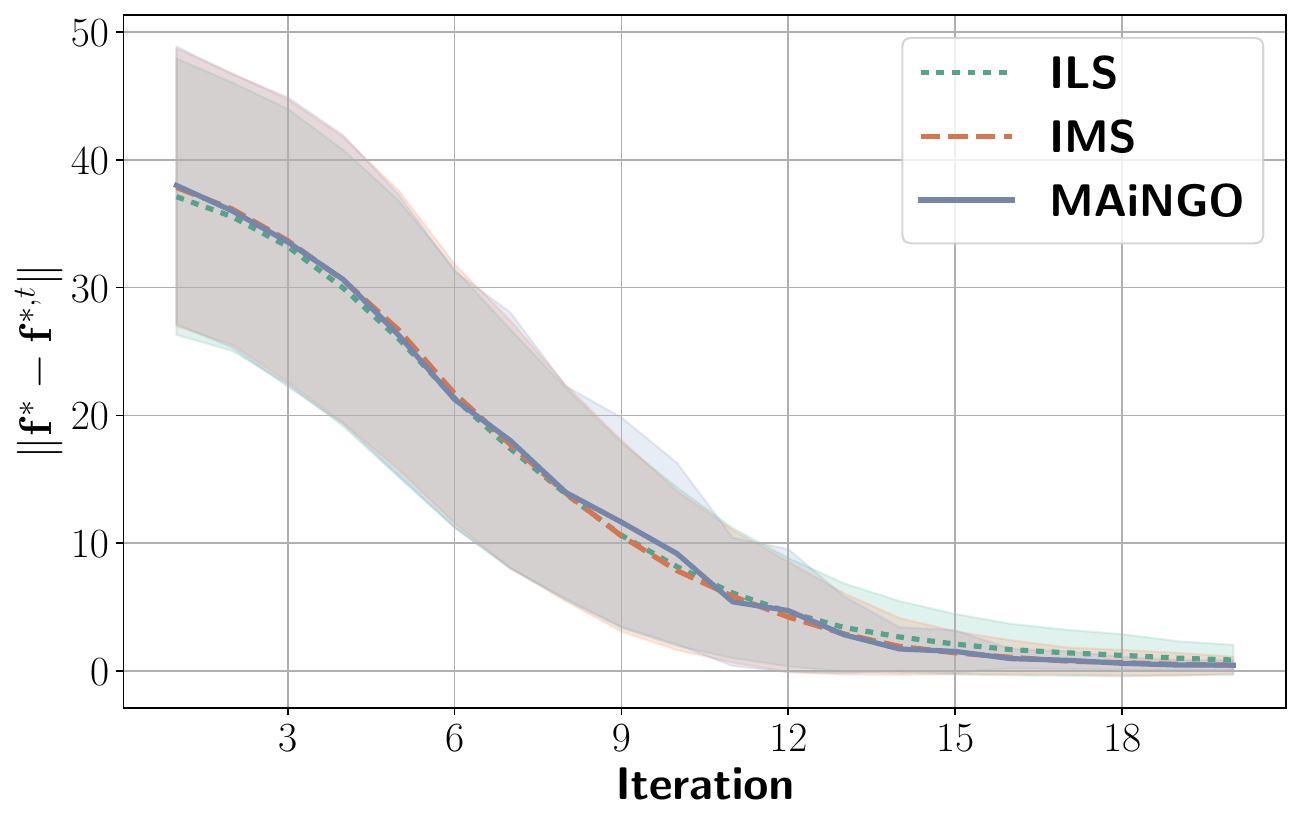}
        \subcaption{Müller-Brown, $\kappa = 2$ and $N = 20$.}
        \label{fig:MB-3t-2k-regret:k2_N20}
    \end{minipage}%
    \caption{
      Simple regret plots (Part 1).
    }
    \label{fig:combined-regret-plots:part:1}
\end{figure}

\begin{figure}[htbp]\ContinuedFloat
    \centering
    \begin{minipage}[t]{0.47\textwidth}
        \centering
        \includegraphics[trim={0.2cm 0.2cm 0.2cm 0cm}, clip, width=\textwidth]
        {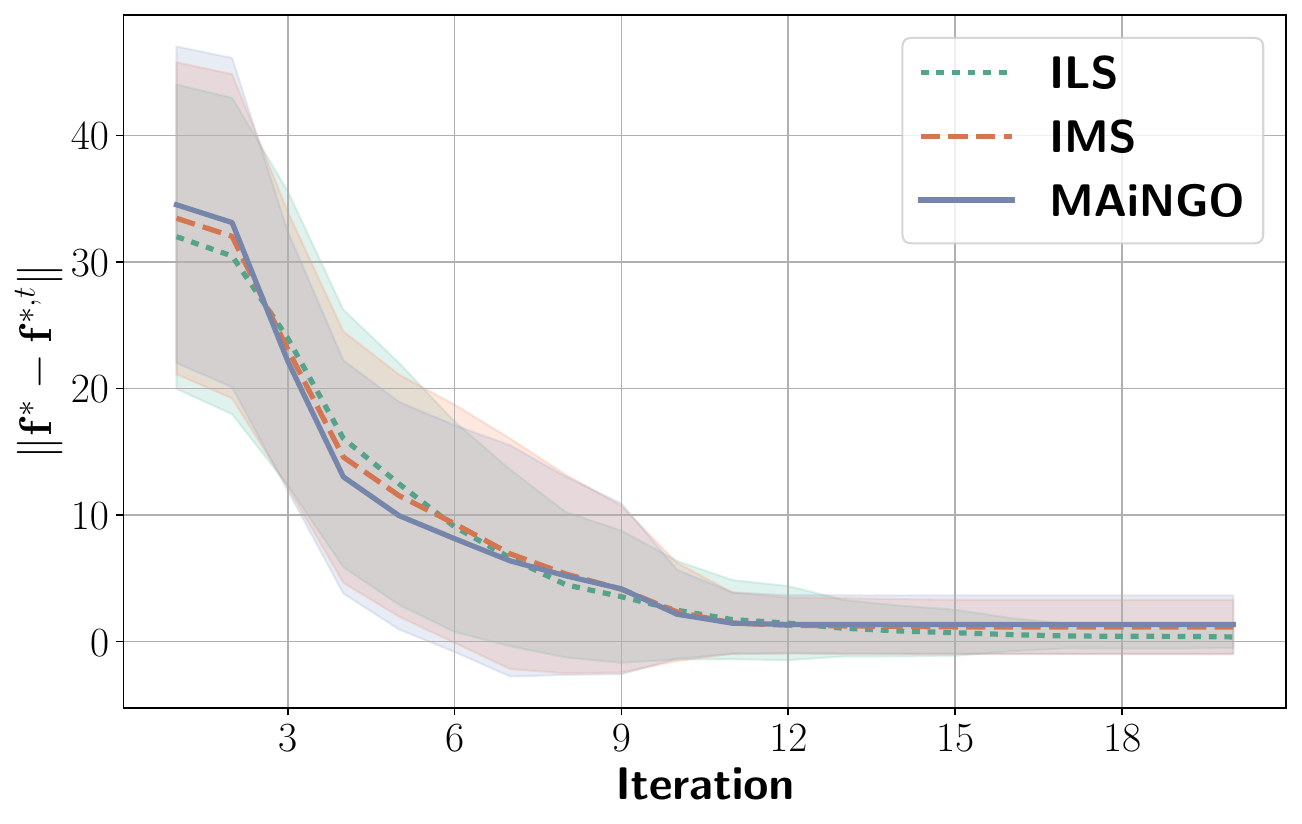}
        \subcaption{Müller-Brown, $\kappa = 0.5$ and $N = 20$.}
        \label{fig:MB-3t-2k-regret:k0.5_N20}
    \end{minipage}%
    \hfill%
    \begin{minipage}[t]{0.47\textwidth}
        \centering
        \includegraphics[trim={0.2cm 0.2cm 0.2cm 0cm}, clip, width=\textwidth]
        {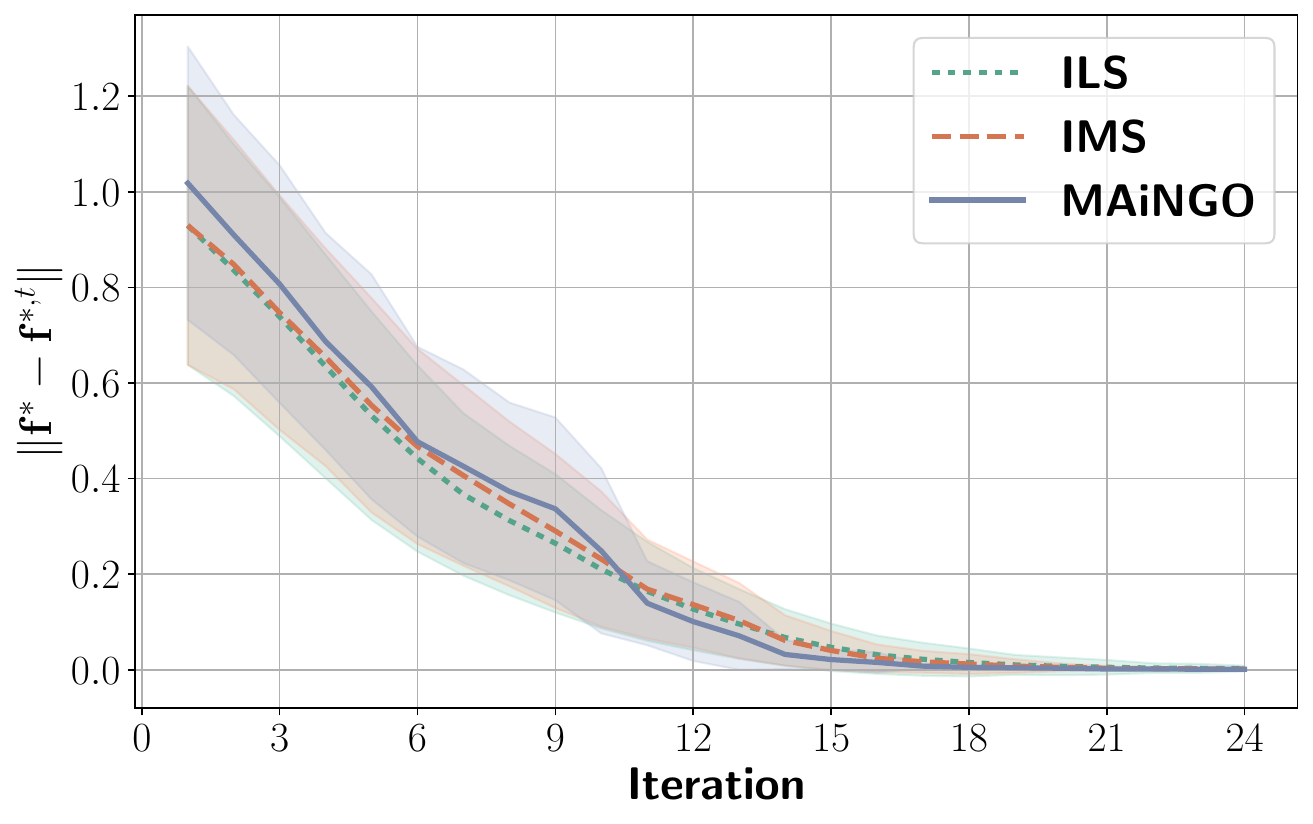}
        \subcaption{2D Camelback, $\kappa = 2$ and $N = 3$.}
        \label{fig:CB-regret}
    \end{minipage}%
    \par\medskip
    \begin{minipage}[t]{0.47\textwidth}
        \centering
        \includegraphics[trim={0.2cm 0.2cm 0.2cm 0cm}, clip, width=\textwidth]
        {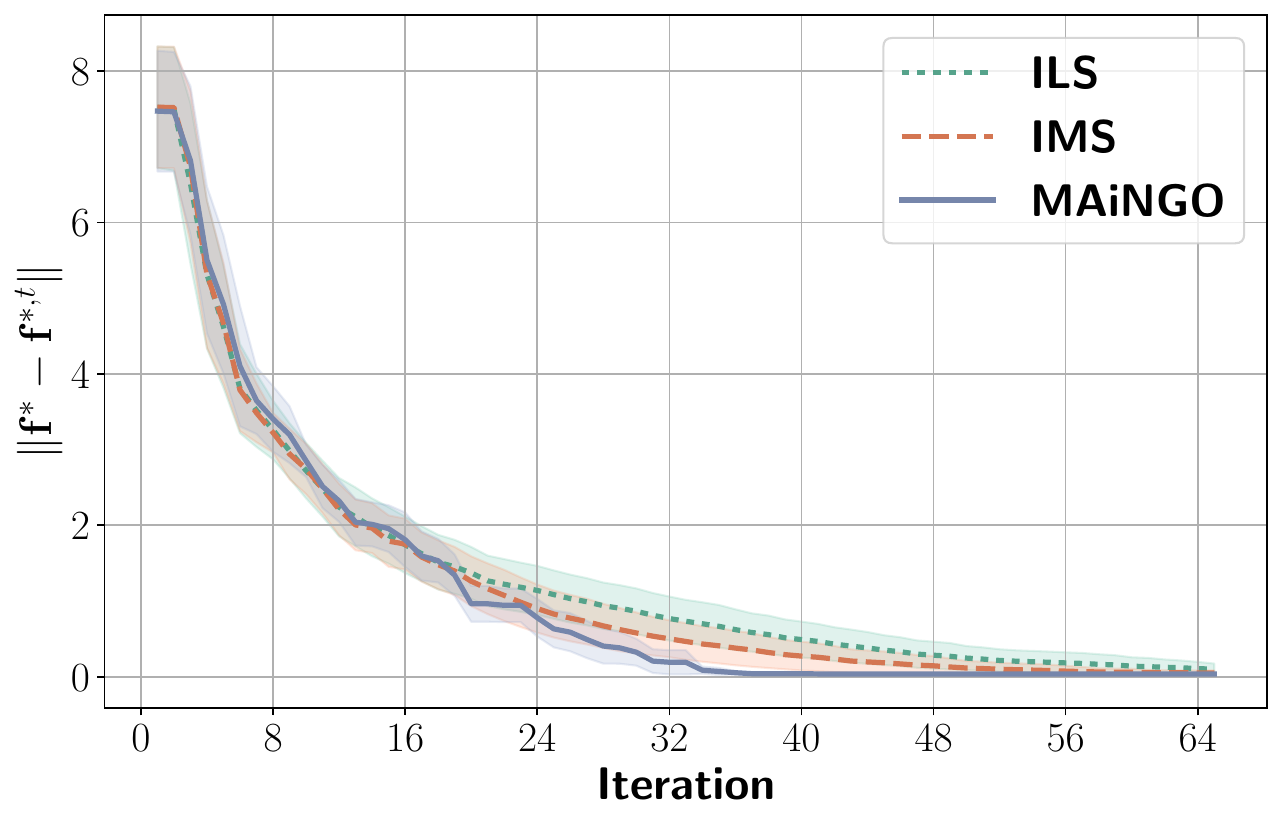}
        \subcaption{3D Ackley, $\kappa = 2$ and $N = 4$.}
        \label{fig:A3D-Regret}
    \end{minipage}%
    \hfill%
    \begin{minipage}[t]{0.47\textwidth}
        \centering
        \includegraphics[trim={0.2cm 0.2cm 0.2cm 0cm}, clip, width=\textwidth]
        {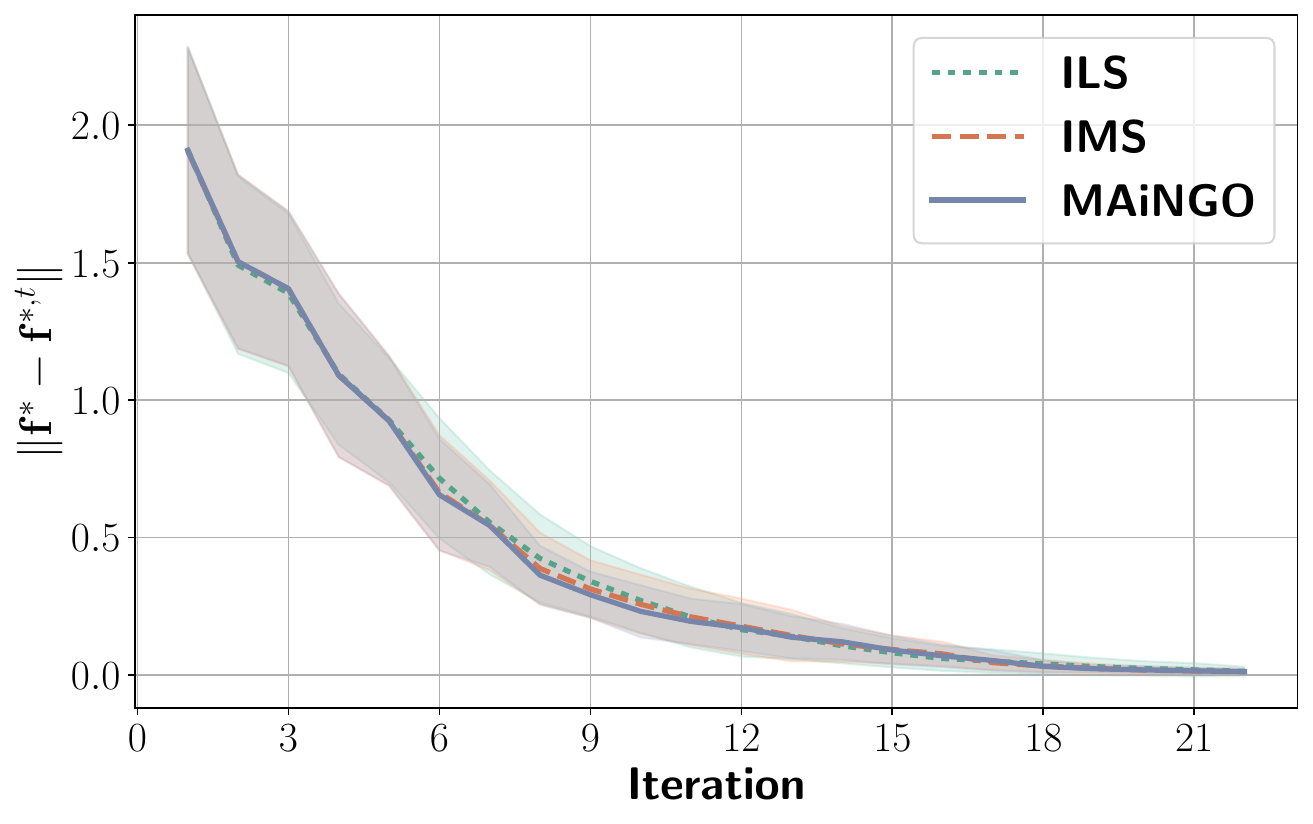}
        \subcaption{4D Hartmann, $\kappa = 2$ and $N = 5$.}
        \label{fig:H4D-Regret}
    \end{minipage}%
    \par\medskip
    \begin{minipage}[t]{0.47\textwidth}
        \centering
        \includegraphics[trim={0.2cm 0.2cm 0.2cm 0cm}, clip, width=\textwidth]
        {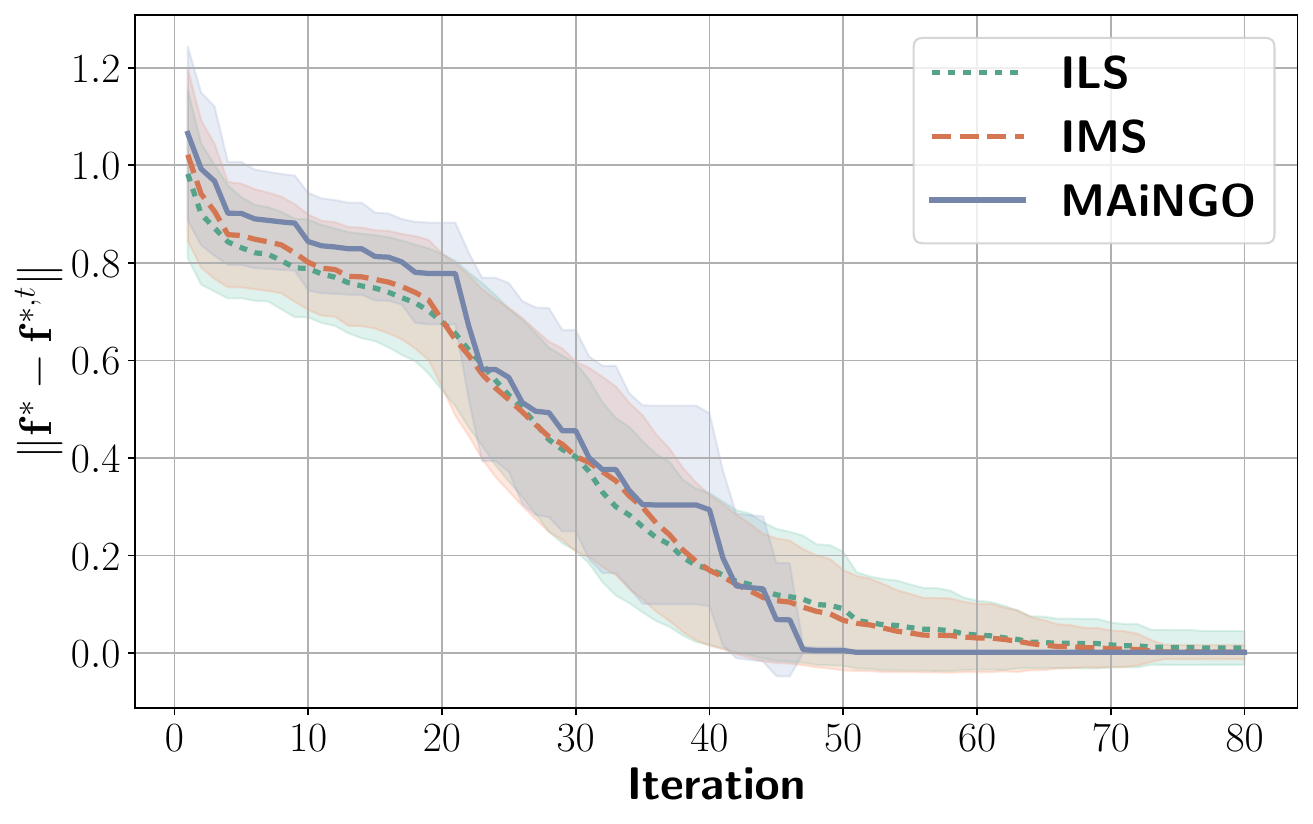}
        \subcaption{2D GKLS, $\kappa = 3$ and $N = 3$.}
        \label{fig:2GKLS-Regret}
    \end{minipage}%
    \hfill%
    \begin{minipage}[t]{0.47\textwidth}
        \centering
        \includegraphics[trim={0.2cm 0.2cm 0.2cm 0cm}, clip, width=\textwidth]
        {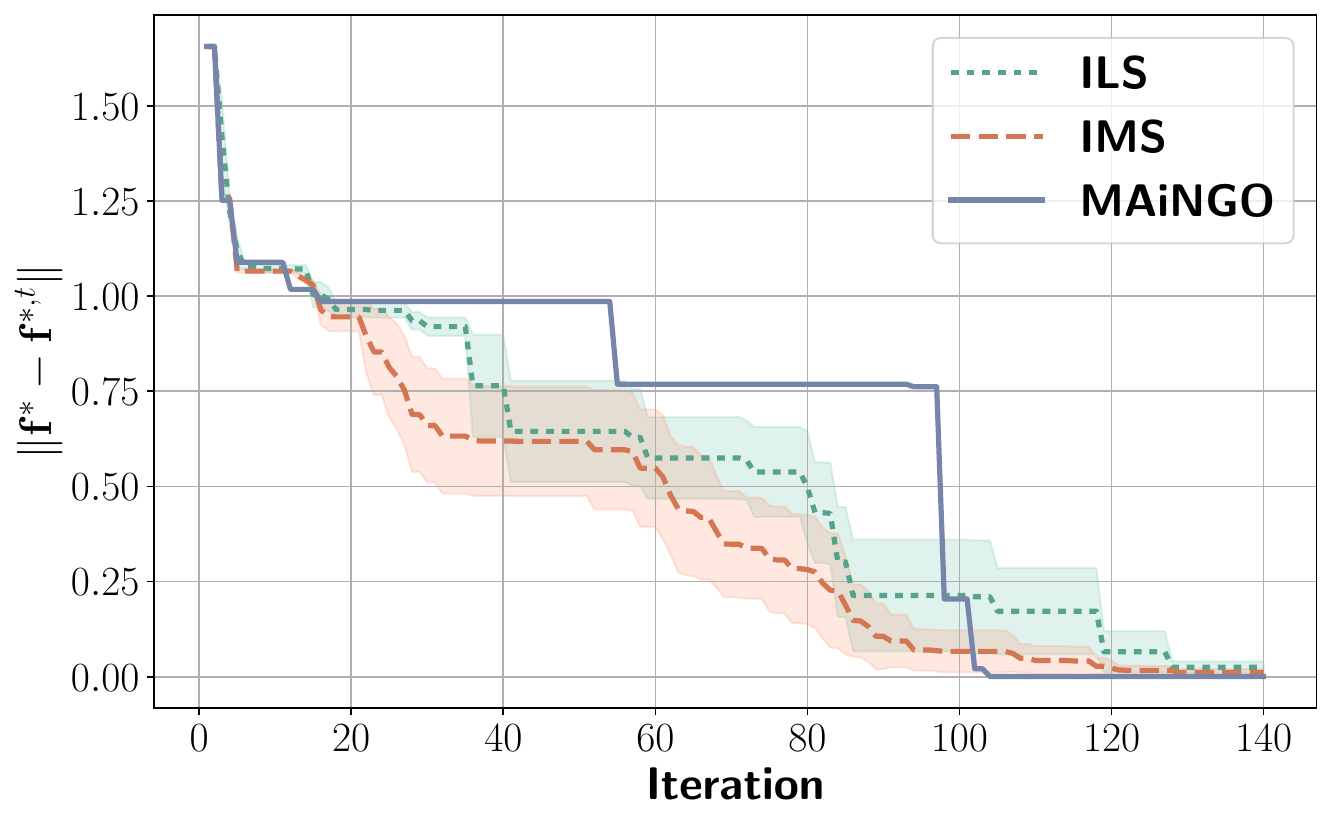}
        \subcaption{3D GKLS, $\kappa = 3$ and $N = 4$.}
        \label{fig:3GKLS-Regret}
    \end{minipage}
    \caption{
        Simple regret plots (Part 2) for runs successfully converging to a globally near-optimal solution, considering only experiments where at least one run reached a globally near-optimal solution for all three solvers using \eqref{eq:lcb}. This subset of data is identical to the ones analyzed in ~\cref{sec:more_case_studies,sec:iterations_to_convergence}.
        Case study formulations are detailed in \cref{sec:test-case-formulation}, with experimental setup, number of runs, and termination criteria in \cref{sec:set-up-details}.
        Results align with the trends observed in \cref{fig:MB-3t-2k-conv} and \cref{fig:other-test-cases-conv}, showing that MAiNGO generally terminates in fewer iterations.
        Lines represents the mean, while the shaded bands indicate $\pm 0.5$ standard deviations.
    }
    \label{fig:combined-regret-plots:part:2}
\end{figure}

\clearpage

\subsection{Considering All Runs, Despite Where They Converged}
\begin{figure}[htbp]
    \centering
    \begin{minipage}[t]{0.47\textwidth}
        \centering
        \includegraphics[trim={0.2cm 0.2cm 0.2cm 0cm}, clip, width=\textwidth]
        {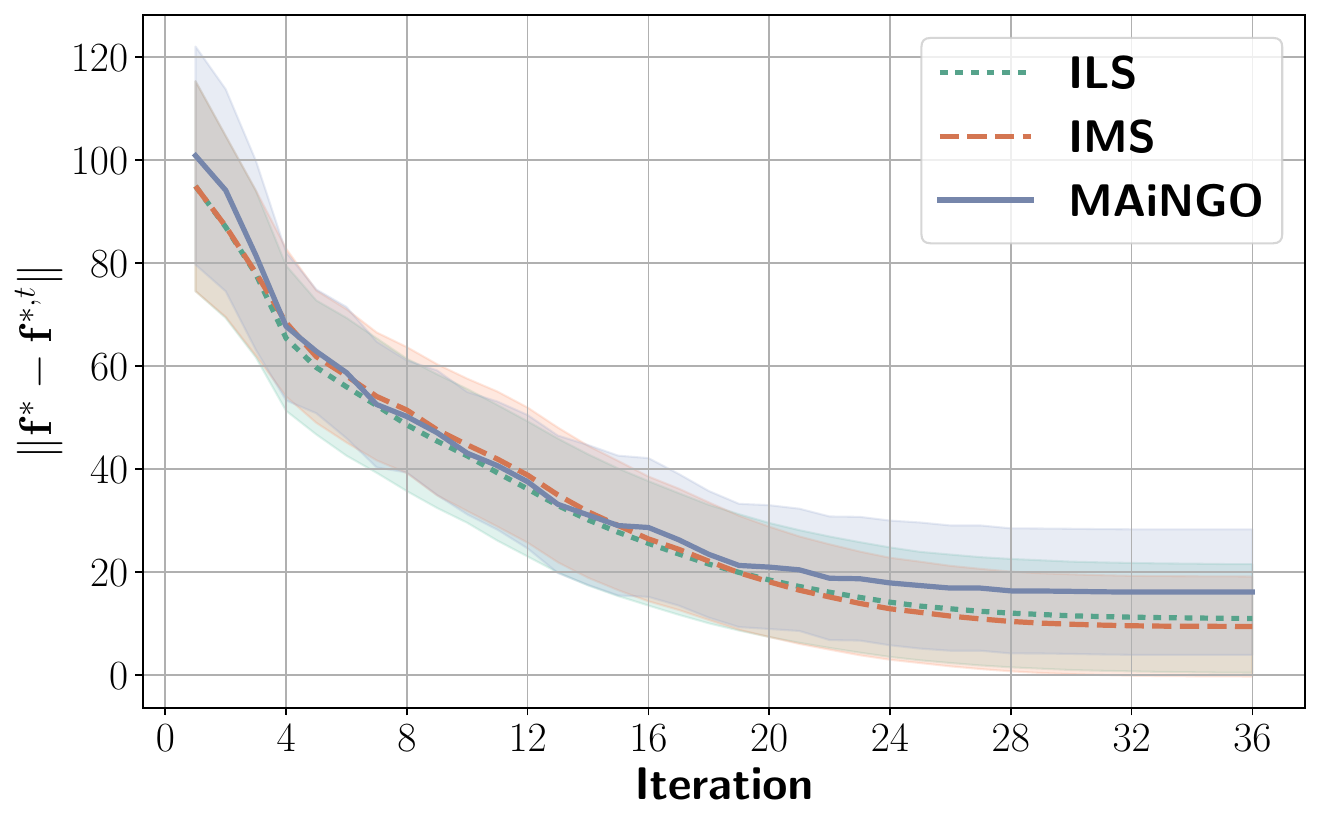}
        \subcaption{Müller-Brown, $\kappa = 2$ and $N = 3$.}
        \label{fig:MB-3t-2k-regret:k2_N2_All}
    \end{minipage}%
    \hfill%
    \begin{minipage}[t]{0.47\textwidth}
        \centering
        \includegraphics[trim={0.2cm 0.2cm 0.2cm 0cm}, clip, width=\textwidth]
        {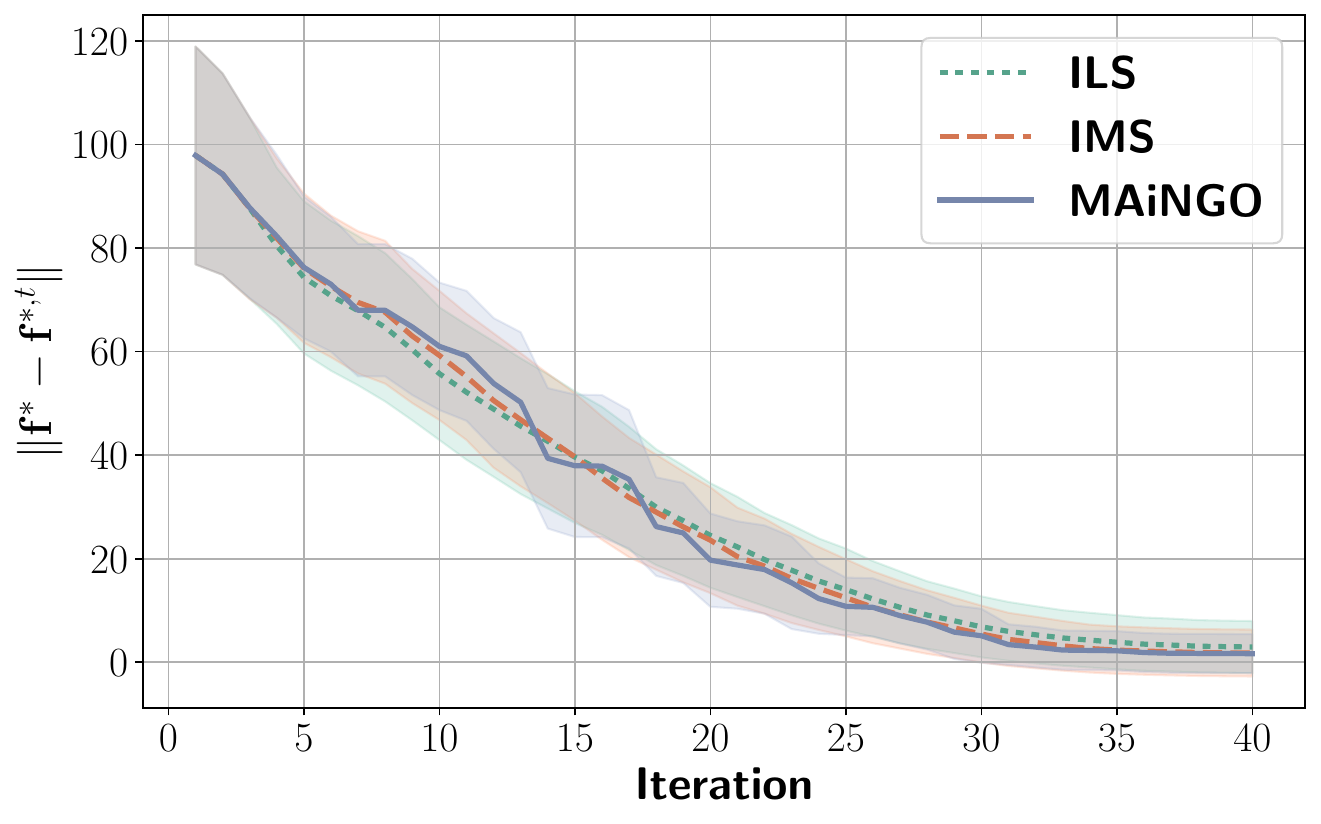}
        \subcaption{Müller-Brown, $\kappa = 3$ and $N = 3$.}
        \label{fig:MB-3t-2k-regret:k3_N3_All}
    \end{minipage}%
    \par\medskip
    \begin{minipage}[t]{0.47\textwidth}
        \centering
        \includegraphics[trim={0.2cm 0.2cm 0.2cm 0cm}, clip, width=\textwidth]
        {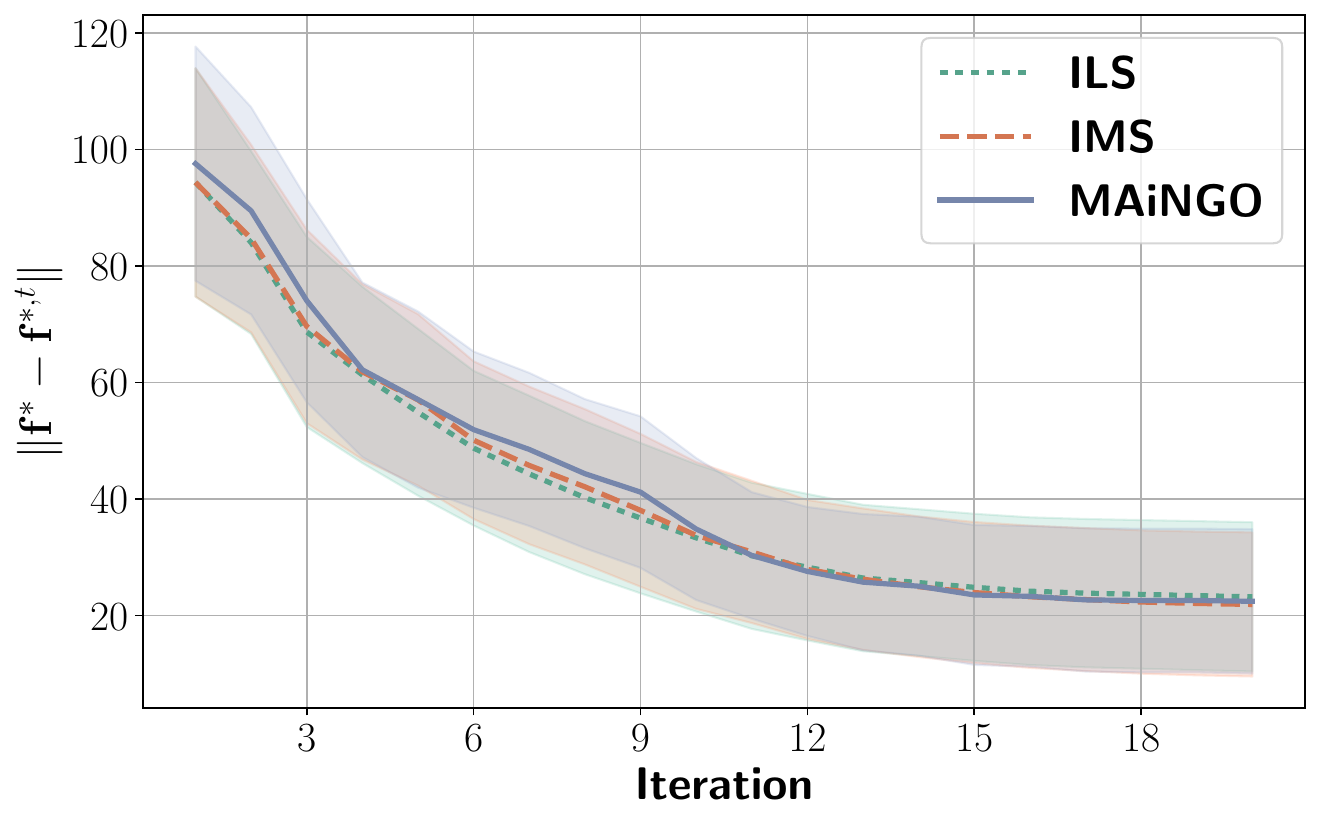}
        \subcaption{Müller-Brown, $\kappa =$ ~\ref{eq:kappa_t_Kandasamy_2015} and $N = 3$.}
        \label{fig:MB-3t-2k-regret:kK_N3_All}
    \end{minipage}%
    \hfill%
    \begin{minipage}[t]{0.47\textwidth}
        \centering
        \includegraphics[trim={0.2cm 0.2cm 0.2cm 0cm}, clip, width=\textwidth]
        {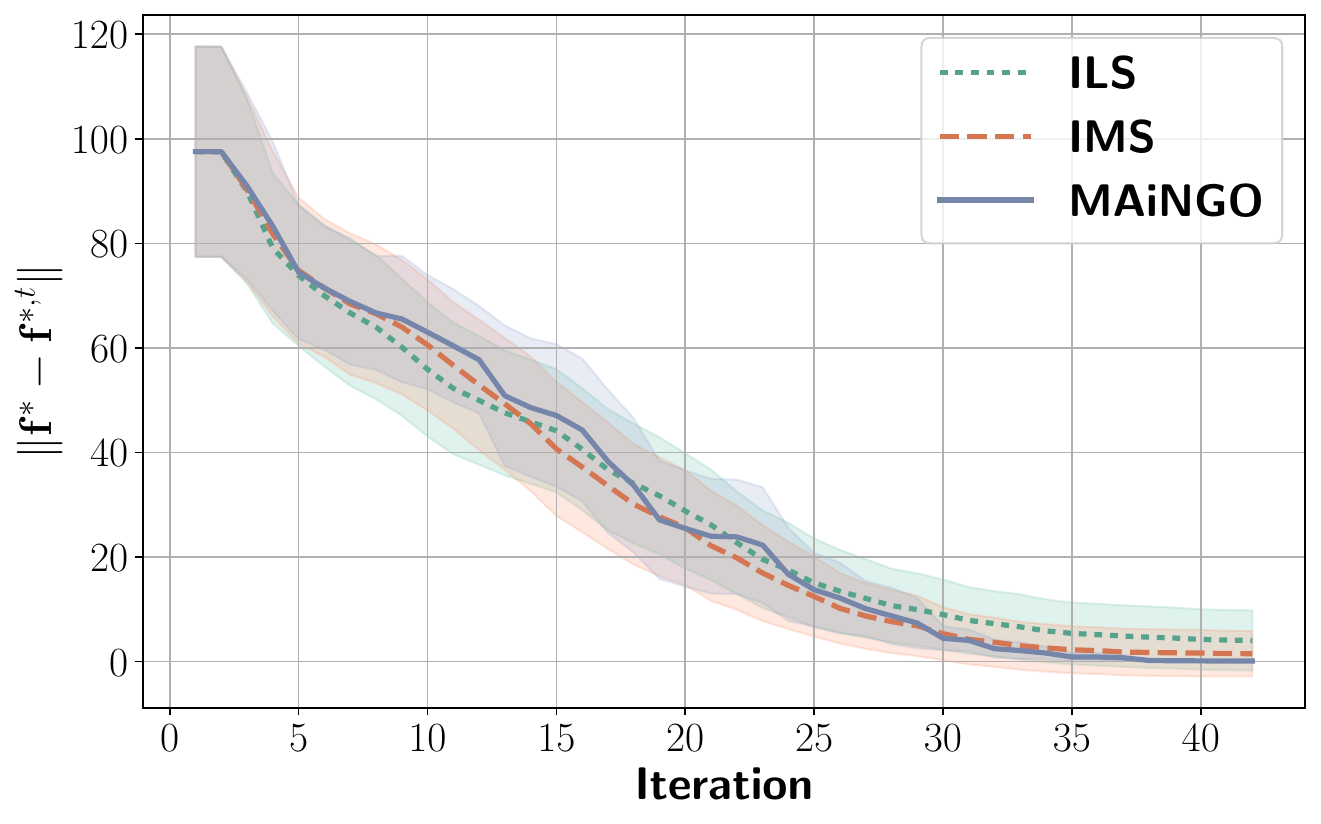}
        \subcaption{Müller-Brown, $\kappa =$ ~\ref{eq:kappa_t_Srinivas_2012} and $N = 3$.}
        \label{fig:MB-3t-2k-regret:kS_N3_All}
    \end{minipage}%
    \par\medskip
    \begin{minipage}[t]{0.47\textwidth}
        \centering
        \includegraphics[trim={0.2cm 0.2cm 0.2cm 0cm}, clip, width=\textwidth]
        {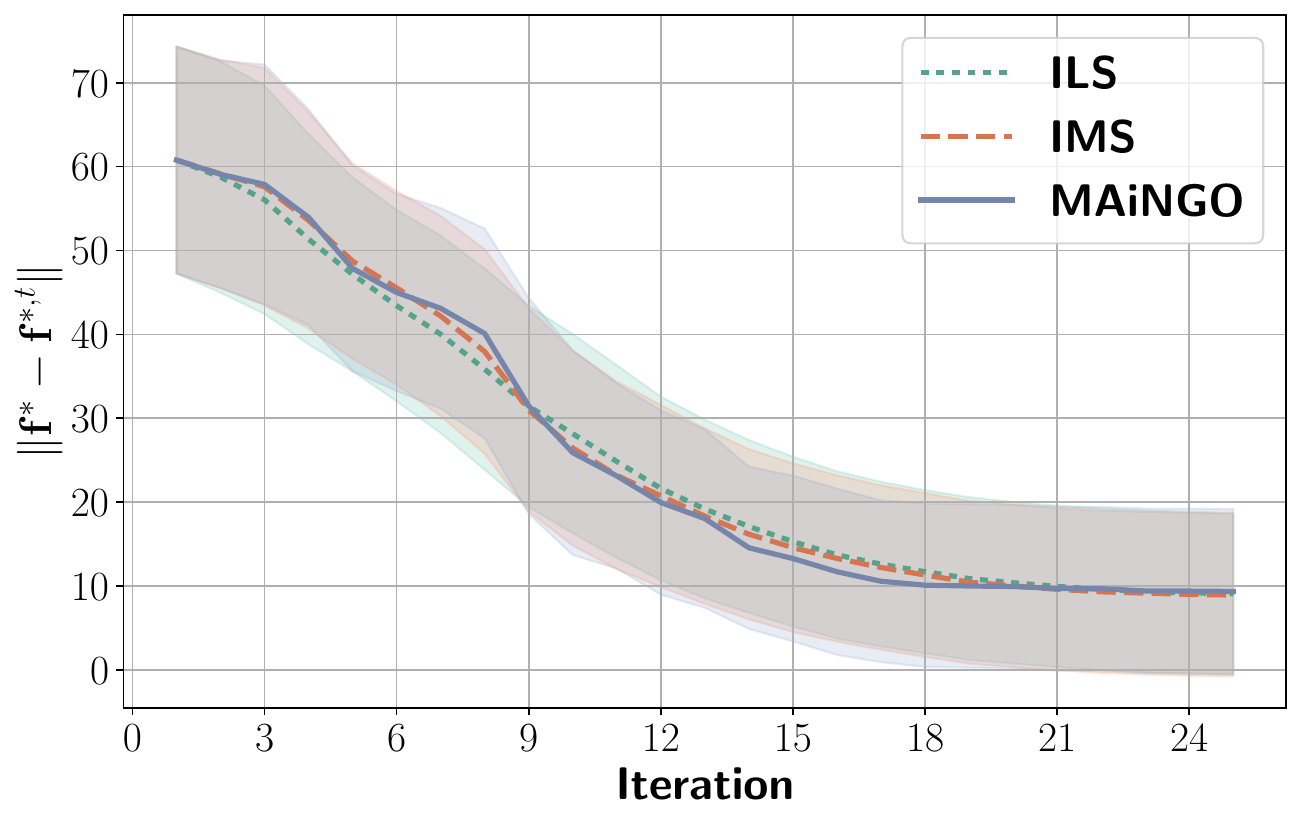}
        \subcaption{Müller-Brown, $\kappa = 2$ and $N = 10$.}
        \label{fig:MB-3t-2k-regret:k2_N10_All}
    \end{minipage}%
    \hfill%
    \begin{minipage}[t]{0.47\textwidth}
        \centering
        \includegraphics[trim={0.2cm 0.2cm 0.2cm 0cm}, clip, width=\textwidth]
        {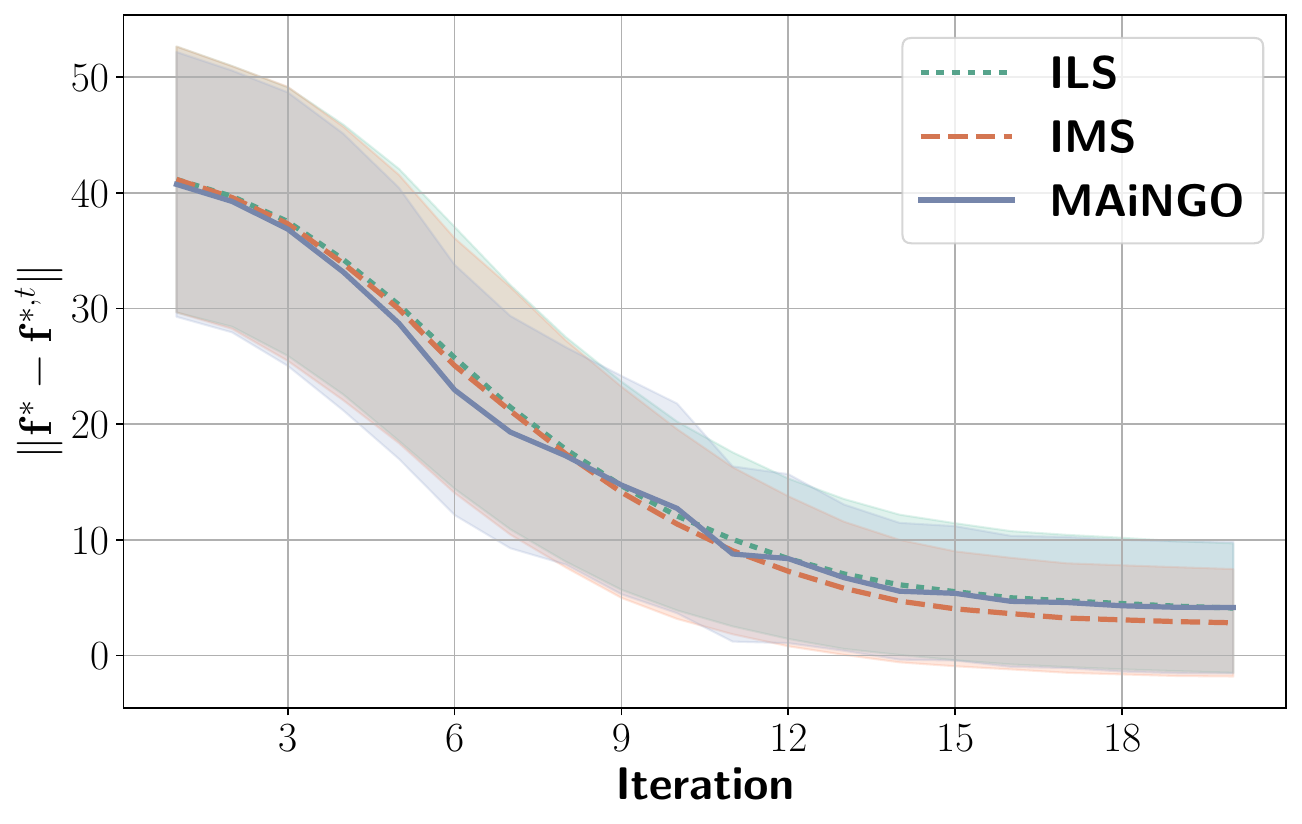}
        \subcaption{Müller-Brown, $\kappa = 2$ and $N = 20$.}
        \label{fig:MB-3t-2k-regret:k2_N20_All}
    \end{minipage}%
    \caption{
      Simple regret plots (Part 1).
    }
    \label{fig:combined-regret-plots_All:part:1}
\end{figure}

\begin{figure}[htbp]\ContinuedFloat
    \centering
    \begin{minipage}[t]{0.47\textwidth}
        \centering
        \includegraphics[trim={0.2cm 0.2cm 0.2cm 0cm}, clip, width=\textwidth]
        {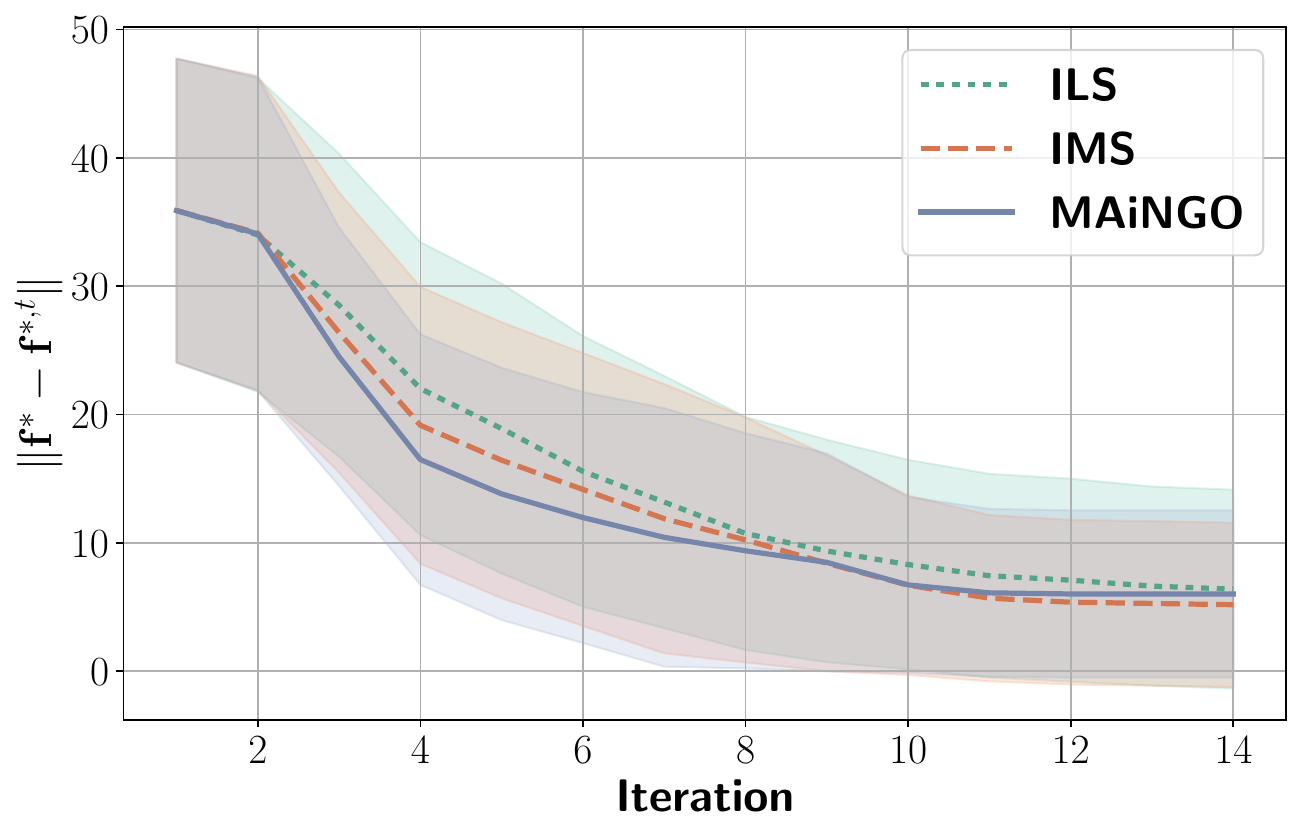}
        \subcaption{Müller-Brown, $\kappa = 0.5$ and $N = 20$.}
        \label{fig:MB-3t-2k-regret:k0.5_N20_All}
    \end{minipage}%
    \hfill%
    \begin{minipage}[t]{0.47\textwidth}
        \centering
        \includegraphics[trim={0.2cm 0.2cm 0.2cm 0cm}, clip, width=\textwidth]
        {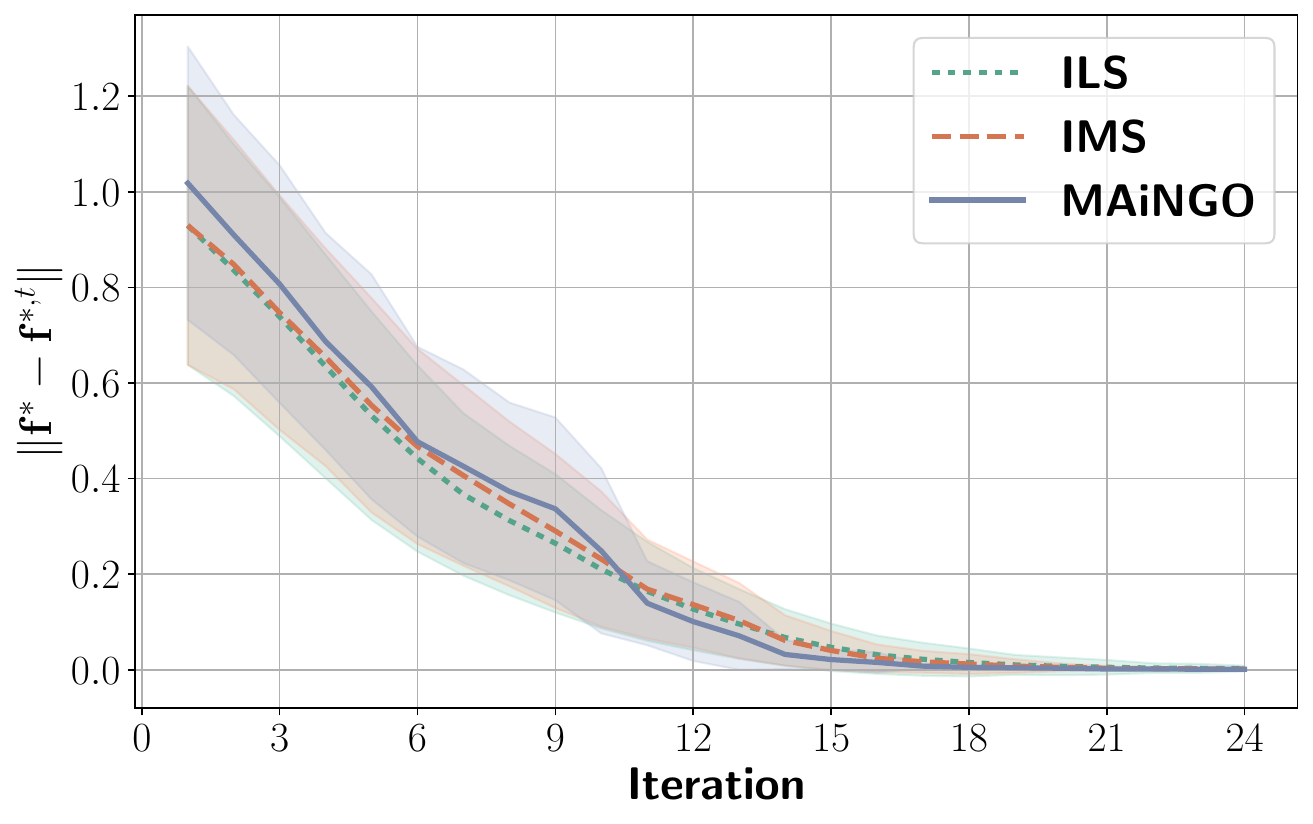}
        \subcaption{2D Camelback, $\kappa = 2$ and $N = 3$.}
        \label{fig:CB-regret_All}
    \end{minipage}%
    \par\medskip
    \begin{minipage}[t]{0.47\textwidth}
        \centering
        \includegraphics[trim={0.2cm 0.2cm 0.2cm 0cm}, clip, width=\textwidth]
        {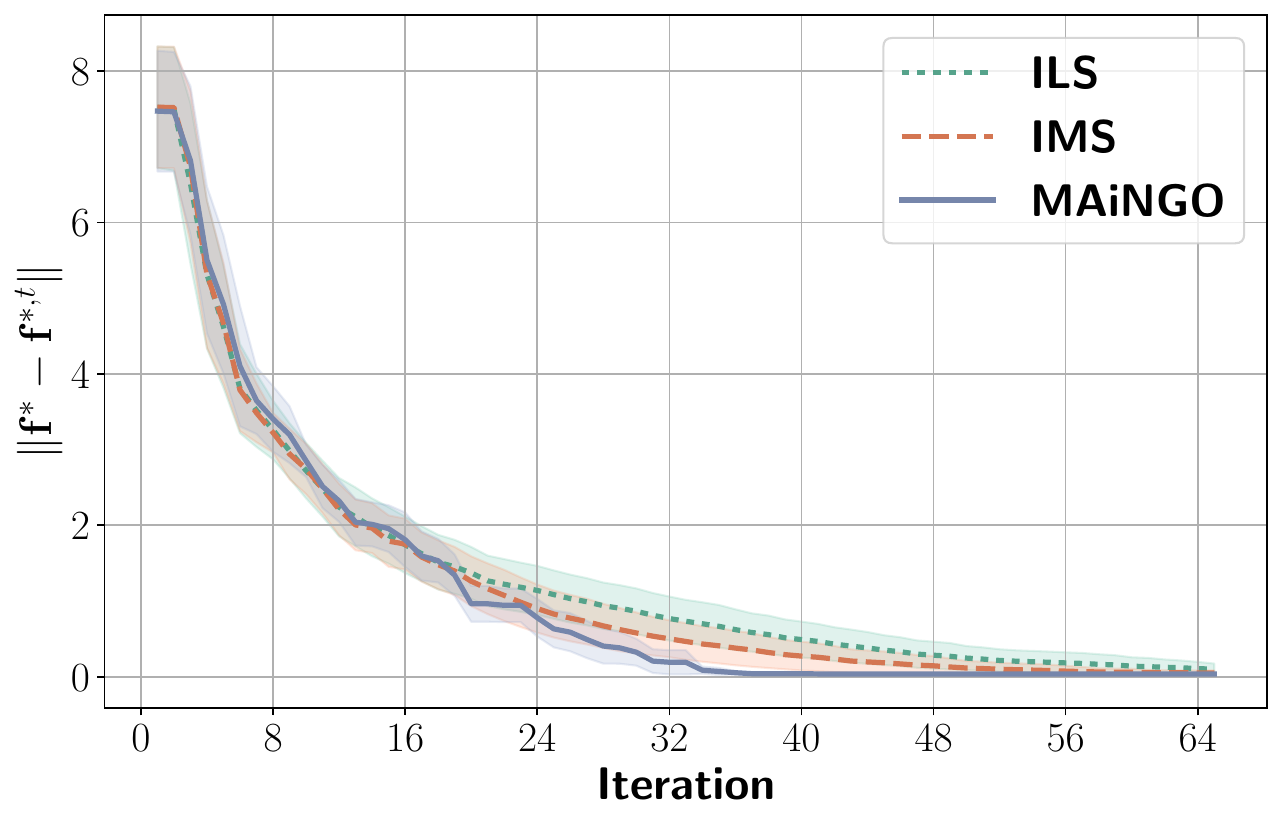}
        \subcaption{3D Ackley, $\kappa = 2$ and $N = 4$.}
        \label{fig:A3D-Regret_All}
    \end{minipage}%
    \hfill%
    \begin{minipage}[t]{0.47\textwidth}
        \centering
        \includegraphics[trim={0.2cm 0.2cm 0.2cm 0cm}, clip, width=\textwidth]
        {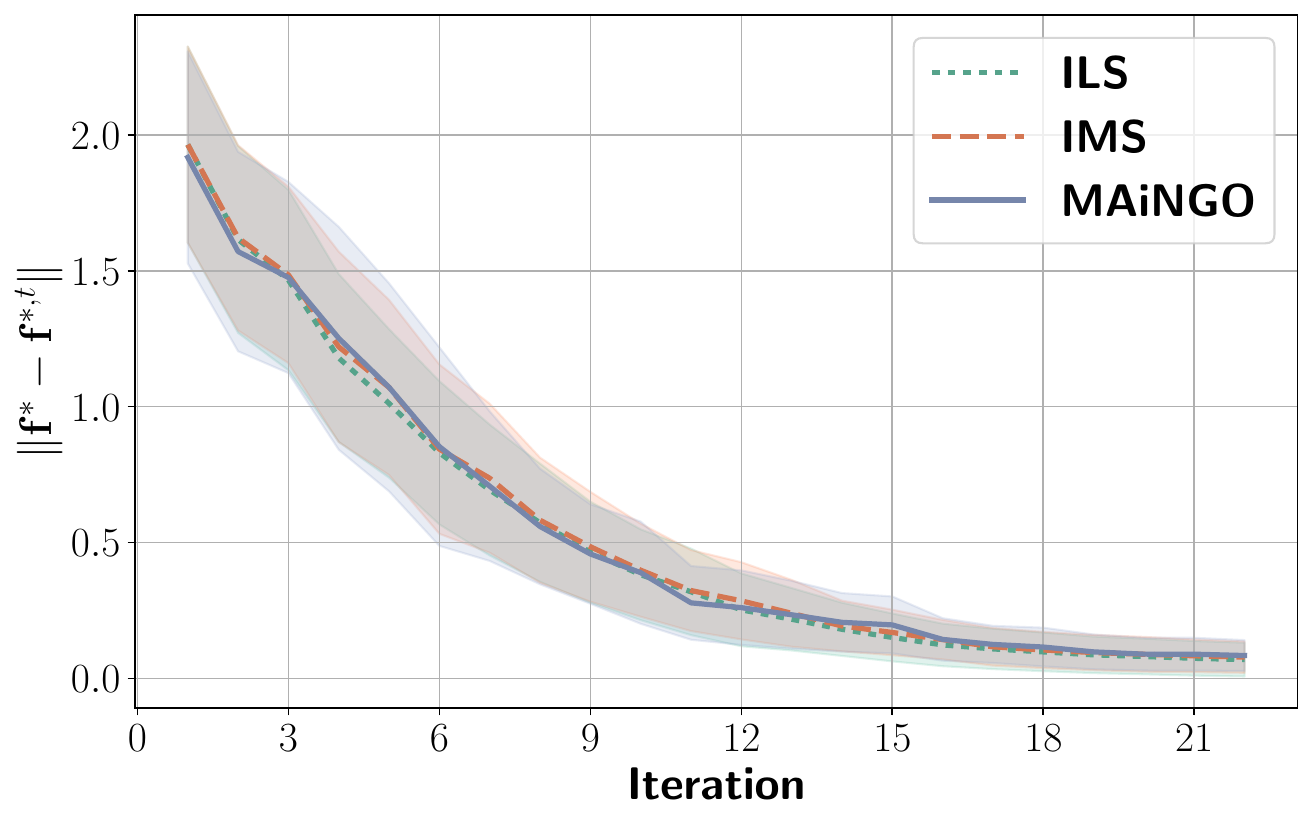}
        \subcaption{4D Hartmann, $\kappa = 2$ and $N = 5$.}
        \label{fig:H4D-Regret_All}
    \end{minipage}%
        \par\medskip
    \begin{minipage}[t]{0.47\textwidth}
        \centering
        \includegraphics[trim={0.2cm 0.2cm 0.2cm 0cm}, clip, width=\textwidth]
        {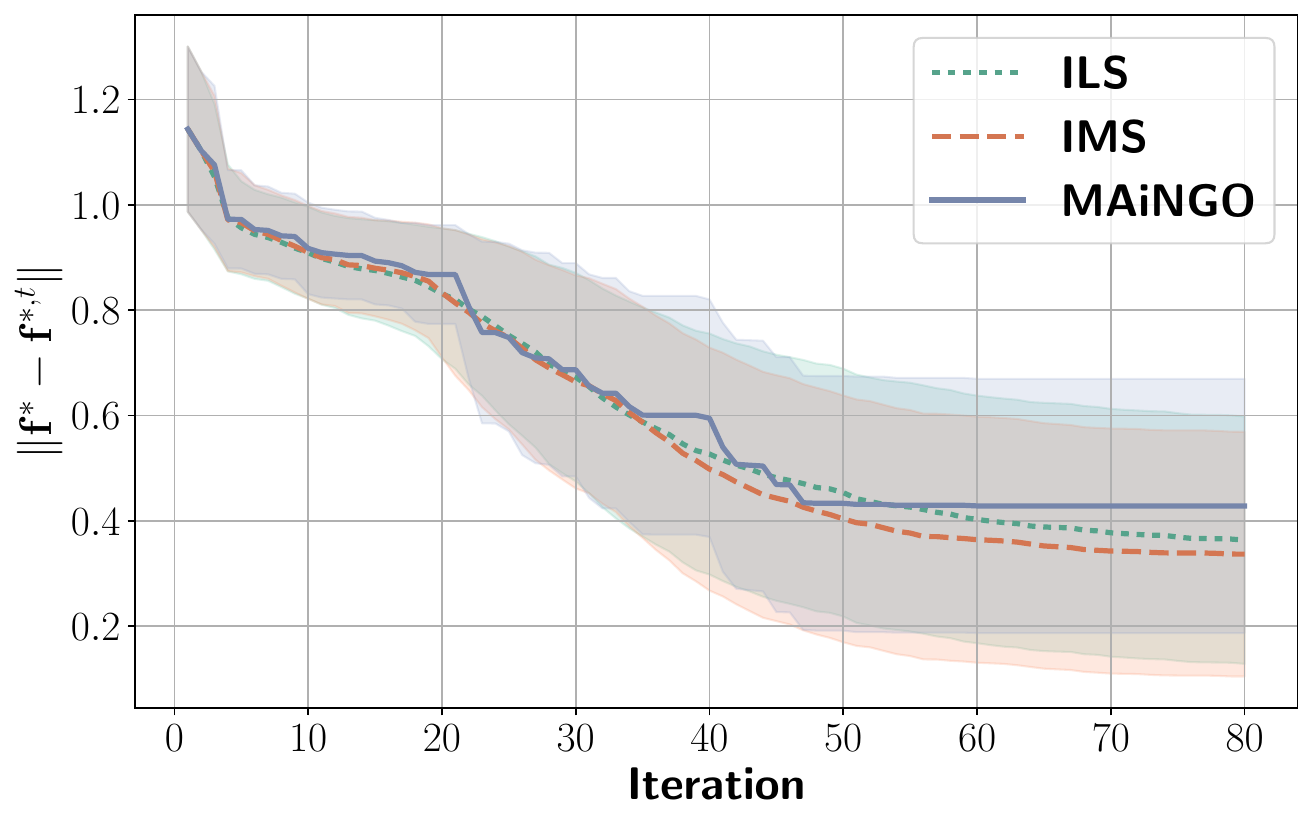}
        \subcaption{2D GKLS, $\kappa = 3$ and $N = 3$.}
        \label{fig:2GKLS-Regret_All}
    \end{minipage}%
    \hfill%
    \begin{minipage}[t]{0.47\textwidth}
        \centering
        \includegraphics[trim={0.2cm 0.2cm 0.2cm 0cm}, clip, width=\textwidth]
        {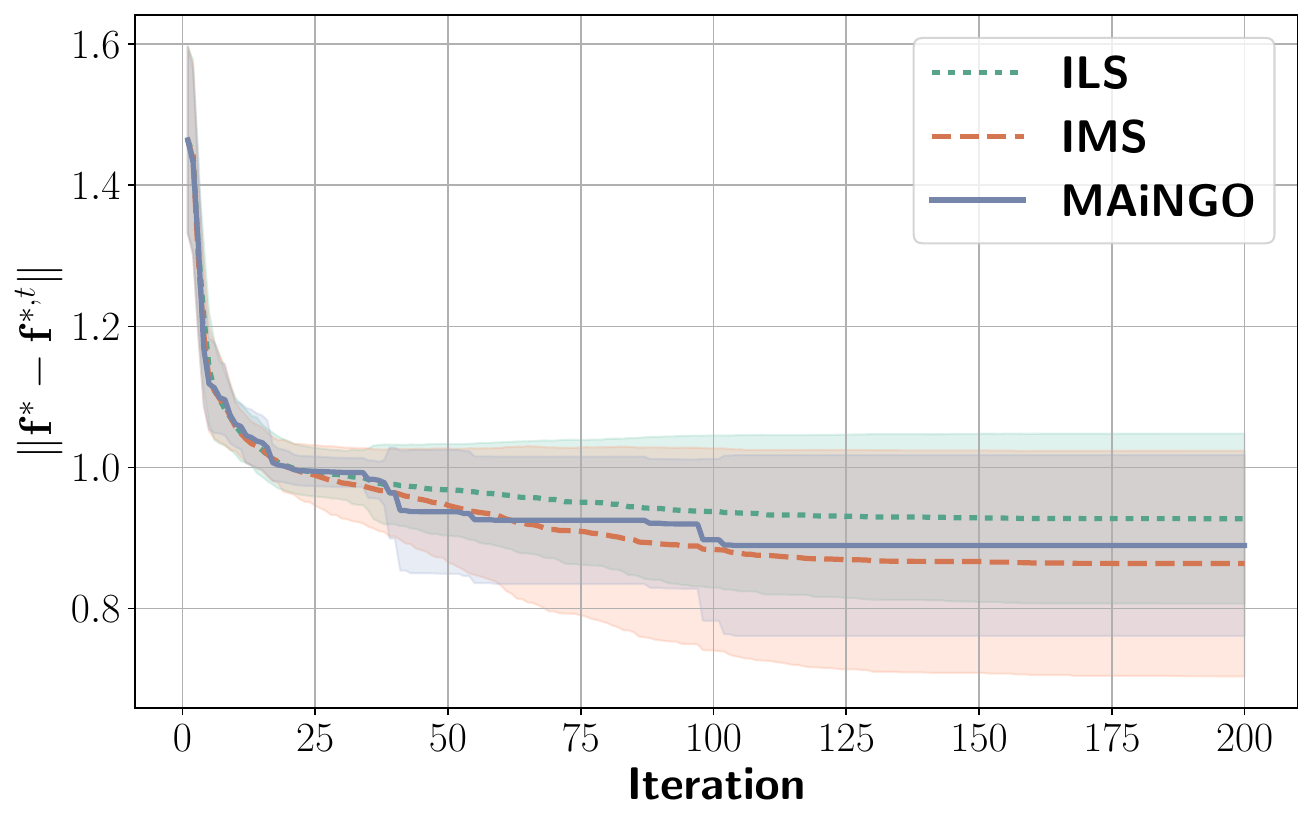}
        \subcaption{3D GKLS, $\kappa = 3$ and $N = 4$.}
        \label{fig:3GKLS-Regret_All}
    \end{minipage}
    \caption{
        Simple regret plots (Part 2) for all conducted experiments and runs using \eqref{eq:lcb}.
        Case study formulations are detailed in \cref{sec:test-case-formulation}, with the experimental setup, number of runs, and termination criteria in \cref{sec:set-up-details}.
        Trajectories that do not approach 0 have a substantial probability of converging to a local minimum of the black box function.
        Lines represent the mean, while the shaded bands indicate $\pm 0.5$ standard deviations.
    }
    \label{fig:combined-regret-plots_All:part:2}
\end{figure}
\clearpage

\section{Iterations to Convergence for All Runs}
\label{sec:stats-with-all-min}
\begin{table}[ht]
    \centering
    \caption{
        Summary of convergence statistics for each case study, considering all runs of all experiments, despite which minimum each converged to.
        Mean, median, and standard deviation of iterations to convergence are reported for this subset of data.
        One-sided t-tests \textit{paired by experiment} (initial dataset) evaluate whether MAiNGO converges in fewer iterations than IMS or ILS.
        MAiNGO converges in fewer iterations than ILS for all case studies and significantly fewer iterations than IMS for all but the case studies using more exploratory acquisition functions ($\kappa=3$, $\kappa=$~\ref{eq:kappa_t_Srinivas_2012}).
        Details on the total number of experiments, runs, and termination criteria for each case study are provided in \cref{sec:set-up-details}.
    }
    \label{tab:mb-fast_All}
    \begin{tabular}{c c c c c c c c c}
        $\kappa$ & $N$ & Solver & Count & Mean & Median & Std & t-test & p-value \\ [0.5ex]
        \hline
        2 & 3   & ILS   & 1736  & 23.6 & 23.0 & 7.55 & 6.53 & \num{1.1E-8} \\
        2 & 3   & IMS   & 1736  & 23.0 & 23.0 & 5.40 & 3.84 & \num{1.9E-7} \\
        2 & 3   & MAiNGO& 1736  & 21.1 & 20.0 & 4.53 & - & -\\ [1ex]
        \hline
        3 & 3   & ILS   & 775   & 34.1 & 33.0 & 6.83 & 1.60 & \num{6.2E-2}\\
        3 & 3   & IMS   & 775   & 33.4 & 33.0 & 5.28 & 0.14 & \num{4.4E-1}\\
        3 & 3   & MAiNGO& 775   & 33.0 & 34.0 & 4.24 & -    & -\\
        \hline
        2 & 10  & ILS   & 961   & 18.5 & 18.0 & 6.10 & 3.21 & \num{1.6E-3}\\
        2 & 10  & IMS   & 961   & 18.0 & 18.0 & 4.83 & 2.29 & \num{1.5E-2}\\
        2 & 10  & MAiNGO& 961   & 17.0 & 17.0 & 3.76 & -    & -\\
       \hline
        2 & 20  & ILS   & 1271  & 14.0 & 13.0 & 5.01 & 3.70 & \num{3.2E-4}\\
        2 & 20  & IMS   & 1271  & 13.2 & 13.0 & 3.85 & 1.32 & \num{9.7E-2}\\
        2 & 20  & MAiNGO& 1271  & 12.9 & 13.0 & 3.37 & -    & -\\
        \hline
        \ref{eq:kappa_t_Srinivas_2012} & 3 & ILS   & 441 & 34.8 & 34.0 & 7.67 & -0.73 & \num{7.6E-1}\\
        \ref{eq:kappa_t_Srinivas_2012} & 3 & IMS   & 441 & 33.6 & 34.0 & 5.83 & -2.21 & \num{9.8E-1}\\
        \ref{eq:kappa_t_Srinivas_2012} & 3 & MAiNGO& 441 & 35.4 & 36.0 & 3.46 & -     & -\\
        \hline
        \ref{eq:kappa_t_Kandasamy_2015} & 3 & ILS   & 1071 & 12.8 & 12.0 & 4.80 & 2.47 & \num{8.5E-3}\\
        \ref{eq:kappa_t_Kandasamy_2015} & 3 & IMS   & 1071 & 12.6 & 12.0 & 4.41 & 3.56 & \num{4.1E-4}\\
        \ref{eq:kappa_t_Kandasamy_2015} & 3 & MAiNGO& 1071 & 11.8 & 12.0 & 4.04 & -    & -\\ [1ex]
    \end{tabular}
\end{table}

\begin{table}[ht]
    \centering
    \caption{
      Additional Case Studies (25 experiments with 25 runs each).
      All runs.
    }
    \label{tab:othergklstestcases-results}
    \begin{tabular}{l c c c c c c c c}
         Test Function & Solver  & Count & Mean & Median & Std & t-test & p-value \\ [0.5ex]
         \hline
         \multirow{3}{*}{2D GKLS}
         & ILS    & 625 & 35.8 & 32.0 & 21.56 & 1.905 & \num{0.034} \\
         & IMS    & 625 & 34.0 & 32.0 & 16.60 & 1.640 & \num{0.057} \\
         & MAiNGO & 625 & 30.8 & 28.0 & 12.08 & -     & - \\ [1ex]
         \hline
         \multirow{3}{*}{3D GKLS}
         & ILS    & 625 & 38.0 & 24.0 & 32.47 & 1.27 & \num{0.108} \\
         & IMS    & 625 & 41.6 & 25.0 & 35.70 & 3.63 & \num{6.7E-4} \\
         & MAiNGO & 625 & 34.2 & 23.0 & 21.48 & -    & - \\ [1ex]
         \hline
         \multirow{3}{*}{4D GKLS}
         & ILS    & 625 & 33.2 & 33.0 & 3.69 & 6.84 & \num{2.25E-7} \\
         & IMS    & 625 & 32.6 & 32.0 & 1.42 & 9.06 & \num{1.64E-9} \\
         & MAiNGO & 625 & 31.1 & 31.0 & 1.24 & -    & - \\ [1ex]
    \end{tabular}
\end{table}

\begin{figure}[htbp]
    \begin{minipage}[t]{0.31\textwidth}
        \centering
        \includegraphics[trim={0cm 0cm 0cm 0cm}, clip,width=\textwidth]
        {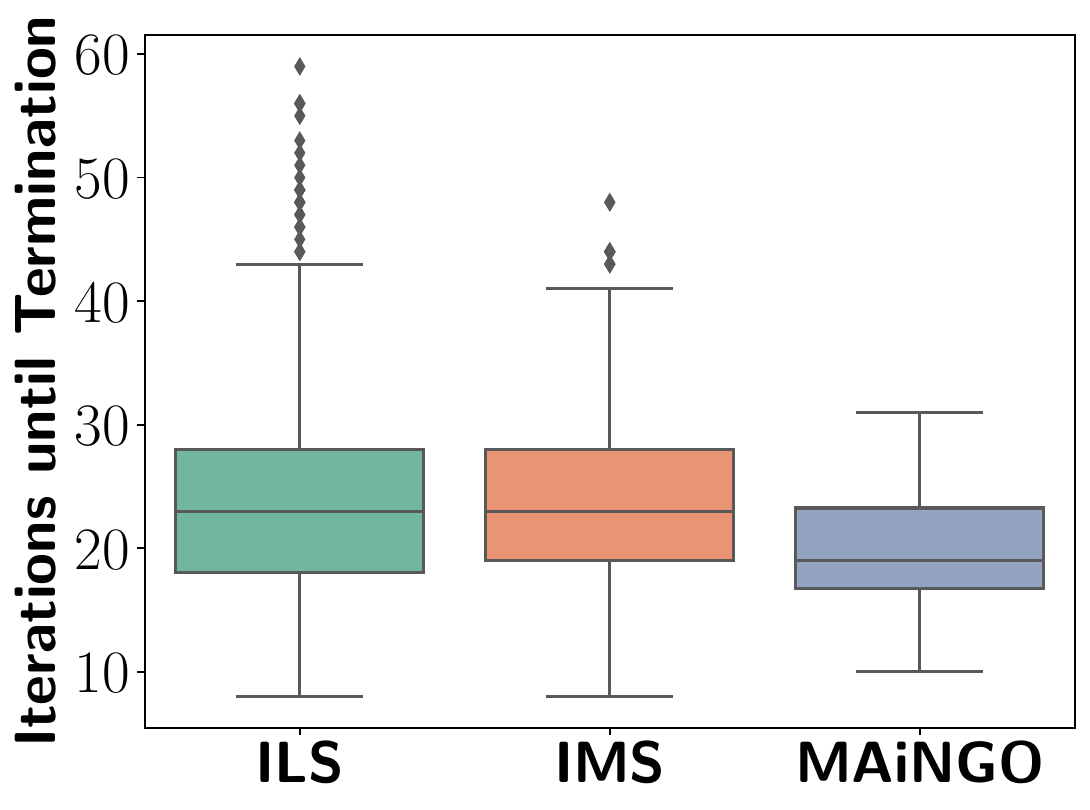}
        \subcaption{
            $\kappa = 2$ and $N = 3$.
        }
    \end{minipage}%
    \hspace{\fill}%
    \begin{minipage}[t]{0.31\textwidth}
        \centering
        \includegraphics[trim={0cm 0cm 0cm 0cm}, clip,width=\linewidth]
        {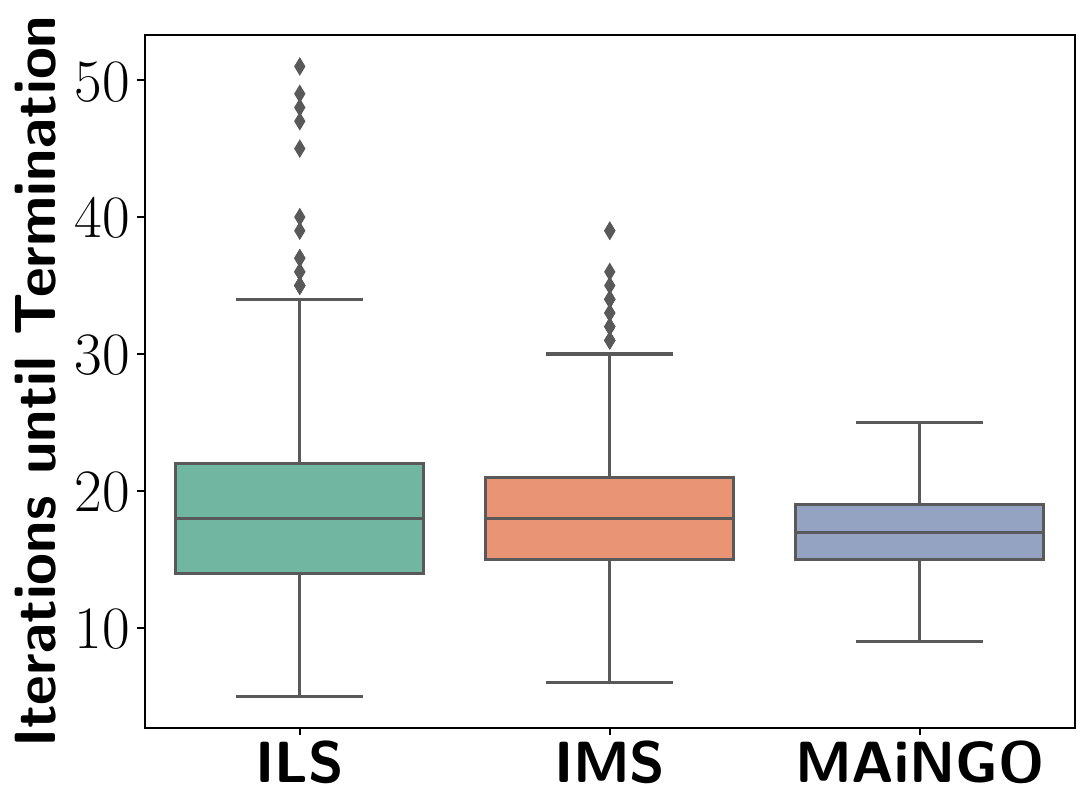}
        \subcaption{
            $\kappa = 2$ and $N = 10$.
        }
    \end{minipage}%
    \hspace{\fill}%
    \begin{minipage}[t]{0.31\textwidth}
        \centering
        \includegraphics[trim={0cm 0cm 0cm 0cm}, clip,width=\linewidth]
        {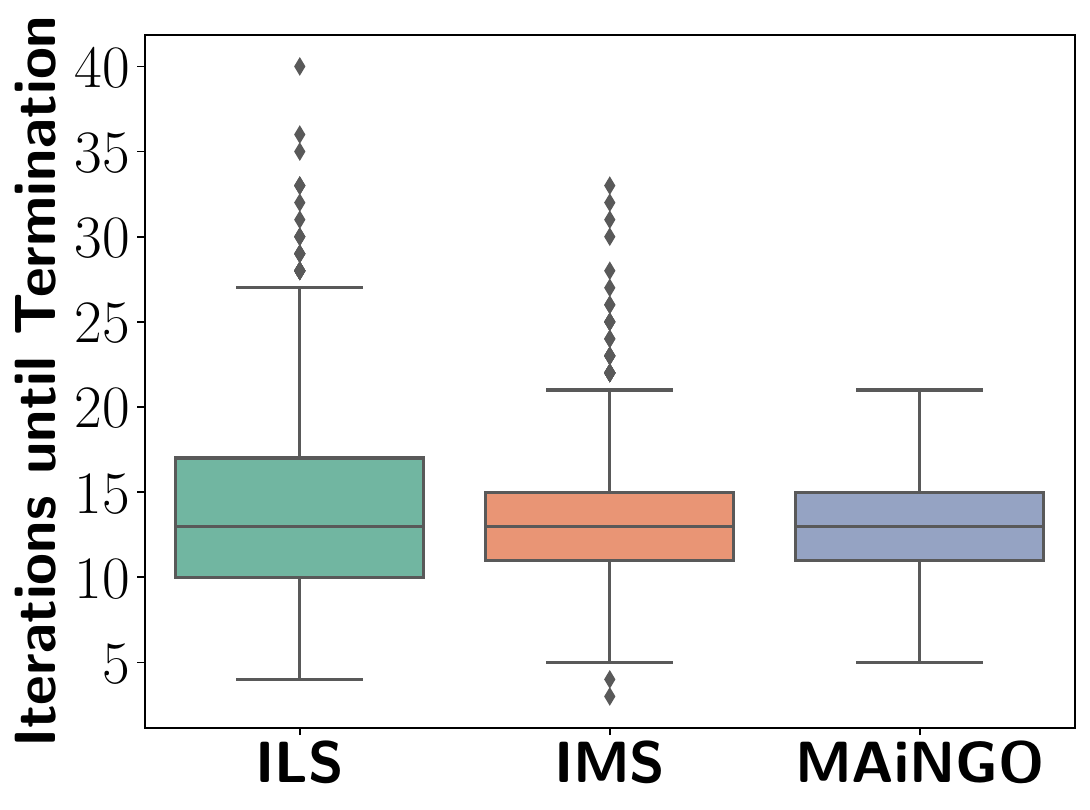}
        \subcaption{
            $\kappa = 2$ and $N = 20$.
        }
    \end{minipage}%
    \par\medskip
    \begin{minipage}[t]{0.31\textwidth}
        \centering
        \includegraphics[trim={0cm 0cm 0cm 0cm}, clip,width=\linewidth]
        {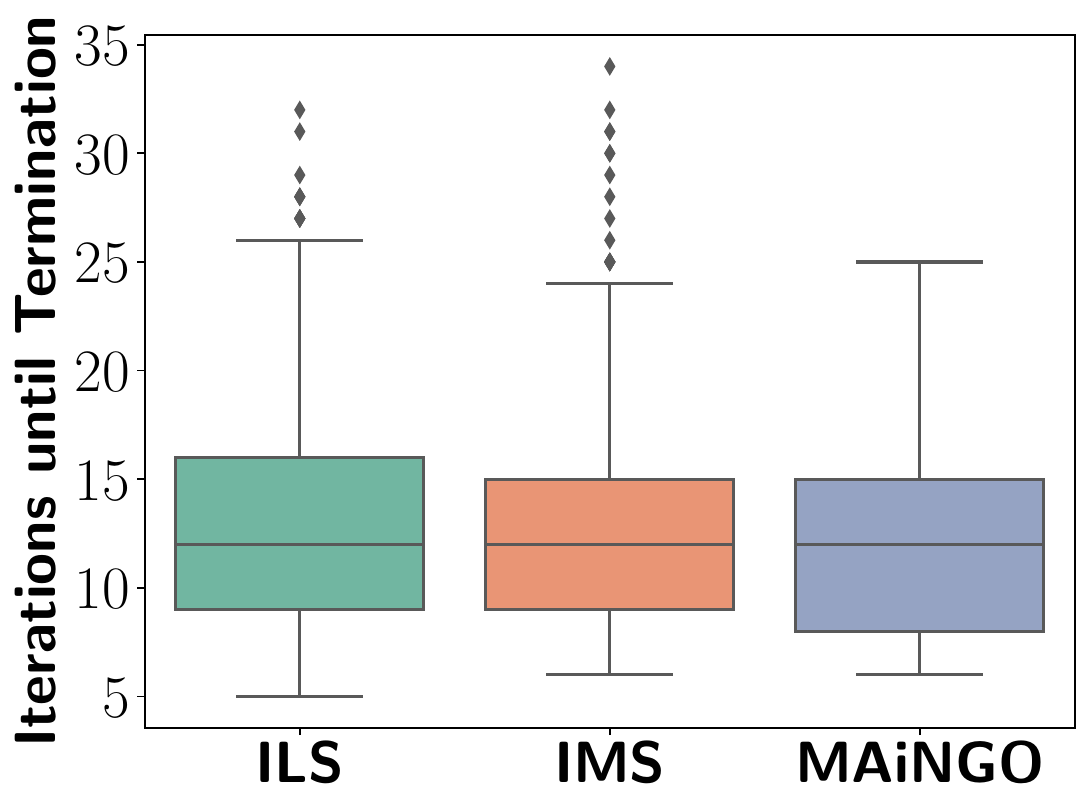}
        \subcaption{
            $\kappa = $~\ref{eq:kappa_t_Kandasamy_2015} and $N = 3$.
        }
    \end{minipage}%
    \hspace{\fill}%
    \begin{minipage}[t]{0.31\textwidth}
        \centering
        \includegraphics[trim={0cm 0cm 0cm 0cm}, clip,width=\linewidth]
        {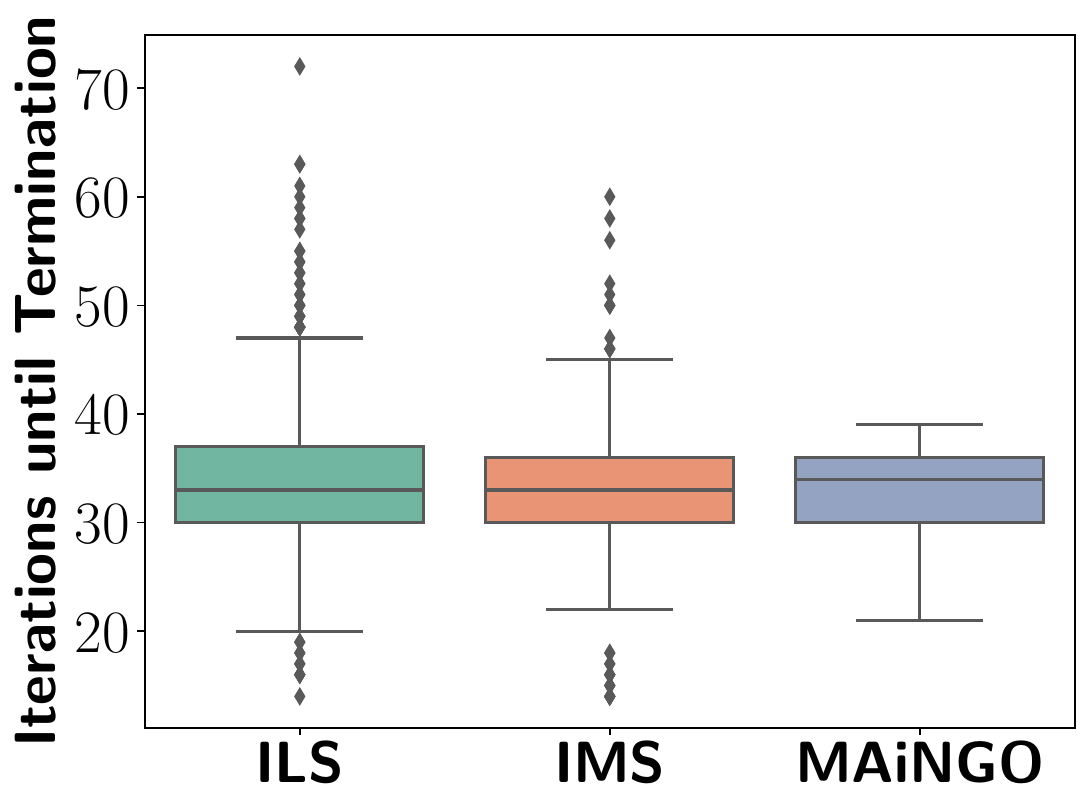}
        \subcaption{
            $\kappa = 3$ and $N = 3$.
        }
    \end{minipage}%
    \hspace{\fill}%
    \begin{minipage}[t]{0.31\textwidth}
        \centering
        \includegraphics[trim={0cm 0cm 0cm 0cm}, clip,width=\linewidth]
        {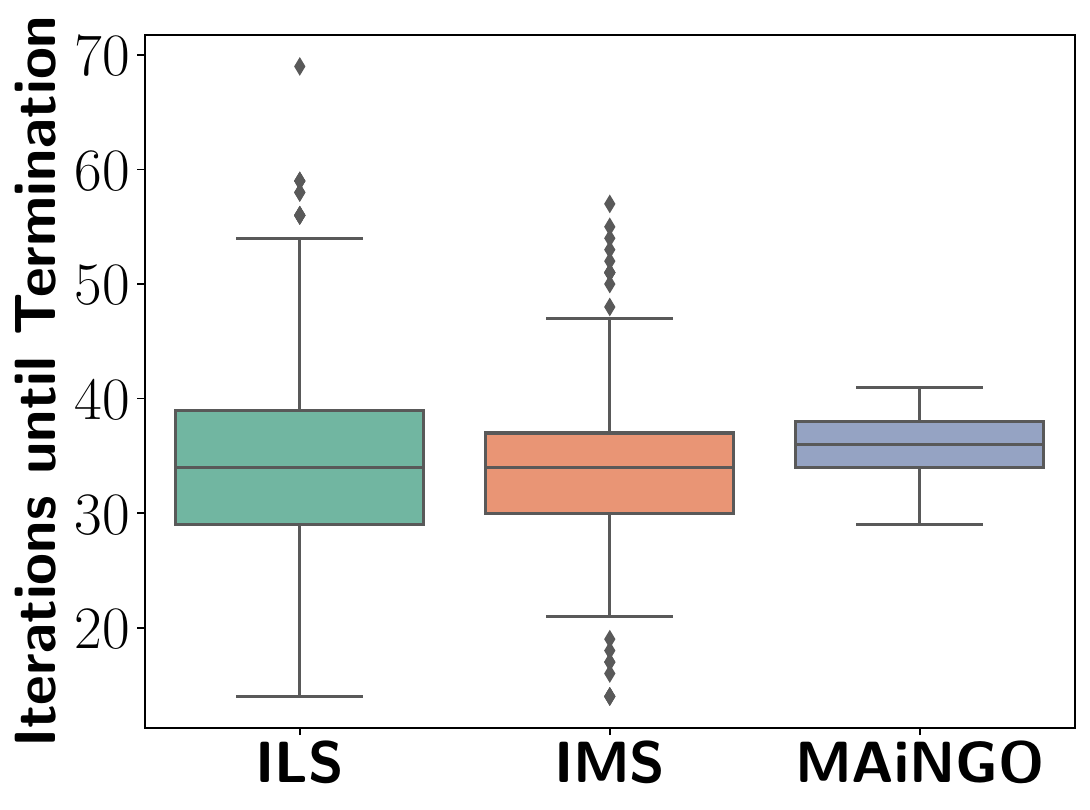}
        \subcaption{
            $\kappa =$ ~\ref{eq:kappa_t_Srinivas_2012} and $N = 3$.
        }
    \end{minipage}%
    \caption{
        Number of iterations until termination for all runs, despite which minimum they converged to, for a series of M\"{u}ller-Brown case studies.
        Details on the number of experiments, runs, and termination criteria are provided in \cref{sec:set-up-details}.
    }
\end{figure}

\begin{figure}[htbp]
    \begin{minipage}[t]{0.31\textwidth}
        \centering
        \includegraphics[trim={0cm 0cm 0cm 0cm}, clip,width=\textwidth]
        {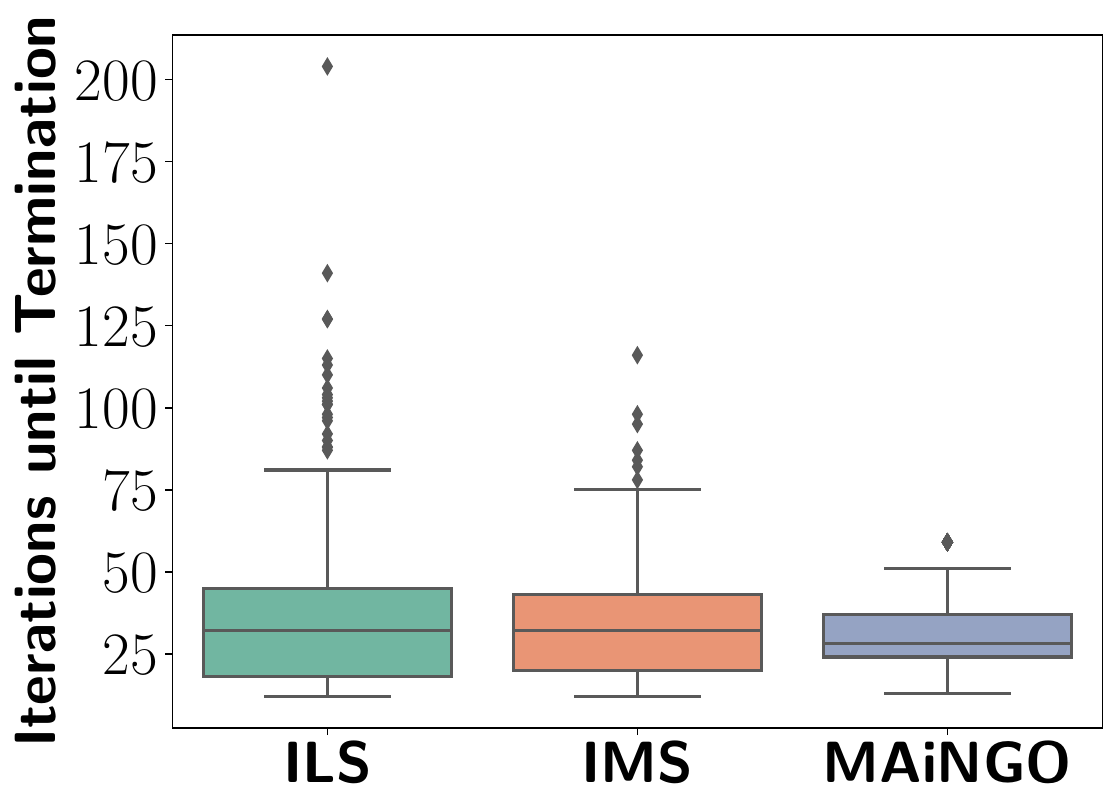}
        \subcaption{
            2D GKLS.
        }
    \end{minipage}%
    \hspace{\fill}%
    \begin{minipage}[t]{0.31\textwidth}
        \centering
        \includegraphics[trim={0cm 0cm 0cm 0cm}, clip,width=\linewidth]
        {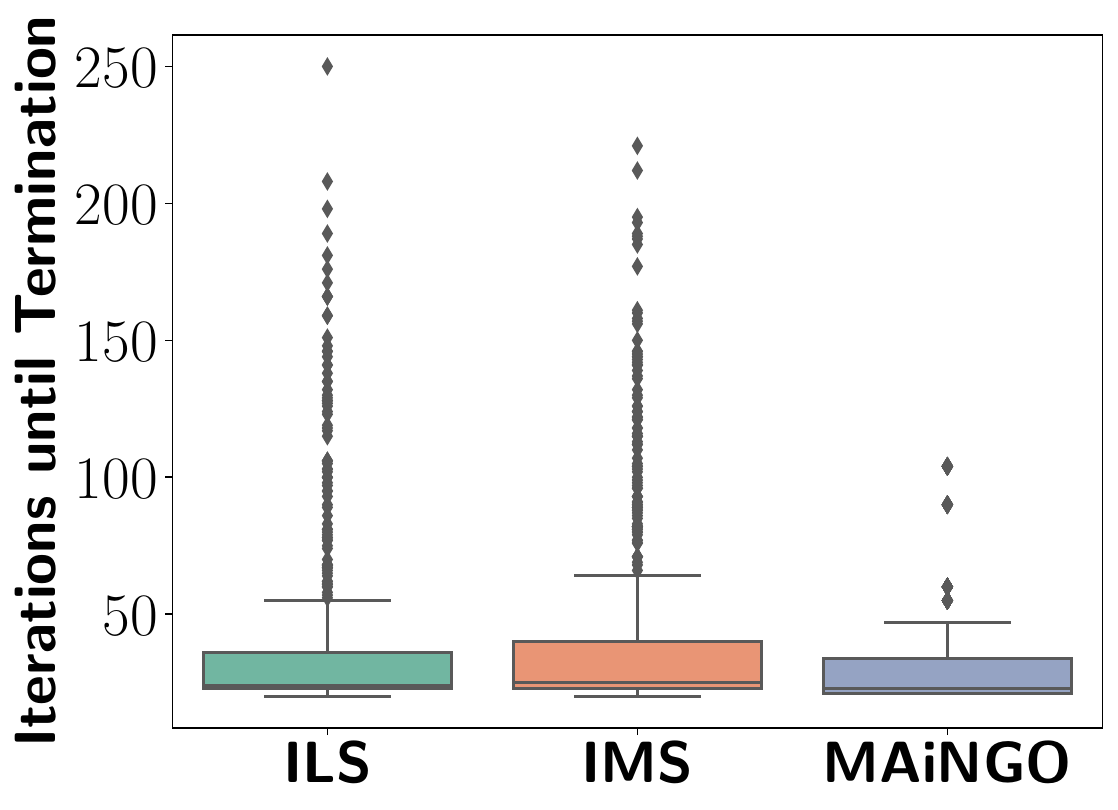}
        \subcaption{
            3D GKLS.
        }
    \end{minipage}%
    \hspace{\fill}%
    \begin{minipage}[t]{0.31\textwidth}
        \centering
        \includegraphics[trim={0cm 0cm 0cm 0cm}, clip,width=\linewidth]
        {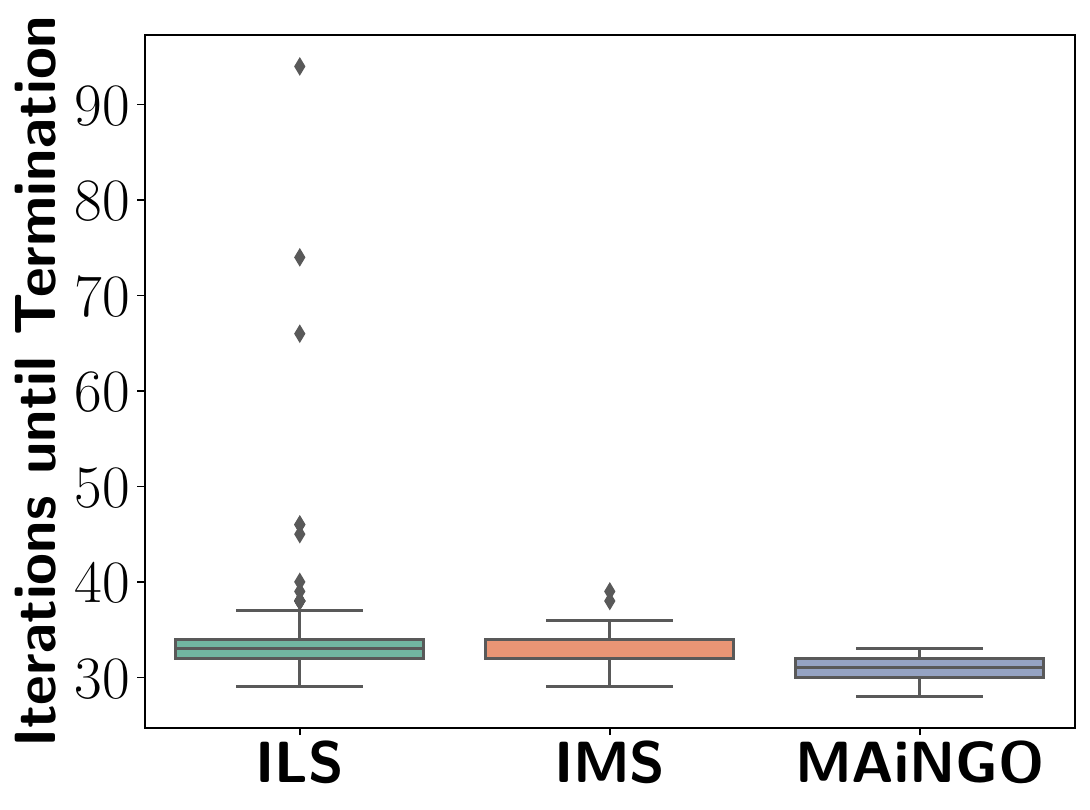}
        \subcaption{
            4D GKLS.
        }
    \end{minipage}
    \caption{
        Number of iterations until termination for all runs, despite which minimum they converged to, for a series of GKLS case studies.
        Details on the number of experiments, runs, and termination criteria are provided in \cref{sec:set-up-details}.
    }
\end{figure}

\clearpage

\section{Conditional Logistic Regression for Comparing Solver Probabilities of Convergence Averaged Over Datasets}
\label{sec:stats_proof}
We show the following.
\vspace*{.3cm}

\noindent {\bf Result 1.} Suppose the three solvers have the same convergence probability when averaged over datasets,
\begin{equation}
     \label{eq:constancya}
     p^{\text{ILS}} = p^{\text{MAiNGO}} = p^{\text{IMS}},
\end{equation}
but the true odds ratios of convergence between solvers \textit{differ across datasets}; that is, the model of conditional logistic regression \eqref{eq:ormodel} is incorrect.
Then the cMLE of $\alpha, \beta$ as a test statistic for \eqref{eq:constancya} -- centers at zero.
\bigskip

The result allows the cMLE to test for equality of overall probabilities between solvers.
This is particularly important because odds ratios generally are not collapsible, i.e., the log ratio of odds of two averaged probabilities over a matching factor, generally differs from the averaged log ratio of odds of the probabilities conditionally on the matching factor.
\bigskip

For ease, we will focus only in comparing solvers $\bar{s}$ (MAiNGO) (one run per dataset) to $\hat{s}$ (ILS) (31 runs per dataset).
Below, we denote by $pr^{mod}(\cdot )$ and $pr( \cdot)$ the model-based and true probabilities of events $(\cdot)$; in particular we refer to $p^{\bar{s}}_{dat}, p^{\hat{s}}_{dat}$ as the true average probabilities of convergence for dataset $dat$, and $p^{mod,\bar{s}}_{dat}, p^{mod, \hat{s}}_{dat}$ as the corresponding model-based probabilities.

\textit{The logistic regression likelihood.}
This likelihood is formed as the product of dataset-specific likelihoods.
Each of those is the likelihood of the data $y^{\bar{s}}_{dat}, y^{\hat{s}}_{dat, (+)}$ conditionally on the total number $T_{dat}$ of converged runs over both solvers, and conditionally on the model-based probabilities given $dat$.
Because, given $T_{dat}$ only $y^{\bar{s}}_{dat}$ is free (as $y^{\hat{s}}_{dat, (+)}$ is determined as $T_{dat}-y^{\bar{s}}_{dat}$), the likelihood is well known to be
\begin{align}
    pr^{mod}(y^{\bar{s}}_{dat} \mid T_{dat}, p^{mod,\bar{s}}_{dat}, p^{mod, \hat{s}}_{dat})& = q(\alpha,T_{dat})^{y^{\bar{s}}_{dat}}\cdot \{1-q(\alpha,T_{dat})\}^{1-y^{\bar{s}}_{dat}} \label{eq:clkdi}\\
    \mbox{ where }q(\alpha,T_{dat})&:=pr^{mod}(y^{\bar{s}}_{dat}=1 \mid T_{dat}, p^{mod,\bar{s}}_{dat}, p^{mod, \hat{s}}_{dat})\nonumber\\[-.4cm] \label{eq:qmodel}\\[-.1cm]&=\frac{\exp{(\alpha)}T_{dat}}{\exp{(\alpha)}T_{dat}+(31+1-T_{dat})}\nonumber
\end{align}

\textit{The conditional Maximum Likelihood Estimator }$\hat{\alpha}^{cMLE}$ solves the score equations $\sum_{dat}S_{dat}(\alpha)$ =0, where $S_{dat}(\alpha)$ is the derivative of the log of \eqref{eq:clkdi}.
When $T_{dat}$ is 0 or $31+1$, $q$ does not depend on $\alpha$ so $S_{dat}(\alpha)$ is 0.
For the other values of $T_i$, and using the fact that $\{\partial q(\alpha,T_{dat})/\partial \alpha \}= \{ q(\alpha,T_{dat})(1-q(\alpha,T_{dat}) )\}$, it is easy to see that

\begin{equation}
    S_{dat}(\alpha)= \begin{cases}  y^{\bar{s}}_{dat} -  q(\alpha,T_{dat})& \mbox{ if } T_{dat} \not\in \{0, (31+1) \} \\
    0 & \mbox{ otherwise } \end{cases} \nonumber
\end{equation}

\textit{Validity of the cMLE for evaluating \eqref{eq:constancya} even when the model is incorrect.} We want to show that if the average solver probabilities are the same, the cMLE of $\alpha$ will be centered at 0, whether the model is correct or not.

We know that, whether or not the model is correct, the cMLE of $\alpha$ is centered at the value, say $\alpha^*$, that solves the average score equations, i.e.,
\begin{equation}
    \label{eq:escore1}
    \mathbb{E}\{  y^{\bar{s}}_{dat} -  q(\alpha^*,T_{dat}) \mid T_{dat} \not\in \{0, (31+1) \}  \} = 0
\end{equation}

If the model is correct, and we have \eqref{eq:constancya}, then the value $\alpha^*=0$ is the correct value of the model, and also solves the above equation, which implies that the cMLE has expected value $\alpha^*=0$.
This can be seen by taking, inside \eqref{eq:escore1}, the expected value, first conditionally on $T_{dat}$ and the model-based $(p^{mod,\bar{s}}_{dat}, p^{mod, \hat{s}}_{dat})$.
There, the model implies the expectation of $y^{\bar{s}}_{dat}$ is precisely $q(\alpha^*,T_{dat})$ (c.f.,  \eqref{eq:qmodel}) and the difference is 0 for every $dat$.

If, however, the model is not correct, then there is no true value of $\alpha$, and this case is not well examined.
In this case, of course, the cMLE of $\alpha$ still exists and is still centered at the value, say $\alpha^{**}$,  that solves \eqref{eq:escore1}, although now the expectation $\mathbb{E}$ is taken with respect to \textit{whatever the unknown correct distribution is}.
In general, there is no closed form expression for $\alpha^{**}$ in terms of the true distributions.
However, what we ask is simpler: to know if $\alpha^{**}=0$ is a solution to the equation \eqref{eq:escore1} when the model \eqref{eq:qmodel} is incorrect but the hypothesis that marginally $p^{\bar{s}}=p^{\hat{s}}$ is true \eqref{eq:constancya}.
We can reverse the question by first simplifying: if $\alpha^{**}$ is first set to zero, is the score equation (\ref{eq:escore1}) solved by the hypothesis $p^{\bar{s}}=p^{\hat{s}}$  ?  At $\alpha^{**}=0$ the LHS of (\ref{eq:escore1}) becomes:
\begin{align}
    \mathbb{E}\{  y^{\bar{s}}_{dat} -  q(0,T_{dat}) \mid T_{dat} \not\in \{0, (31+1) \}  \} =  \mathbb{E}\{  y^{\bar{s}}_{dat} -  \frac{T_{dat}}{31+1} \mid T_{dat} \not\in \{0, (31+1) \}  \} \nonumber
\end{align}
Moreover, because $T_{dat} =0$ or $(31+1)$ implies $y^{\bar{s}}_{dat} -  \frac{T_{dat}}{31+1}=0$ anyway, a zero for the RHS of the last expression is equivalent to
\begin{equation}
    \begin{aligned}
        0\quad=\quad & \mathbb{E}\left\{y^{\bar{s}}_{dat} - \frac{T_{dat}}{31+1} \right\} \mbox{ (unconditionally) } \\
              =\quad & \mathbb{E}\left\{y^{\bar{s}}_{dat} \right\} - \mathbb{E}\left( \frac{ y^{\bar{s}}_{dat} + \sum_{r=1}^{31}y^{\bar{s}}_{dat,r}}{1+31} \right) \\
              =\quad & p^{\bar{s}} - \frac{1}{1+31} \left( p^{\bar{s}} +31 p^{\hat{s}} \right)
    \end{aligned}
\end{equation}
Finally, under the assumption $p^{\bar{s}}=p^{\bar{s}}$, the last expression is indeed zero, even if the model assumptions are false, which shows that \eqref{eq:escore1} is true and so establishes the result.

\end{document}